\documentclass[a4paper,10pt,twoside,onecolumn,final]{article}
\usepackage{amsmath,amsfonts,amssymb,amscd,mathrsfs,xfrac}
\usepackage{extarrows}
\usepackage{epsfig}
\usepackage{enumerate}
\usepackage{color}
\usepackage{graphicx}
\usepackage{psfrag}
\usepackage[nodayofweek]{datetime}
\usepackage[paper=a4paper,hscale=0.7,vscale=0.75,centering,heightrounded]{geometry}
\usepackage[ntheorem]{empheq} % this loads amsmath as well 
\empheqset{box=\bigfbox} 
\usepackage[thmmarks,amsmath,hyperref]{ntheorem} 
\theoremstyle{plain} 
\theoremheaderfont{\normalfont\bfseries} 
\theorembodyfont{\slshape}
\theoremsymbol{\ensuremath{_\Box}}
\theoremnumbering{Alph}
\newtheorem{theorem}{Theorem}[section]

\theoremnumbering{arabic}
\newtheorem{proposition}{Proposition}[section]

\newtheorem{corollary}[proposition]{Corollary}

\theoremstyle{plain}
\theoremheaderfont{\normalfont\bfseries} 
\theorembodyfont{\upshape}
\theoremsymbol{\ensuremath{_\Box}}
\newtheorem{definition}[proposition]{Definition}  % [chapter]
\newtheorem{remark}[proposition]{Remark}  % [chapter]
\newtheorem{example}[proposition]{Example}  % [chapter]

\theoremstyle{nonumberplain}
\theoremheaderfont{\normalfont\bfseries} 
\theorembodyfont{\upshape}
\theoremsymbol{\ensuremath{_\blacksquare}} % _\blacksquare % _Box
\newtheorem{proof}{Proof}
\theoremsymbol{\ensuremath{_\blacksquare}}

\theoremlisttype{allname}
\allowdisplaybreaks[3] 
\definecolor{darkgray}{gray}{0.4}  
\definecolor{ddarkgray}{gray}{0.2} 
\definecolor{redgray}{rgb}{0.5,0.25,0.25}  
\definecolor{bluegray}{rgb}{0.25,0.25,0.5}  
\usepackage[colorlinks=true, pdfstartview=FitV, linkcolor=redgray, citecolor=bluegray, urlcolor=black, pdfpagelabels, naturalnames=true, draft=false]{hyperref}
\usepackage[small,bf]{caption} 

  \renewenvironment{thebibliography}[1]{%
    \begin{oldthebibliography}{#1}%
       \addtolength{\itemsep}{-0.5ex}%
  }%
  {%
    \end{oldthebibliography}%
  }
\usepackage{marginnote}

\newcommand{\ReLU}{\operatorname{ReLU}}
\title{{\Large \textsc{Neural Control of Discrete Weak Formulations:
Galerkin, %, Petrov--Galerkin,
Least-Squares \& Minimal-Residual Methods with Quasi-Optimal
Weights}}}
\author{Ignacio Brevis\footnotemark[2]\,, Ignacio Muga\footnotemark[2]\,, Kristoffer G.~van der Zee\footnotemark[1]}
\date{{\small 18 May 2022}}
\newlength{\bigfboxsep}
\newcommand{\bigfbox}[1]{\setlength{\fboxsep}{\bigfboxsep}\fbox{#1}}

\newcommand{\negquad}{\mspace{-18.0mu}}
\newcommand{\negqquad}{\negquad\negquad}

\newcommand{\norm}[1]{{\|#1\|}}
\newcommand{\bignorm}[1]{\big\|#1\big\|}

\newcommand{\innerprod}[1]{(#1)}
\newcommand{\biginnerprod}[1]{\big(#1\big)}

\newcommand{\bigdual}[1]{\big\langle#1\big\rangle}
\newcommand{\Bigdual}[1]{\Big\langle#1\Big\rangle}

\providecommand{\grad}{\nabla}
\renewcommand{\grad}{\nabla}

\newcommand{\dd}{\mathrm{d}}

\begin{document}
%-----------------------------------------------------------------------------%
%
%=============================================================================%
% Title
%=============================================================================%
%
%-----------------------------------------------------------------------------%
\maketitle
%-----------------------------------------------------------------------------%
%
%=============================================================================%
% Author footnotes
%=============================================================================%
%
%-----------------------------------------------------------------------------%
\renewcommand{\thefootnote}{\fnsymbol{footnote}}
\footnotetext[2]{Instituto de Matemáticas, Pontificia Universidad Católica de Valparaíso, Chile}
\footnotetext[1]{School of Mathematical Sciences, University of Nottingham, UK; Corresponding author, kg.vanderzee@nottingham.ac.uk\\ \hspace*{1.5em} {\scriptsize Dedicated to J. Tinsley Oden.}}
\renewcommand{\thefootnote}{\arabic{footnote}}
%-----------------------------------------------------------------------------%
%
%=============================================================================%
% Abstract
%=============================================================================%
%
%-----------------------------------------------------------------------------%
\begin{abstract}
There is tremendous potential in using neural networks to optimize numerical methods. 
%In recent years there has been tremendous interest in using neural networks to optimize numerical methods
In this paper, we introduce and analyse a framework for the \emph{neural optimization of discrete weak formulations}, suitable for finite element methods.
%
%We introduce and analyse % \emph{Neural Weak Formulations}
%\emph{finite element weak formulations with neural control}.
%the \emph{neural optimization of discrete weak formulations}.
The main idea of the framework is to include a neural-network function acting as a \emph{control} variable in the weak form. 
% hat the weak form contains a neural-network function that acts as a control variable.
%that the weak form contains a 
Finding the neural control that (quasi-) minimizes a suitable cost (or loss) functional, then yields a numerical approximation with desirable attributes. In particular, the framework allows in a natural way the incorporation of known data of the exact solution, or the incorporation of stabilization mechanisms (e.g., to remove spurious oscillations). 
%By finding the minimizing a suitable cost (or loss) functional, the numerical approximation  
% Finding the neural-network function that quasi-minimizes a suitable cost (or loss) functional, the numerical approximation  
%
%. there is a control parameter  underlying parameter is a \emph{neural network} (or other suitably-expressive approximant).
% The purpose of the neural control is to allow for the training of an \emph{optimal} discretization method. Optimality is defined in the sense of minimizing a suitable functional containing, e.g., (partial) data on the exact solution, or desired features of the resulting approximations . 
\par
The main result of our analysis pertains to the well-posedness and convergence of the associated constrained-optimization problem. In particular, we prove under certain conditions, that the discrete weak forms are stable, and that  quasi-minimizing neural controls exist, which converge quasi-optimally.
%, and for which the resulting discrete weak form is stable,  the discrete problem is well-posed and t
%
%discretely-minimized Neural Weak Formulations converge quasi-optimally to the continuous minimizer.
% \par
We specialize the analysis results 
% to conforming finite element methods based on  in particular, 
to Galerkin, % Petrov--Galerkin,
least-squares and minimal-residual formulations, where the neural-network dependence appears in the form of suitable weights.  
% functions.
Elementary numerical experiments support our findings and demonstrate the potential of the framework.
\end{abstract}
%-----------------------------------------------------------------------------%
%
%=============================================================================%
% Tableofcontents
%=============================================================================%
%
%-----------------------------------------------------------------------------%
\clearpage
\tableofcontents
%-----------------------------------------------------------------------------%
%
%=============================================================================%
% Headings
%=============================================================================%
%
%-----------------------------------------------------------------------------%
\pagestyle{myheadings}
\thispagestyle{plain}
\markboth{\small Neural control of discrete weak formulations}{\small Neural control of discrete weak formulations}
%-----------------------------------------------------------------------------%
%
%=============================================================================%
\clearpage
\section{Introduction}
In recent years there has been tremendous interest in the merging of neural networks and machine-learning algorithms with traditional methods in scientific computing and computational science~\cite{HigHigSIREV2019, ECCP2020, KarKevLuPerWanYanNRP2021, PenAlbBugCanDeDurGarKarLytPerPetKuhACME2021}. In this paper we demonstrate how neural networks can be utilized to optimize finite element methods.
\par
%In one of its most familiar forms,
In one of its most familiar mathematical forms, the finite element method is a discretization technique for partial differential equations (PDEs) based on a weak formulation using discrete subspaces, i.e., the exact solution $u\in \mathbb{U}$ is approximated by~$u_h\in \mathbb{U}_h$, which is the unique solution of the discrete problem:
\begin{alignat}{2}
\notag
 &\text{Find } u_h \in \mathbb{U}_h:
 \\
\label{eq:introPG}
 &\qquad  b(u_h,v_h) = f(v_h)\;,
 \qquad \forall v_h\in \mathbb{V}_h\,,
\end{alignat}
where~$\mathbb{U}_h$ is a discrete subspace of the infinite-dimensional Hilbert or Banach space~$\mathbb{U}$ (typically a Sobolev space on a domain~$\Omega\subset\mathbb{R}^d$), $\mathbb{V}_h$ is a subspace of a Hilbert or Banach space~$\mathbb{V}$ with $\dim \mathbb{V}_h = \dim \mathbb{U}_h$, $b:\mathbb{U}\times \mathbb{V}\rightarrow \mathbb{R}$ is a continuous bilinear form, $f:\mathbb{V}\rightarrow \mathbb{R}$ a continuous linear form, and the exact solution~$u$ satisfies $b(u,v) = f(v)$ for all~$v\in \mathbb{V}$.%
\footnote{When $\mathbb{U}_h = \mathbb{V}_h$, this is a \emph{Galerkin} method, otherwise it is a \emph{Petrov--Galerkin} method.}
\par
It is well-known that the accuracy of~$u_h$ can be improved by enlarging~$\mathbb{U}_h$ (e.g., by refining the underlying finite element mesh).%
\footnote{Indeed, a~priori error analysis reveals that $\norm{u-u_h}_{\mathbb{U}} \le C \inf_{w_h\in \mathbb{U}_h}\norm{u-w_h}_{\mathbb{U}}$, provided $b(\cdot,\cdot)$ satisfies a discrete inf--sup condition on $\mathbb{U}_h\times \mathbb{V}_h$;
%and where~$u\in \mathbb{U}$ solves $b(u,v) = f(v)$ for all~$v\in \mathbb{V}$; 
see e.g.,~\cite{OdeRedBOOK2011, ErnGueBOOK2021b}.
}
However, \emph{for a fixed value of~$h$}, the particular~$u_h$ defined by~\eqref{eq:introPG} may be very \emph{unsatisfactory}.
% as measured by certain output quantities of interest. 
In fact, there is no reason why a certain quantity of interest of~$u_h$ is accurate at all,%
\footnote{E.g., the value~$u_h(x_0)$ for some point~$x_0\in\Omega$ is generally quite distinct from~$u(x_0)$.}
or why the approximation inherits certain qualitative features of the exact solution.%
\footnote{E.g., $u_h$ may exhibit spurious oscillations, while~$u$ is monotone.}
%
% , because those quantities are not appearing in th.
% For example, the value~$u_h(x_0)$ of the approximation at some point~$x_0$ in the domain~$\Omega$ is generally inaccurate, unless $x_0$ happens to  . 
Indeed, the discrete problem~\eqref{eq:introPG} is a \emph{rigid} statement in the sense that it identifies a \emph{single} element in~$\mathbb{U}_h$, irrespective of desired attributes, whereas there could be many other elements in~$\mathbb{U}_h$ that are far superior. 
\subsection{Neural optimization of discrete weak forms}
The objective of this work is to propose and analyse a framework 
% that allows for neural control of
for the \emph{neural optimization of discrete weak formulations} to significantly improve quantitative and qualitative attributes of discrete approximations. In particular, we consider \emph{Galerkin}, %\emph{Petrov--Galerkin},
\emph{least-squares}, and \emph{minimal-residual} formulations.
\par
% In essence, 
The main idea of the framework is that it incorporates a neural-network function~$\xi$ as a control variable in the discrete test space~$\mathbb{V}_h(\xi)$. That is, the approximation~$u_h = u_{h,\xi}$ now depends on~$\xi$ and solves the discrete problem:%
\begin{alignat}{2}
\notag
 &\text{Find } u_h=u_{h,\xi} \in \mathbb{U}_h:
 \\
\label{eq:introPGxi}
 &\qquad  b(u_{h,\xi},v_h) = f(v_h)\;,
 \qquad \forall v_h\in \mathbb{V}_h(\xi)\,.
\end{alignat}
Then, in order to obtain a desired approximation~$u_{h,\bar{\xi}}$, we aim to find a neural-network function~$\bar{\xi}$ that \emph{quasi}-minimizes a desired cost (or loss) functional:%
\footnote{We also allow for the inclusion of a regularization term in the cost functional; see Section~\ref{ssec:Problem}.
}
\begin{alignat}{2}
\label{intro:qminJ}
 J(u_{h,\bar{\xi}}) \longrightarrow \operatorname{quasi-min}\,.
\end{alignat}
%
% \par
%
The notion of quasi-minimization is critical when aiming to minimize over a set of neural-network functions (i.e., the set of functions implemented by neural networks of a fixed architecture); see Section~\ref{ssec:quasimin} for further details (in particular, Definitions~\ref{def:quasimin} and \ref{def:qminprob}).
\par
The quasi-minimization problem~\eqref{intro:qminJ} is essentially a \emph{nonstandard} PDE-constrained optimization, with the nonstandard part being the dependence of the state problem~\eqref{eq:introPGxi} on~$\xi$ via the discrete test space~$\mathbb{V}_h(\xi)$. Importantly, $\mathbb{V}_h(\xi)$~will be parameterized by~$\xi$ in such a way so as to ensure \emph{stability} of the discrete problem~\eqref{eq:introPGxi}. Moreover, as will become clear in the following sections, the basis functions in $\mathbb{V}_h(\xi)$ need not be computed explicitly, but equivalent formulations to~\eqref{eq:introPGxi} can be used,
% devoid of $\xi$-dependent spaces, 
which instead incorporate~$\xi$ by means of suitable \emph{weight functions}. These formulations essentially lead to a PDE-constrained optimization with a nonlinear control-to-state map.
\subsection{Potential of the methodology}
There are two main benefits of having neural control of discrete weak forms:
\begin{itemize}
\item \emph{Incorporation of data}: Knowledge of quantities of the exact solution can be taken into account in a natural way by setting, for example, 
\begin{alignat*}{2}
  J(u_{h,\xi}) =
 \frac{1}{2}\big| q(u_{h,\xi}) - \bar{q} \big|^2\,,
\end{alignat*}
where $q:\mathbb{U}\rightarrow \mathbb{R}$ is a functional measuring the quantity of interest and $\bar{q}\in \mathbb{R}$ is known data.%
\footnote{The data $\bar{q}$ represents~$q(u)$, and it could be obtained through experiments, high-fidelity computation, or otherwise.}
Minimizing such a~$J(\cdot)$ ensures that the discrete solution~$u_h$ to~\eqref{eq:introPGxi} is \emph{data-driven} in the sense that $u_h$ becomes constrained by the data.%
\footnote{This is somewhat similar in spirit to \emph{physics-informed neural networks} (PINN)~\cite{RaiPerKarJCP2019}, where however a single neural-network function minimizes a combination of the residual and data misfit.}
We note that multiple quantities can be taken into account using, for example,
\begin{alignat*}{2}
  J(u_{h,\xi}) = \frac{1}{N_{\mathrm{data}}}\sum_{i=1}^{N_{\mathrm{data}}}
 \frac{1}{2}\Big| q_i(u_{h,\xi}) - \bar{q}_i \Big|^2\,,
\end{alignat*}
or, more generally, using some operator~$Q:\mathbb{U}\rightarrow \mathbb{Z}$; see Section~\ref{sec:Abstract}.
%where $q_i$ and $\bar{q}_i$, for $i=1,\ldots,N$ are various  and corresponding known data.
%similar to PINN
%
\item \emph{Incorporation of stabilization mechanisms}: Qualitative attributes of the discrete solution can be enhanced by minimizing a suitably-chosen~$J(\cdot)$. In this way discrete solutions can be enforced to, e.g., satisfy an a~priori known maximum principle, have monotone (or spurious oscillation free) behavior around discontinuities and layers, or have a certain discrete wave number (i.e., free from pollution). In the past decades, many different stabilized finite element methods have been proposed (and analyzed) that impose such attributes~\cite{GueSINUM2004, BurErnMOC2005, JohKnoCMAME2007, EvaHugSanCMAME2009, DemGopMugZitCMAME2012, PetMOC2017}. Within our framework such a method is naturally obtained after (quasi-) minimization (i.e., method~\eqref{eq:introPGxi} with~$\xi = \bar{\xi}$). As an example, Guermond~\cite{GueSINUM2004} advocates the $L^1$-minimization of the residual; in other words, within our framework one would choose:
\begin{alignat*}{2}
  J(u_{h,\xi}) = \bignorm{f - B u_{h,\xi} }_{L^1(\Omega)}\,,
\end{alignat*}
where $f-B u_{h,\xi}$ is the strong form of the residual.
\end{itemize}
\par
The idea of using neural networks to parameterize the test space was initially proposed in our earlier work~\cite{BreMugZeeCAMWA2021}, where it was restricted to minimal-residual formulations within a parametric PDE setting. The current work presents significantly more general settings and formulations as well as analyses of their well-posedness and convergence.
\par
While the above shows examples of $J(\cdot)$ corresponding to \emph{unsupervised} learning (i.e., there is no need to know the exact solution~$u$), when the original problem is \emph{parametric} itself (e.g., a parametric PDE), \emph{supervised} learning becomes meaningful. Indeed, in that case, the data may be the exact solution $u_{\lambda_i}$ for certain parameters~$\lambda_i$, $i=1,\ldots,N_{\mathrm{data}}$. This then allows for the training of finite element discretizations with superior accuracy in quantities of interest even on very coarse meshes. We refer to our earlier work~\cite{BreMugZeeCAMWA2021} for the methodology and illustrative examples in that case.
%
%
% We notice that where the definition of quasi-minimization is explained in Section~.
% and solve the quasi-minimization problem:
% %
% \begin{alignat}{2}
% \notag
%  &\text{Find } \xi_n\in \mathcal{M}_n:
%  \\
% \label{eq:introMin}
%  &\qquad  J(u_{h,\xi_n},\xi_n) \le  \inf_{\eta_n \in \mathcal{M}_n} J(u_{h,\eta_n},\eta_n) + \delta_n\;,
% \end{alignat}
%
% a quasi-minimizer of a functional in order to find a quasi-optimal  cost function (or loss function)~$J(u_{h},\xi)$
% Within a set of neural network functions associated to 
%
\par
\subsection{Main contributions: Well-posedness, convergent quasi-minimizers, weighted conforming formulations}
Let us briefly outline the main contributions of this work. 
The first main contribution is the analysis of an abstract constrained-optimization problem associated to~\eqref{intro:qminJ}; see Section~\ref{sec:Abstract}. In particular, we consider an abstract state problem equivalent to~\eqref{eq:introPGxi}, but in the form of a \emph{mixed} system with a~$\xi$-dependent bilinear form.%
\footnote{The mixed system is motivated by residual-minimization theory~\cite{DemGopBOOK-CH2017, MugZeeSINUM2020}: Minimal residual formulations are equivalent to mixed systems, which in turn are equivalent to Petrov--Galerkin formulations.}
%
% This system depends on~$\xi$ via a bilinear form, and implicitly defines a test space~$\mathbb{V}_h(\xi)$ for which equivalence holds with~\eqref{}
%  where the $\xi$-dependence appears within a bilinear form. 
% We demonstrate that it is equivalent to~\eqref{eq:introPGxi}, and 
We prove, under suitable conditions, that the state problem is \emph{well-posed} (uniformly with respect to~$\xi$); see Proposition~\ref{prop:stateEq}. Furthermore, we present differentiability conditions (on the $\xi$-dependence) that allow us to prove the \emph{existence} of quasi-minimizers (within sets of neural-network functions, of some size~$n$) to the associated constrained optimization~\eqref{intro:qminJ}, which converge \emph{quasi-optimally} (upon $n\rightarrow \infty$); see Corollary~\ref{cor:constrainedProb} for details.
% (also for a precise meaning of quasi-optimal convergence for quasi-minimizers). 
%
\par
We note that our analysis is based on a fundamental result for the quasi-minimization of strongly-convex and differentiable functionals (see Theorem~\ref{thm:quasi}), which is of independent interest and applies, e.g., to the analysis of deep Ritz methods~\cite{XuCCP2020, PouJCMDS2022, MulZei2021} and PINN methods~\cite{ShiZhaKar2020, MisMolIMANA2022b, CaiCheLiuJCP2021}.
\par
The second main contribution of this work is the application of our framework to certain weak formulations used by conforming finite element methods; see Section~\ref{sec:formulations}. In these applications, the neural-network control variable~$\xi$ will appear by means of suitable weights in the bilinear forms. In particular, we will analyse weighted least-squares, weighted Galerkin,
%  (and weighted Petrov--Galerkin), 
and weighted minimal-residual formulations. 
\par
For weighted least-squares and weighted minimal-residual formulations, suitable conditions on the weights imply (via the abstract result of the first main contribution) stability of the discrete problem (uniformly in $\xi$). Furthermore, suitable differentiability conditions on the weights imply existence of (quasi-optimally) convergent quasi-minimizers of the associated constrained minimization. 
\par
On the other hand, for weighted Galerkin,
%(and weighted Petrov--Galerkin), 
it turns out that stability is \emph{not} immediate, and may require constraints on~$\xi$ depending on the problem at hand.%
\footnote{In essence, the reason for instability relates to a discrete inf-sup condition of a weighted bilinear form.}
Therefore, neural control is far more convenient
for least-squares and minimal-residual formulations, the fundamental reason being the inherent stability that comes with their underlying minimization principle.
% It is the inherent stability, which comes with the minimization principle underlying least-squares and minimal-residual formulations, that make neural control ideal for these that and we advocate the use of minimal residual principles for 

% In particular, we show that   (see Theorem~\ref{thm:reduced})
%
\par
We support our findings with numerical experiments in Section~\ref{sec:Numerics}. While our theoretical results directly apply to any linear operator, we choose the advection-reaction PDE to illustrate various numerical aspects, viz., the incorporation of data (Section~\ref{ssec:numExp:QoI}), the quasi-optimal convergence of quasi-minimizers (Section~\ref{ssec:numExp:convANN}), and the incorporation of $L^1$-type stabilization (Section~\ref{ssec:numExp:L1}). 
%
% First, Abstract framework, prove discrete existence of quasi-minimisers, quasi-optimal convergence.
% Second, Formulations.
% We also illustrate the capabilities with some elementary numerical examples. We also refer to~\cite{BreMugZeeCAMWA2021} for numerical examples for MinRes in a parametric setting and upcoming work.
%
% The purpose of this framework is to ssignificantly improve the accuracy of discrete approximations. This is achieved by means of a control variable~$\xi$ by allowing .  
%
% The advantages of neural 
%
% (hence fixed~$\mathbb{U}_h$ and~$\mathbb{V}_h$), the setting  for a fixed value of~ subspaces~$\mathbb{U}$
%
%  (upon which $\mathbb{U}_h$ and~$\mathbb{V}_h$ become dense subspaces of~$\mathbb{U}$ and~$\mathbb{V}$, respectively)
%
% n approximation to some operator equation
%
% \begin{alignat*}{2}
%  B u = f
% \end{alignat*}
% %
% where~$B:\mathbb{V} \rightarrow \mathbb{V}^*$ continuous linear, $f\in \mathbb{V}^*$, in the form
% %
% \begin{alignat*}{2}
%  u_h = B_h^{-1} f_h
% \end{alignat*}
% %
% with $B_h:\mathbb{V}_h \rightarrow \mathbb{V_h}^*$
% %
% \par
% %
% Neural weak formulations:
% %
% \begin{alignat*}{2}
%  B_{h,\omega} u_{h,\omega} = f_{h,\omega}
% \end{alignat*}
% %
% where $\omega$ has been learned from data.

% The aim of this paper is to propose a general framework for the neural control of weak formulations
%

\subsection{Related work}
There are a number of works related to ours.
\par
%\subsubsection{Optimizing methods}
\emph{Optimizing numerical methods}: Traditionally, the incorporation of known data or other desired attributes in numerical PDE approximations is achieved via the method of Lagrange multipliers, see e.g., Evans, Hughes \&~Sangalli~\cite{EvaHugSanCMAME2009}, Kergrene, Prudhomme, Chamoin \&~Laforest~\cite{KerPruChaLafCMAME2017}, and references therein. 
More recently, neural networks have been proposed to learn the parameters that define a numerical method; see Ray \&~Hesthaven~\cite{RayHesJCP2018}, Mishra~\cite{MisMINE2018} and others~\cite{BarHoyHicBrePNAS2019, DisHesRayJCP2020, WanSheLonDonCICP2020, SchRayHesJCP2021}. Interestingly, a recent learning methodology for adaptive mesh refinement has been proposed that ensures optimal convergence; see Bohn \&~Feischl~\cite{BohFeiCAMWA2021}.
Within the context of optimizing finite-element formulations, a minimal-residual framework that ensures stability was proposed in our previous work~\cite{BreMugZeeCAMWA2021}. 
Our current work contributes to these developments by providing the analysis of a general framework for neural optimization of finite element methods. 
\par
%
% \subsubsection{Neural networks and PDEs / FEM}
\emph{Neural networks for PDEs}: The use of neural networks for approximating directly the solution to PDEs has received wide-spread interest since the works by E~\&~Yu~\cite{EYuCMS2018}, Sirignano \&~Spiliopoulos~\cite{SirSpiJCP2018}, Berg \&~Nystr\"om~\cite{BerNysNC2018} and Raissi, Perdikaris \&~Karniadakis~\cite{RaiPerKarJCP2019}, amongst others. Recently, there have been a number of ideas that propose an adaptive construction of neural-network approximations; see Ainsworth \&~Dong~\cite{AinDonSISC2021}, Liu, Cai \&~Chen~\cite{LiuCaiCheCAMWA2021} and Uriarte, Pardo \&~Omella~\cite{UriParOmeCMAME2022}. Neural networks can also be used to obtain the coefficients of the basis expansion used by a standard (linear) approximation~\cite{HesUbiJCP2018, KhaBalJosSarHegKriGan2021}. 
\par
\emph{Neural networks for inverse PDEs}:
In the context of inverse problems involving PDEs, the use of neural networks to represent unknown PDE~coefficients (fields) and constitutive models has been explored by, e.g., Teichert, Natarajan, Van der Ven \&~Garikipati~\cite{TeiNatVenGarCMAME2019}, Berg \&~Nystr\"om~\cite{BerNysJCMDS2021} and Xu \&~Darve~\cite{XuDarJCP2022}. These works are similar to the current work in the sense that standard (finite element) methods are used to solve the PDE, while a neural network is embedded within the discrete formulation. We note that the analysis provided by our current work can be extended to those inverse problems.
\par
\emph{Error analysis for neural-network approximations}:
There are a number of works containing a~priori error analysis for neural-network based PDE approximations. For those related to the deep Ritz method; see
Xu~\cite[Section~5]{XuCCP2020}, Pousin~\cite[Section~3]{PouJCMDS2022}, and M\"uller \&~Zeinhofer~\cite{MulZei2021}.  
For those related to physics-informed neural networks (PINN) and least-squares methods; see Sirignano \&~Spiliopoulos~\cite[Section~7]{SirSpiJCP2018}, Mishra \&~Molinaro~\cite{MisMolIMANA2022b, MisMolIMANA2022a}, Pousin~\cite[Section~4]{PouJCMDS2022} and Cai, Chen \&~Liu~\cite{CaiCheLiuJCP2021}. Recently, a~posteriori error analysis has also been studied, in particular goal-oriented analysis using the dual-weighted residual (DWR) methodology; see, e.g., 
Roth, Schr\"oder and Wick~\cite{RotSchWicSNAS2022}, 
Minakowski \&~Richter~\cite{MinRic2021} and Chakraborty, Wick, Zhuang \&~Rabczuk~\cite{ChaWicZhuRab2021}. We note that in our current work, while we have in mind the error analysis for neural-control approximations, the abstract analysis presented in Section~\ref{sec:Abstract} is essentially an extension of the above-mentioned a~priori analysis to a certain class of problems involving a convex and differentiable cost functional. 
\section{Abstract framework}
\label{sec:Abstract}
In this section we present the analysis of the abstract state equation (in the form of a mixed system) and the associated optimization problem.  
% parameterized by a neural-control variable  where the neural
%
We essentially follow the classical theory of optimal control (PDE-constrained optimization) by Lions~\cite{LioBOOK1971}; see also, \cite{HinPinUlbUlbBOOK2009, TroBOOK2010, BorSchBOOK2012}. Our resulting optimization problem bears similarity to that of parameter identification of PDE coefficients; see Rannacher \&~Vexler~\cite{RanVexSICON2005} and references therein for its error analysis. While we present our abstract framework within Hilbert spaces (and using a quadratic cost), we note that extensions to Banach spaces are feasible, but not within the scope of the current work.
\subsection{Discrete state problem and associated cost functional}
%The constrained optimization problem}
%Continuous and neural-discrete control problem
\label{ssec:Problem}
Let~$\mathbb{X}$ be a Hilbert space for the control variable,  $\mathbb{U}$ and $\mathbb{V}$ be 
Hilbert spaces for trial and test functions, respectively,  
$\mathbb{U}_h \subset \mathbb{U}$ be a discrete (finite element) subspace, and $\hat{\mathbb{V}}\subseteq \mathbb{V}$.%
\footnote{Later on, when considering minimal residual formulations, $\hat{\mathbb{V}}$ will be a discrete (finite element) subspace of  $\mathbb{V}$, but for the other formulations $\hat{\mathbb{V}} = \mathbb{V}$.}
In all that follows, we think of~$h$ (hence $\mathbb{U}_h$) as being fixed.
%While $\mathbb U_h$ is suppose to be a finite element space, the space $\hat{\mathbb V}$ is suppose to be either a finite element space, or the full space $\mathbb V$.   
%
Given~$\xi \in \mathbb{X}$ and $f\in\mathbb V^*$ (the dual of~$\mathbb{V}$), we consider the discrete state problem given by:
%be defined by the mixed system consider the mixed problem to  such that
%
\begin{subequations}
\label{eq:StateEq}
\begin{empheq}[left=\left\{,right=\right.,box=]{alignat=3}
\notag
& \text{Find } (r,u_h) \in \hat{\mathbb{V}} \times \mathbb{U}_h:\negqquad\negqquad\negqquad\negqquad
\\
\label{eq:StateEq_a}
 &\qquad a(\xi; r, v) + b(u_{h}, v) && = f( v), && \qquad\forall v\in \hat{\mathbb{V}},\\
\label{eq:StateEq_b}
 & \qquad b(w_h, r) && = 0, && \qquad \forall w_h\in \mathbb{U}_h\,,
\end{empheq}
\end{subequations}
where $b(\cdot,\cdot)$ is a continuous bilinear form on~$\mathbb{U}\times \mathbb{V}$, i.e., $b(\cdot,\cdot) \in \mathcal{L}(\mathbb{U} \times \mathbb{V};\mathbb{R})$, and 
for each~$\xi\in \mathbb{X}$,  $a(\xi;\cdot,\cdot)$ is a continuous bilinear form on~$\mathbb{V}\times \mathbb{V}$, i.e., 
$a(\xi;\cdot,\cdot)\in\mathcal{L}(\mathbb{V} \times \mathbb{V};\mathbb{R})$.
% ,
%
%\footnote{$\mathcal{L}(\mathbb{W}_1;\mathbb{W}_2)$ denotes the space of linear continuous maps from $\mathbb{W}_1$ to $\mathbb{W}_2$.
% }
%
%and . 
To explicitly indicate the dependence of $r$ and $u_h$ on $\xi$, we use the notation:
\begin{alignat*}{2}
  (r_\xi,u_{h,\xi}) = \text{solution of \eqref{eq:StateEq_a}--\eqref{eq:StateEq_b} for a given~}\xi\,.
\end{alignat*}
\par
In Section~\ref{sec:analysis}, we demonstrate that~\eqref{eq:StateEq_a}--\eqref{eq:StateEq_b} is equivalent to~\eqref{eq:introPGxi} for a particular choice of~$\mathbb{V}_h(\xi)$; see Proposition~\ref{prop:equiv}.
%
%\par
%
The discrete problem in~\eqref{eq:StateEq_a}--\eqref{eq:StateEq_b} is essentially a general formulation, which for a specific choice of $a(\cdot\,;\cdot,\cdot)$ and~$\hat{\mathbb{V}}$ reduces to a (weighted) Galerkin, least-squares or minimal residual method; see Section~\ref{sec:formulations}.  
\par
Next, let $\mathbb{Z}$ be a 
% \emph{observable}
Hilbert space, and let $Q : \mathbb{U} \rightarrow \mathbb{Z}$ be a  
linear continuous (observation) operator.
Then, given an observation $z_o\in \mathbb{Z}$ and regularization parameter $\alpha\ge 0$, we consider the cost (or loss) functional~$J:\mathbb{U}_h\times \mathbb{X}\rightarrow \mathbb{R}$ defined by:
\begin{alignat}{2}
\label{eq:defBigJ}
  J(w_h,\xi) := J_1(w_h) + \alpha\, j_2(\xi)\,,
  %\qquad\forall (w_h,\xi)\in\mathbb U_h\times\mathbb X\,, 
\end{alignat}
where
\begin{subequations}
\begin{empheq}[left=\left.,right=\right.,box=]{alignat*=3} 
%\label{eq:J_1}
J_1(w_h) &:= 
% g_1 \Big( \bignorm{ Q(u_{h,\xi}) - z_d }_{\mathbb{Z}} \Big)
% =  
\frac{1}{2} \bignorm{ Q(w_h) - z_o}_{\mathbb{Z}}^2 \,,
% \,, && \quad \forall w_h\in\mathbb U_h\,,
\\  
%\label{eq:j_2}
j_2(\xi) &:=
% g_2\Big(\norm{ \xi  }_{\mathbb{X}}\Big) 
%= 
\frac{1}{2} \norm{ \xi  }_{\mathbb{X}}^2 \,.
%\,, && \quad \forall \xi\in\mathbb X\,.
\end{empheq}
\end{subequations}
The associated \emph{reduced} cost functional $j:\mathbb{X}\rightarrow\mathbb{R}$ is then given by:
\begin{alignat}{2}
\label{eq:defj}
  j(\xi) := % J(u_{h,\xi},\xi) = J_1(u_{h,\xi}) 
  j_1(\xi) + \alpha\, j_2(\xi)\,,
\end{alignat}
where $j_1:\mathbb{X}\rightarrow \mathbb{R}$ is defined by:
\begin{alignat*}{2}
% \label{eq:j_1}
j_1(\xi) &:= J_1(u_{h,\xi}) = 
\frac{1}{2} \bignorm{ Q(u_{h,\xi}) - z_o}_{\mathbb{Z}}^2 \,,
\end{alignat*}
%
% consider the following control problem:
%
%subsection{The discrete control problem}
%
While ideally we would like to minimize~$j(\cdot)$ over (the infinite-dimensional)~$\mathbb{X}$, we proceed by considering neural-network approximations.
\subsection{Neural quasi-minimization}
\label{ssec:quasimin}
To accommodate neural optimization, we consider the subset $\mathcal{M}_n\subset\mathbb X$ consisting of all functions implemented by neural networks of a fixed architecture parameterized by~$n$.%
\footnote{In the terminology of Petersen, Raslan and Voigt~\cite{PetRasVoiFOCM2021}, the set~$\mathcal{M}_n$ consists of the \emph{realisations} of all possible neural networks of some fixed architecture (and some given activation function). While a neural network is identified with the set of weight and bias parameters, its realisation is the \emph{function} implemented by the network.}
We shall simply refer to $\mathcal{M}_n$ as a set of \emph{neural-network functions}, and we think of~$n$ as a measure of the size of the architecture (e.g., the total number of neurons, or total number of parameters). 
% The set $\mathcal{M}_n$ which are to be thought the functions generated by neural network with a (fixed) architecture (layers) that is parameterized by~$n$.%
%
\par
When aiming to minimize~$j(\cdot)$, a significant complication is that the set~$\mathcal M_n$ \emph{may not be closed} (topologically) in $\mathbb{X}$.%
\footnote{For example, \cite[Theorem~3.1]{PetRasVoiFOCM2021} shows that, under mild conditions on the architecture and activation function, $\mathcal{M}_n$ is not a closed subset of~$L^2(\Omega)$ (or, more generally, $L^p(\Omega)$, with $0<p<\infty$), unless, e.g., an upper bound is imposed on the weight parameters~\cite[Proposition~3.7]{PetRasVoiFOCM2021}.
}
Hence, even though $j(\cdot)$~may have an infimum on~$\mathcal{M}_n$, there may not be a minimizer in~$\mathcal{M}_n$. Therefore, 
% it is not advisable
one should not aim to \emph{completely} minimize $j(\cdot)$, % (over~$\mathcal{M}_n$)
but instead use a relaxed notion of \emph{quasi-minimization} as used by Shin, Zhang \&~Karniadakis~\cite{ShiZhaKar2020}%
\footnote{
%While 
Quasi-minimization can also be thought of as solving the minimization problem up to some optimization accuracy, cf.~\cite{MulZei2021}.
% its significance when minimizing over open subsets (such as when dealing with neural-network functions) should not be understated.
} 
(for which the existence of an infimum implies the existence of a quasi-minimizer): % cf.~\cite{ShiZhaKar2020, MulZei2021}):
\begin{definition}[Quasi-minimizers and quasi-minimizing sequences]
\label{def:quasimin}
~\\
Let $j:\mathbb{X}\rightarrow \mathbb{R}$ be a cost functional.
\begin{itemize}
    \item[(i)]
 Let~$\delta_n>0$ and $\mathcal{M}_n\subset \mathbb{X}$ be a subset of~$\mathbb{X}$ (not necessarily closed in~$\mathbb{X}$).  A function~$\bar{\xi}_n \in \mathcal{M}_n$ is said to be a \emph{quasi-minimizer} of~$j(\cdot)$ if the following holds true:%
\footnote{Observe that if $j(\cdot)$ has an infimum on~$\mathcal{M}_n$, then immediately a quasi-minimizer exists (in $\mathcal{M}_n$). This is true simply by the definition of the infimum.}
\begin{alignat}{2}
\label{eq:qmin}
 j(\bar{\xi}_n) \le \inf_{\xi_n\in \mathcal{M}_n} j(\xi_n) + \frac{\delta_n}{2}\,.
\end{alignat}
\item[(ii)]
Consider a sequence of subsets~$(\mathcal{M}_n)_{n\in \mathcal{N}}$ of~$\mathbb{X}$, with~$\mathcal{N}$ being a strictly-increasing sequence of natural numbers. A sequence $(\bar{\xi}_n)_n$, with $\bar{\xi}_n\in \mathcal{M}_n$, is said to be a \emph{quasi-minimizing sequence} if \eqref{eq:NCqminReduced}~holds true for all~$n\in \mathcal{N}$ with~$\delta_n>0$ such that:
% , and additionally, $(\delta_n)_{n\in\mathcal N}$ is an associated sequence of positive real numbers such that: 
%
\begin{alignat*}{2}
\delta_n \rightarrow 0 \quad \text{ as } \quad n\rightarrow \infty\,.
\end{alignat*}
%
% \begin{alignat*}{2}
%  j(\bar{\xi}_n) \le \inf_{\xi_n \in \mathcal{M}_n } j(\xi_n) + \frac{\delta_n}{2} \,.
% \end{alignat*}
%
%Problem~\eqref{eq:NCqminReduced} is said to be a \emph{quasi-minimization} of~$j(\xi_n)$ for~$\xi_n\in \mathcal{M}_n$,
\end{itemize}
\end{definition}
%
%
%a cost on a set of neural-network functions. 
% Instead, one way 
\par
In summary, the neural optimization problem that we consider is the following:
\begin{definition}[The quasi-minimizing control problem]
\label{def:qminprob}
~\\
The following statements are equivalent.
\\
%
% Let~$\delta_n > 0$. We are now ready to define 
\emph{Reduced quasi-minimizing control problem}: 
For $j(\cdot)$ given by~\eqref{eq:defj}, we aim to quasi-minimize~$j(\cdot)$, i.e., given $\delta_n>0$,
\begin{empheq}[left=\left\{,right=\right.,box=]{alignat=3}
\notag
 &\text{Find } \bar{\xi}_n\in \mathcal{M}_n:
  \\
 \label{eq:NCqminReduced}
&\qquad j(\bar{\xi}_n) \le \inf_{\xi_n\in \mathcal{M}_n} j(\xi_n) + \frac{\delta_n}{2}\,.
\end{empheq}
%
%find~$\bar{\xi}_n \in \mathcal{M}_n$ such that~\eqref{eq:NCqminReduced} holds, 
% where $j(\cdot)$ is defined in~\eqref{eq:defj}, or in other words, 
% (recalling the definition of~$j(\cdot)$ in~\eqref{eq:defj}), 
\emph{Constrained quasi-minimizing control problem}: For $J(\cdot,\cdot)$ given by~\eqref{eq:defBigJ},
we aim to quasi-minimize $J(u_h, \xi)$ subject to~\eqref{eq:StateEq_a}--\eqref{eq:StateEq_b}, i.e., given $\delta_n>0$,
\begin{empheq}[left=\left\{,right=\right.,box=]{alignat=3}
\notag
 &\text{Find } \bar{\xi}_n\in \mathcal{M}_n:
 \\
\label{eq:NCqmin}
 &\qquad  J(u_{h,\bar{\xi}_n},\xi_n) \le  \inf_{\xi_n \in \mathcal{M}_n} J(u_{h,\xi_n},\eta_n) + \frac{\delta_n}{2}\;.
\end{empheq}
%
%where $J(\cdot)$ is given in~\eqref{eq:defj}.
\end{definition}

\begin{example}[Need for quasi-minimizers]
% for PINN]
Let us discuss a simple example illustrating the non-existence of minimizers, hence the need for quasi-minimizers.%
\footnote{This is essentially an example of a PINN problem, i.e., minimizing a strong residual and boundary condition in least-squares sense. It is not difficult to construct a similar example for a neural control problem.}
\par
Let~$x=(x_1,x_2)\in\Omega = (0,1)^2 \subset \mathbb{R}^2$. Given $z\in (0,1)$, let~$\chi_{[z,1]}$ denote the characteristic function
of the subset~$[z,1]$.%
\footnote{That is, $\chi_{[z,1]}(x_1) = 1$ if $x_1 \in [z,1]$ and $=0$ otherwise.}
Consider the following cost functional:
%
% Discrete minimizers need not exists, e.g., consider for $z\in (0,1)$ the minimization of
\begin{alignat*}{2}
j(\xi) =  \frac{1}{2}\int_0^1\int_{0}^1 \bigg(\frac{\partial\xi}{\partial x_2}\bigg)^2\,\dd x_1 \,\dd x_2
+ \int_{0}^1 
\big( \xi- \chi_{[z,1]}\big)^2 \,\dd x_1
\end{alignat*} 
for $\xi\in \mathbb{X} = \Big\{ \eta \in L^2(\Omega) \,\big|\,\frac{\partial\eta}{\partial x_2} \in L^2(\Omega) \Big\}$. Minimizing $j(\cdot)$ over~$\mathbb{X}$ solves a first-order PDE (constant advection in the direction of the $x_2$-axis) with discontinuous data given by $\chi_{[z,1]}$, which is a well-posed problem~\cite{BocGunBOOK-CH2016}.
\par
Let $\mathcal{M}_n$ be the set of two-layer neural-network functions $\Omega\mapsto \mathbb{R}^2\mapsto \mathbb{R}$ using two neurons and ReLU activation in the hidden layer, i.e., %
\begin{alignat*}{2}
\mathcal{M}_n = \bigg\{ \xi_n: \Omega\rightarrow \mathbb{R} \,\Big|\, \xi_n(x) = \sum_{i=1}^2a_i \operatorname{ReLU}(w_i \cdot x - b_i)\,,\, a_i,b_i\in \mathbb{R}, w_i\in \mathbb{R}^2\bigg\}\,.
\end{alignat*}
Note that an infimizing sequence of~$j(\cdot)$ in~$\mathcal{M}_n$ is given by:
\begin{alignat*}{2}
 \xi_m(x) =\begin{cases}
 0\quad & 0\le x_1< z_m:=(1-\tfrac{1}{m})z\,,
 \\
 \dfrac{x_1-z_m}{z-z_m}\quad & z_m\le x_1< z\,,
 \\
 1 & z\le x_1 \le 1\,,
 \end{cases}
\end{alignat*}
for $m=1,2,3,\ldots,$ 
but whose limit~$\xi_m\rightarrow \bar{\xi}$ in $\mathbb{X}$ as $m\rightarrow \infty$ is a \emph{discontinuous} function (with $j(\bar{\xi}) = 0$). Therefore the infimizer~$\bar{\xi}$ does \emph{not} exist in~$\mathcal{M}_n \subset C(\overline{\Omega})$.
\par
On the other hand, quasi-minimizers $\bar{\xi}_n$
% for which $j(\bar{\xi}_n) \le \inf_{\eta_n\in \mathcal{M}_n} j(\eta_n) + \delta$ for any $\delta>0$, 
do exist in~$\mathcal{M}_n$, in particular, $\xi_m$ as defined above is a quasi-minimizer for~$m$ large enough.%
\footnote{Indeed, one can verify by direct calculation that $m$ must be such that
$
 \frac{1}{3}(z-z_m) \le \frac{\delta_n}{2}
 $, i.e.,  $ m \ge \frac{2}{3} z\delta_n^{-1}$.
}
\end{example}
\subsection{Analysis of reduced control problem}
We first proceed with the analysis of the reduced control problem~\eqref{eq:NCqminReduced}. % We first consider the reduced problem, which is obtained by eliminating $u_h$ from~$J$, and return to the constrained problem in Section~\ref{sec:analysis}.
% \par
Let the state operators~$R_h:\mathbb{X}\rightarrow \hat{\mathbb{V}}$ and $S_h:\mathbb{X} \rightarrow \mathbb{U}_h$ be defined by:
\begin{subequations}
\label{eq:state_operator}
\begin{alignat}{2}
 R_h(\xi) := r_{h,\xi}, \qquad \forall \xi\in \mathbb{X}\,,
\\
 S_h(\xi) := u_{h,\xi}, \qquad \forall \xi\in \mathbb{X}\,,
\end{alignat}
\end{subequations}
where $r_{h,\xi}$ and $u_{h,\xi}$ are the first and second component, respectively, of the solution to the mixed system~\eqref{eq:StateEq}. 
Then the reduced cost~$j(\cdot)$ given in~\eqref{eq:defj} can be written as follows:
%
% \begin{alignat}{2}\label{eq:j_1}
%  j_1(\cdot) = J_1(S_h(\cdot),\cdot) 
% \end{alignat}
% %
% our interest turns to the quasi-minimization of the reduced cost:
\begin{alignat}{2}
\notag
  j(\xi) = j_1(\xi) + \alpha\, j_2(\xi)
  &= J_1\big(S_h(\xi),\xi\big)
+  \alpha j_2(\xi)
\\
\label{eq:reducedcost}
 &= \frac{1}{2}
 \bignorm{ Q \circ S_h(\xi)  - z_o }_{\mathbb{Z}}^2
 + \frac{ \alpha }{2}
 \norm{ \xi }_{\mathbb{X}}^2
 \,.
\end{alignat}
%then problems~\eqref{eq:JminX}~and~\eqref{eq:JminMn} are equivalent (respectively) to the unconstrained minimizations
% we then have the reduced quasi-minimization problem:
% our interest turns into the following unconstrained minimization problems:
%
% \begin{subequations}
% \begin{empheq}[left=\left\{,right=\right.,box=]{alignat=3} 
% \label{eq:jminX}
% & \min_{\xi\in\mathbb X} j(\xi) := j_1(\xi) + \alpha\, j_2(\xi),\\
% \label{eq:jminMn}
% & \min_{\xi_n\in\mathcal M_n} j(\xi_n).
% \end{empheq}
% \end{subequations}
%
%
% \marginnote{What if S is near-linear: i.e., linear + perturbation. Small $\alpha$  works?}
\par
Our main result depends on the following fundamental theorem, which is of independent interest:
\begin{theorem}[Differentiable, strongly-convex quasi-minimization]
\label{thm:quasi}
~\\
Let $j:\mathbb{X}\rightarrow \mathbb{R}$ be a cost functional. Assume that $j(\cdot)$ is G\^ateaux differentiable with derivative $j':\mathbb{X}\rightarrow \mathbb{X}^*$ being Lipschitz continuous, i.e., there is a constant~$L>0$ such that
\begin{alignat*}{2}
\bignorm{j'(\xi) - j'(\eta)}_{\mathbb{X}^*}
\le L \bignorm{\xi-\eta}_{\mathbb{X}}\,, 
\qquad \forall \xi,\eta\in\mathbb{X}\,,
\end{alignat*}
Furthermore, assume that $j(\cdot)$ is strongly convex, i.e., there is a constant~$\gamma>0$ such that
\begin{alignat}{2}
\label{eq:strongConvex}
\Bigdual{j'(\xi) - j'(\eta)\,,\,\xi-\eta}_{\mathbb{X}^*,\mathbb{X}}
\ge \gamma \bignorm{\xi-\eta}_{\mathbb{X}}^2\,, 
\qquad \forall \xi,\eta\in\mathbb{X}\,.
\end{alignat}
Then the following hold true:
\begin{itemize}
    \item[(i)] $j(\cdot)$ has a unique minimizer~$\bar{\xi} \in \mathbb{X}$, 
%to problem~\eqref{eq:jminX}, 
which satisfies:
\begin{alignat*}{2}
 j'(\bar{\xi}) = 0 \qquad \text{in } \mathbb{X^*}\,.
\end{alignat*}
\item[(ii)] For any subset~$\mathcal{M}_n\subset \mathbb{X}$, $j(\cdot)$ has a quasi-minimizer~$\bar{\xi}_n\in \mathcal{M}_n$ that satisfies~\eqref{eq:qmin}.
%problem~\eqref{eq:NCqminReduced} has a quasi-minimizer~$\bar{\xi}_n \in \mathcal{M}_n$.
% to problem~\eqref{eq:jminMn}. % satisfying
%\marginnote{18 April, @Ignacio: I have made a subtle, but important(!) change to~(iv).}
\item[(iii)]
     Any quasi-minimizer~$\bar{\xi}_n$ in $\mathcal{M}_n$ satisfies the following quasi-optimal error estimate:
\begin{equation}\label{eq:quasi_opt_abstract}
\bignorm{\bar\xi-\bar\xi_n}_{\mathbb X}\le
\bigg( \frac{L}{\gamma}\inf_{\xi_n\in\mathcal M_n}
\bignorm{\bar\xi-\xi_n}_{\mathbb X}^2
+ \frac{\delta_n}{\gamma} \bigg)^{1/2}\,.
\end{equation}
\end{itemize}
\end{theorem}
\begin{proof}
See Appendix~\ref{sec:quasi_proof}.
\end{proof}
\par
We now analyse when our~$j(\cdot)$ satisfies the assumptions of Theorem~\ref{thm:quasi}. 
\begin{theorem}[Reduced control problem: Differentiability \& strong convexity]
\label{thm:reduced}%
~\\
%Consider an observation operator 
Let $\alpha > 0$ and $j(\cdot) = j_1(\cdot) + \alpha\, j_2(\cdot)$ be as in~\eqref{eq:reducedcost}. Let $Q\in \mathcal{L}(\mathbb{U},\mathbb{Z})$.
%together with 
% a differentiable state operator 
Assume $S_h:\mathbb X\to\mathbb U_h$ is differentiable,  $S_h(\cdot)$ and $S_h'(\cdot)$ are uniformly bounded on $\mathbb X$, and $S_h'(\cdot)$ is Lipschitz continuous. Then: 
%Consider the $j_1$ functional defined in~\eqref{eq:j_1}, the $j_2$ functional defined in~\eqref{eq:j_2}, and $j:=j_1+\alpha j_2$.
% \marginnote{Or S unif bounded and Lipschitz (implies S' is uniformly bounded), and S' Lipschitz}
% Assume $j_1$ is continuous, $j_1'$ Lipschitz continuous.
% , $j_2$ continuous and strongly convex.
\begin{itemize}
\item[(i)] $j_1,j_2,j:\mathbb{X}\rightarrow \mathbb{R}$ are G\^ateaux differentiable with
$j_1', j_2', j':\mathbb{X} \rightarrow \mathbb{X}^*$
% $j_1'(\cdot)$ $j_2'(\cdot)$ and $j(\cdot)$ 
Lipschitz continuous.
% (hence also $j'$). 
\end{itemize}
Additionally, assume $\alpha$ is sufficiently large. 
Then:
% Then the following statement holds true:
\begin{itemize}
\item[(ii)] $j:\mathbb{X}
\rightarrow \mathbb{R}$ is strongly convex, i.e., there is a constant~$\gamma>0$ such that~\eqref{eq:strongConvex} holds true.%
\footnote{In particular, when $\alpha > L_1$, where~$L_1$ is the Lipschitz constant of~$j_1'(\cdot)$, then $\gamma = \alpha - L_1$.}
%
% \begin{alignat*}{2}
% \Bigdual{j'(\xi) - j'(\eta)\,,\,\xi-\eta}_{\mathbb{X}^*,\mathbb{X}}
% \ge \gamma \bignorm{\xi-\eta}_{\mathbb{X}}^2\,.
% \end{alignat*}
%
\end{itemize}
\end{theorem}
\begin{proof}
See Appendix~\ref{sec:reduced_proof}.
\end{proof}
\begin{corollary}[Reduced control problem: (Quasi-)minimizers \& quasi-optimality]
\label{cor:reducedProblem}
~\\
Under the conditions of Theorem~\ref{thm:reduced}, the statements~(i), (ii) and (iii) of Theorem~\ref{thm:quasi} hold true.
\end{corollary}
\begin{proof}
The results of Theorem~\ref{thm:reduced} are the assumptions of Theorem~\ref{thm:quasi}.
\end{proof}

\begin{remark}[Quasi-optimal rates]
The first part on the right-hand side of the quasi-optimality result~\eqref{eq:quasi_opt_abstract} can be estimated in terms of~$n$ using results from neural-network approximation theory; see, e.g., Yarotsky~\cite{YarNN2017}, G\"uhring, Kutyniok and Petersen~\cite{GuhKutPetAA2020}, and references therein. Such a result may be useful in finding a proper balance of~$\delta_n$ as~$n\rightarrow \infty$. Alternatively, the choice of~$\delta_n$ may be found through a proper \emph{a~posteriori} estimator, which seems to be an open problem. %  as far as we are aware.
\end{remark}
\begin{remark}[Condition on~$\alpha$]
The proof of Theorem~\ref{thm:reduced} reveals that the condition that $\alpha$ is sufficiently large may be weakened if $j_1$ has additional structure (e.g., convexity).
%; see,  e.g., ~\cite{BarBOOK2016,BreBOOK2011,DacBOOK2008}.
Indeed, convexity of $j_1$ guarantees that $j$ will be strongly convex, with strongly convexity constant equal to $\alpha>0$. If the case, there is no need of Lipschitzness of $j_1'$ in order to prove statement~(iii) of Theorem~\ref{thm:reduced}, only $\alpha>0$ will be enough. % However, Lipschitzness of $j_1'$ will be required to show 
Furthermore, statement~(v) of Theorem~\ref{thm:reduced} becomes:
% , in which case the estimate~\eqref{eq:quasi_opt} translates to
$$
\|\bar\xi-\bar\xi_n\|_{\mathbb X}<
\bigg( \frac{\alpha+L_1}{\alpha}\inf_{\eta_n\in\mathcal M_n}
\bignorm{\bar\xi-\eta_n}_{\mathbb X}^2
+ \frac{\delta_n}{\alpha} \bigg)^{1/2}.
$$
\end{remark}
\begin{remark}[Physics-informed neural networks (PINN)]
% Connections with other neural-network approaches]
~\\
Theorem~\ref{thm:quasi} can be applied to PINN~\cite{RaiPerKarJCP2019} (for neural-network approximations to PDEs). Indeed, consider
% $\mathbb X=\mathbb U = \mathbb{U}_h$, $S_h=\operatorname{id}$, 
% % $\mathbb Z= L^2(\Omega)$, 
% $Q=B$, and $z_o=f$, 
% %\in L^2(\Omega)$, 
% one obtains the residual in~$j_1(\cdot)$:
$$
j(\xi) 
= \frac{1}{2}\bignorm{f- B\xi}_{\mathbb{L}}^2\,,
$$
where $f- B\xi$ is an abstract residual in some abstract Hilbert space~$\mathbb{L}$ (which may include the PDE residual, initial condition and boundary conditions, as in~\cite{MisMolIMANA2022b}, as well as a data residual, as in~\cite{MisMolIMANA2022a}). If~$B:\mathbb{X}\rightarrow \mathbb{L}$ is a linear operator, then the assumptions of Theorem~\ref{thm:quasi} (Lipschitz continuity and strong convexity) hold true.
%\par
%
% Not only parametrized finite-element approaches fit into the framework of Theorem~\ref{thm:reduced}, and in particular, into the definition of $j_1(\cdot)$ given in~\eqref{eq:j_1}. For instance, if $B:\mathbb U_B\to L^2(\Omega)$ is a differential operator defined on the graph space
% $$
% \mathbb U_B:=\{u\in L^2(\Omega)\, :\, Bu\in L^2(\Omega) + \hbox{ boundary conditions}\}\,,
% $$
% then, making $\mathbb X=\mathbb U=\mathbb U_B$, $S=\operatorname{id}$, $\mathbb Z= L^2(\Omega)$, $Q=B$, and $z_o=f\in L^2(\Omega)$, we get
% $$
% j_1(\xi) = \frac{1}{2}\|Q\circ S(\xi) - z_o \|_{\mathbb Z}^2
% = \frac{1}{2}\|B\xi - f \|_{L^2}^2\,.
% $$
% Observe that Theorem~\ref{thm:reduced} immediately applies for such a $j_1$-functional, which can be interpreted as an abstract form of the PINN method~\cite{RaiPerKarJCP2019}, or the DGM method~\cite{SirSpiJCP2018}, when homogeneous boundary/initial conditions are imposed strongly.   
% Of course, we can complexify the formulation in terms of $j_1(\cdot)$ and the underlying functional spaces to handle non-homogeneous boundary/initial conditions as well.
\end{remark}
\begin{remark}[Deep Ritz method]
~\\
Theorem~\ref{thm:quasi} can also be applied to the Deep Ritz method~\cite{EYuCMS2018}. Indeed, consider 
%With minor modifications, Theorem~\ref{thm:reduced} can be applied to the Deep Ritz method~\cite{EYuCMS2018}. Indeed, define $j_1:\mathbb X\to \mathbb R$ as follows:
$$
j(\xi) = \frac{1}{2} b(\xi,\xi) - f(\xi)\,,
$$
where $b\in\mathcal L(\mathbb X\times\mathbb X;\mathbb R)$ is a coercive bilinear form and $f\in\mathbb X^*$. For such a $j(\cdot)$, the assumptions of Theorem~\ref{thm:quasi} (Lipschitz continuity and strong convexity) hold true.
%
%is differentiable (with Lipschitz-continuous $j'_1$) and strongly convex.
%and coercive. 
%Hence, Theorem~\ref{thm:reduced}(iii)-(v) applies.
% for such a $j(\cdot)$ by making $\mathbb X=\mathbb U$.
%
\end{remark}

\subsection{Analysis of constrained control problem}
\label{sec:analysis}
We now proceed with the analysis of the \emph{constrained} control problem~\eqref{eq:NCqmin}. We begin by providing conditions that guarantee the well-posedness of the state problem.
\begin{proposition}[Stability of the state problem]
\label{prop:stateEq}
Let $a(\xi;\cdot,\cdot)\in \mathcal{L}(\mathbb{V} \times \mathbb{V};\mathbb{R})$ for each~$\xi\in \mathbb{X}$, 
% Let us consider a parameter-dependent bilinear form 
% $$\mathbb X\ni \xi \mapsto a(\xi;\cdot,\cdot)\in \mathcal{L}(\mathbb{V} \times \mathbb{V};\mathbb{R})\,,$$ and let 
and let $b(\cdot,\cdot) \in \mathcal{L}(\mathbb{U} \times \mathbb{V};\mathbb{R})$. 
% Given conforming subspaces 
For $\mathbb U_h\subset \mathbb U$ and $\hat{\mathbb V}\subseteq\mathbb V$, let the kernel subspace
$
\hat{\mathbb{K}}:=\{v\in\hat{\mathbb V}: b(w_h,v)=0, \forall w_h \in\mathbb U_h\}.
$
% Let $\mathbb{U}_h,\mathbb{V}$ reflexive Banach. 
%Given~$\xi\in \mathbb{X}$, let~$a(\xi;\cdot,\cdot)$ and $b(\cdot,\cdot)$ be continuous. 
%Let $\hat{\mathbb{K}}:=\{v\in\hat{\mathbb V}: b(w_h,v)=0, \forall w_h \in\mathbb U_h\}$.
Then, the following statements hold true:
\begin{itemize}
    \item[(i)] For each~$\xi\in\mathbb X$, problem~\eqref{eq:StateEq} is well-posed (for any~$f\in \mathbb{V}^*$) if and only if
there exist constants $\alpha_h\equiv\alpha_h(\xi)>0$ and $\beta_h>0$ such that:%
\footnote{%
\label{ftnt:a-surjectivity}
Only when $\hat{\mathbb V}$ is infinite-dimensional, one needs the extra hypothesis in~\eqref{eq:inf-sup_a}$_2$. 
% (cf.~Remark~\ref{rem:a-surjectivity}):
% $\big\{v_2\in\hat{\mathbb{K}}: a(\xi;v_1,v_2)=0,\, \forall v_1\in\hat{\mathbb{K}}\big\}=\{0\}$.
Whenever $a(\xi,\cdot,\cdot)$ is an equivalent inner product on $\mathbb V$, then this condition 
is actually automatically satisfied. Indeed, zero is the only element in $\mathbb V$ which is orthogonal to itself.
}
% 
%\begin{itemize}
    %\item 
    \begin{subequations}
    \label{eq:inf-sup}
    \begin{alignat}{2}
    \label{eq:inf-sup_a}
\left.
\begin{alignedat}{2}
    \displaystyle\inf_{v_1\in\hat{\mathbb{K}}}\sup_{v_2\in\hat{\mathbb{K}}} \frac{a(\xi;v_1,v_2)}{\|v_1\|_\mathbb V\|v_2\|_\mathbb V} \geq \alpha_h\,,
    \\
\big\{v_2\in\hat{\mathbb{K}}: a(\xi;v_1,v_2)=0,\, \forall v_1\in\hat{\mathbb{K}}\big\}=\{0\}\,,
\end{alignedat}
\qquad \right\}
    %\item  
    \\
    \label{eq:inf-sup_b}
    % \qquad\hbox{ and }\qquad
    \displaystyle\inf_{w_h\in\mathbb U_h}\sup_{v\in\hat{\mathbb V}}  \frac{b(w_h,v)}{\|w_h\|_{\mathbb U}\|v\|_\mathbb V} \geq \beta_h\,. \qquad 
    \end{alignat}
    \end{subequations}
%\end{itemize}
\item[(ii)]
%Assume $0<\underline{a} \le \norm{a}\le \bar{a}$ uniformly in~$\xi$. 
If~\eqref{eq:inf-sup} is satisfied, then the following a priori bound holds true for the solution $u_{h}\in\mathbb U_h$ of problem~\eqref{eq:StateEq}: 
$$
  \norm{u_{h}}_\mathbb{U}= \le \frac{1}{\beta_h}\left(
  1+\frac{\|a(\xi;\cdot,\cdot)\|_{\mathcal{L}(\mathbb{V} \times \mathbb{V};\mathbb{R})}}{\alpha_h}\right)\norm{f}_{\mathbb{V}^*}\,. 
$$
\item[(iii)] Furthermore, if $a(\xi,\cdot,\cdot)$ is an equivalent inner-product on $\mathbb V$, with associated norm $\|\cdot\|_{\mathbb V,\xi}:=\sqrt{\smash[b]{a(\xi;\cdot,\cdot)}}$, i.e., for some $C_{1,\xi},C_{2,\xi}>0$, 
% \footnote{For some equivalence constants $C_{1,\xi}>0$ and $C_{2,\xi}>0$.}
\begin{equation}\label{eq:norm_equiv}
C_{1,\xi} \|v\|_{\mathbb V}\leq \|v\|_{\mathbb V,\xi}
\leq C_{2,\xi} \|v\|_{\mathbb V}, \quad\forall v\in\mathbb V,
\end{equation}
then %we can take 
$\alpha_h=(C_{1,\xi})^2$ in~\eqref{eq:inf-sup_a}, and additionally, the following improved a~prior bound holds true:
\begin{equation}\label{eq:apriori}
  \norm{u_{h}}_\mathbb{U} \le \frac{C_{2,\xi}}{ C_{1,\xi}}\frac{1}{\beta_h}\norm{f}_{\mathbb{V}^*}\,. 
\end{equation}
\end{itemize}
\end{proposition}
\begin{proof}
See Appendix~\ref{sec:wellposed_proof}.
\end{proof}
\par
To establish the equivalence between the mixed system~\eqref{eq:StateEq} and the Petrov--Galerkin statement~\eqref{eq:introPGxi}, 
% it will be helpful to introduce the following operators.
% We introduce operator notation to facilitate the presentation of the next result.
let us define the operators $A:\mathbb X\to \mathcal L(\hat{\mathbb V},\hat{\mathbb V}^*)$ and $B\equiv B_h\in \mathcal{L}(\mathbb U_h;\hat{\mathbb V}^*)$ by:
% the operators such that
\begin{subequations}
\begin{alignat}{2}
 \label{eq:A_operator} 
A(\xi) \hat v := \, & a(\xi;\hat v,\cdot)\in\hat{\mathbb V}^*, & \quad\forall \xi\in\mathbb X,\, \forall \hat v\in\hat{\mathbb V};
\\
 \label{eq:B_operator} 
B w_h :=\, & b(w_h,\cdot) \in\hat{\mathbb V}^*, & \forall w_h\in \mathbb U_h.
\end{alignat}
\end{subequations}
Note that the state equations~\eqref{eq:StateEq_a}--\eqref{eq:StateEq_b} can then be written as follows:
% Thus, Proposition~\ref{prop:stateEq} establishes the well-posedness of the problem to find $\hat r\equiv~R(\xi)\in\hat{\mathbb V}$ and $u_h\equiv S_h(\xi)\in\mathbb U_h$ such that
\begin{subequations}\label{eq:mixed_operator}
\begin{alignat}{3}
%\left\{
%\begin{array}{clll}
\label{eq:mixed_operator_a}
A(\xi)  r & + Bu_h  && = f  \qquad && \text{in } \hat{\mathbb V}^*\,,
\\
\label{eq:mixed_operator_b}
B^*  r &  && =0 && \text{in } (\mathbb{U}_h)^*\,.
%\end{array}\right.
\end{alignat}
\end{subequations}
\begin{proposition}[Equivalent Petrov--Galerkin problem]
\label{prop:equiv}
~\\
Assume the conditions of Proposition~\ref{prop:stateEq}, including the well-posedness condition~\eqref{eq:inf-sup_b}. Instead of~\eqref{eq:inf-sup_a}, assume the stronger hypothesis (full inf-sup, instead of just on the kernel):
\begin{subequations}
\label{eq:inf-sup_a2}
\begin{alignat}{2}
    \displaystyle\inf_{v_1\in\hat{\mathbb{V}}}\sup_{v_2\in\hat{\mathbb{V}}} \frac{a(\xi;v_1,v_2)}{\|v_1\|_\mathbb V\|v_2\|_\mathbb V} \geq \alpha_h\,,
    \\
\Big\{v_2\in\hat{\mathbb{V}}\,:\, a(\xi;v_1,v_2)=0,\, \forall v_1\in\hat{\mathbb{V}}\Big\}=\{0\}\,,
\end{alignat}
\end{subequations}
%Under the conditions of Proposition~\ref{prop:stateEq}, assume there are constants $\alpha_h>0$ and $\beta_h>0$ such that \eqref{eq:inf-sup_a}--\eqref{eq:inf-sup_b} hold true. 
Let the test space $\mathbb{V}_h(\xi)$ be given by:
\begin{alignat}{2}
\label{eq:Vhxi}
  \mathbb{V}_h(\xi) 
  %&=  A(\xi)^{-*} B \mathbb{U}_h
  % \\
  &= \Big\{ v \in \mathbb{V} \,\Big|\, A(\xi)^* v =  B w_h\text{ for some } w_h \in \mathbb{U}_h \Big\}\,.
\end{alignat}
Then the state problem~\eqref{eq:StateEq} is equivalent to the Petrov--Galerkin problem~\eqref{eq:introPGxi} with~$\mathbb{V}_h(\xi)$ given by~\eqref{eq:Vhxi}.
\end{proposition}
\begin{proof}
See Appendix~\ref{sec:equiv_proof}.
\end{proof}
\par
Finally, we now present (differentiability) conditions on $\xi\mapsto A(\xi)$ that guarantee the (differentiability) requirements on~$\xi \mapsto S_h(\xi)$ in Theorem~\ref{thm:reduced} and Corollary~\ref{cor:reducedProblem}. Once in place, existence of (quasi)-minimizers and quasi-optimal convergence follow immediately for the constrained control problem.
\par
To anticipate the connection between derivatives~$A'$ and $S_h'$ (as well as $R_h'$),% 
\footnote{Recall that the G\^ateaux derivative of, e.g.,~$A$ at~$\xi\in \mathbb{X}$ in the direction~$\eta\in \mathbb{X}$ is given by 
%
%\begin{alignat*}{2}
$
 A'(\xi)\eta 
  = \displaystyle{ \lim_{t\rightarrow 0}
   \frac{A(\xi+t\eta) - A(\xi)}{t} }
%  \,,
$,
%\end{alignat*}
%
provided the limit exists in~$\mathcal L(\hat{\mathbb V},\hat{\mathbb V}^*)$.  If the map $\eta\mapsto A'(\xi)\eta$ is linear and continuous from $\mathbb X$ to $\mathcal L(\hat{\mathbb V},\hat{\mathbb V}^*)$, then $A$ is G\^ateaux differentiable at $\xi\in\mathbb X$.
}
note that a formal differentiation of~\eqref{eq:mixed_operator} (with $r=R_h(\xi)$ and $u_h = S_h(\xi)$) with respect to~$\xi$ in the direction~$\eta\in \mathbb{X}$ yields:
\begin{alignat*}{3}
A(\xi)  R_h'(\xi)\eta & + BS_h'(\xi)\eta  && = -A'(\xi)\eta \,R_h(\xi)  \qquad && \text{in } \hat{\mathbb V}^*\,,
\\
B^*  R_h'(\xi)\eta &  && =0 && \text{in } (\mathbb{U}_h)^*\,.
%\end{array}\right.
\end{alignat*}
One may therefore expect that suitable conditions on~$A(\cdot)$ will imply desired conditions on~$S_h(\cdot)$ (and $R_h(\cdot)$):
% \par
% %
% To better appreciate the connection between the derivatives 
% %One can verify that formally the pair~$(r',u_h') = \big(r'(\xi)(\eta), u_h'(\xi)(\eta)\big)$ satisfies
% \begin{subequations}
% \label{eq:StateEqDeriv}
% \begin{alignat}{2}
%   a(\xi;r',v) + b(u_{h}',v) &= -a'(\xi;r,v)(\eta)
%   \qquad && \forall v\in \hat{\mathbb{V}}
%   \\
%   b(w_h,r') &= 0 && \forall w_h\in \mathbb{U}_h
% \end{alignat}
% \end{subequations}
%
% where
% \begin{alignat*}{2}
% a'(\xi;r,v)(\eta) := \lim_{t\rightarrow 0}
% \frac{ a(\xi + t\eta;r,v) - a(\xi;r,v) }{t}
% \end{alignat*}
%
\begin{proposition}[State differentiability] 
\label{prop:state_deriv}
~\\
Let~$R_h(\cdot)$ and $S_h(\cdot)$ be the state operators as defined in~\eqref{eq:state_operator}, and let~$A(\cdot)$ be as defined in~\eqref{eq:A_operator}.  
Assume the conditions of Proposition~\ref{prop:stateEq}, including the well-posedness conditions~\eqref{eq:inf-sup}. 
% Given~$\xi\in\mathbb X$, let $\hat r\equiv R(\xi)\in\hat{\mathbb V}$ and $u_h \equiv S_h(\xi)\in\mathbb U_h$ be the solutions of the mixed system~\eqref{eq:mixed_operator}.
Then, the following statements hold true:
\begin{itemize}
\item[(i)]
If $A(\cdot)$ has a G\^ateaux derivative at~$\xi\in \mathbb{X}$ in the direction~$\eta\in \mathbb{X}$, then $R_h(\cdot)$ and $S_h(\cdot)$ have a G\^ateaux derivative at~$\xi$ in the direction~$\eta$. 
%Moreover, the G\^ateaux derivative pair 
%$\big(\hat r'(\xi)\eta,u_h'(\xi)\eta\big)\in\hat{\mathbb V}\times\mathbb U_h$ is the unique solution of the mixed 
%system~\eqref{eq:mixed_operator} with right hand side 
%$f=-[A'(\xi)\eta]\hat r\in\hat{\mathbb V}^*$.
\item[(ii)] If $A(\cdot)$ is G\^ateaux-differentiable at $\xi$, then so are $R_h(\cdot)$ and $S_h(\cdot)$. 
\item[(iii)] If $A(\cdot)$, $A'(\cdot)$ and $\alpha_h^{-1}(\cdot)$ are uniformly bounded on $\mathbb X$, then $R_h'(\cdot)$ and $S_h'(\cdot)$ are also uniformly bounded on $\mathbb X$. 
\item[(iv)] Additionally, if $A'(\cdot)$ is Lipschitz continuous, then $R_h'(\cdot)$ and $S_h'(\cdot)$ are Lipschitz continuous as well.
%\item[(iv old)] 
%If $\hat{A}_\xi$ is uniformly invertible and $\hat{A}'$ is uniformly bounded, then $S'(\xi)$ is uniformly bounded.
\end{itemize}
\end{proposition}

\begin{proof}
See Appendix~\ref{sec:state_deriv_proof}.
\end{proof}

\begin{corollary}[Constrained problem: (Quasi-)minimizers \& quasi-optimality]
\label{cor:constrainedProb}
~\\
Let $J(w_h,\xi) = J_1(w_h) + \alpha \,j_2(\xi)$ as in~\eqref{eq:defBigJ} with~$Q\in \mathcal{L}(\mathbb{U};\mathbb{Z})$. Let the associated~$j(\cdot)$ be as in~\eqref{eq:reducedcost}. 
% $J(w_h,\xi) = J_1(w_h) + \alpha\, j_2(\cdot)$ be as~\eqref{eq:defBigJ}.
Under the conditions of Propositions~\ref{prop:stateEq} and~\ref{prop:state_deriv}, and assuming $\alpha$ is sufficiently large, the statements~(i), (ii) and (iii) of Theorem~\ref{thm:quasi} hold true. 
\\
In other words,
the constrained control problem~\eqref{eq:NCqmin} has a quasi-minimizer in~$\mathcal{M}_n$ that converges quasi-optimally to the unique minimizer in~$\mathbb{X}$.
%
%Assume formulation~\eqref{eq:StateEq} is well-posed (i.e., it satisfies Proposition~\ref{prop:stateEq}(i)). Assume that the induced operator 
% $A:\mathbb X\to\mathcal L(\hat{\mathbb V},\hat{\mathbb V}^*)$ defined in~\eqref{eq:A_operator} satisfies all the hypothesis (i)-(iv) of Proposition~\ref{prop:state_deriv}. Then, the state operator $S_h:\mathbb X\to\mathbb U_h$ (defined in~\eqref{eq:state_operator}) satisfy all the hypothesis of Theorem~\ref{thm:reduced}. In particular, the minimization problem~\eqref{eq:jminX} has a unique solution if $\alpha$ is large enough, and there exists quasi-minimizers of the minimization problem~\eqref{eq:jminMn}. Moreover, quasi-minimizers are quasi-optimal in the sense of equation~\eqref{eq:quasi_opt}. 
\end{corollary}
\begin{proof}
The results of Propositions~\ref{prop:stateEq} and~\ref{prop:state_deriv}, together with $\alpha$ sufficiently large, are the assumptions of Theorem~\ref{thm:reduced}, whose results are the assumptions of Theorem~\ref{thm:quasi}.
\end{proof}

\section{Conforming weak formulations with suitable control} 
% in Hilbert spaces}
\label{sec:formulations}
In this section, we study various weighted versions of conforming weak formulations, viz., least-squares, Galerkin and minimal-residual formulations. The aim is to propose suitable $\xi$-dependent weighting within the weak forms, in order to be able to prove the assumptions of Propositions~\ref{prop:stateEq} and~\ref{prop:state_deriv}. By Corollary~\ref{cor:constrainedProb}, we can then conclude that the corresponding constrained neural-control problem has desired properties (existence of quasi-minimizers and quasi-optimal convergence).
\par
%
% Duality pairings will be just denoted by $\left<\,\cdot\,,\,\cdot\,\right>$. Their precise meaning will be clear by the context in which they will be used.
%
In what follows, we often consider a positive weight function~$\omega$. We shall use the notation $\varpi:=1/\omega$ to indicate the (multiplicative) inverse of $\omega$.
\subsection{Weighted least-squares formulations}
%
% To make things precise, we
Let $d\in\mathbb N$ and $\Omega\subset\mathbb R^d$ be an open bounded domain. 
% Consider the Hilbert space $L^2(\Omega)$ of square-integrable measurable functions on $\Omega$, whose inner product will be denoted by $(\cdot,\cdot)_{L^2(\Omega)}$, and 
% corresponding norm by $\|\cdot\|_{L^2(\Omega)}$.
% Let us denote the $L^2(\Omega)$ inner product by
%
%
% \par
%
Let $B:\mathbb{H}_B\to L^2(\Omega)$ be a linear differential operator in strong form, where $\mathbb{U} = \mathbb{H}_B$ denotes the graph space
$$
\mathbb{H}_B:=\Big\{w\in L^2(\Omega)\,\Big|\, Bw\in L^2(\Omega)+\hbox{boundary conditions}\Big\}.
$$
%Additionally, $\mathbb{H}_B$ may also consider homogeneous boundary conditions.
We further assume that $\mathbb{H}_B$ is a Hilbert space when endowed with the inner product 
$$
\biginnerprod{w_1, w_2}_{\mathbb{H}_B}
:=\biginnerprod{w_1,w_2}_{L^2(\Omega)}
+ \biginnerprod{Bw_1,Bw_2}_{L^2(\Omega)}\, , \qquad \forall w_1,w_2\in\mathbb{H}_B\,,
$$
and that $B$ is boundedly invertible from $\mathbb{H}_B$ onto $\mathbb V^*:=L^2(\Omega)=:\mathbb V = \hat{\mathbb{V}}$. 

Given $f\in L^2(\Omega)$, a positive weight function $\omega:L^2(\Omega)\to L^\infty(\Omega)$, a control $\xi\in \mathbb{X} = L^2(\Omega)$, and a conforming discrete finite element space $\mathbb U_h\subset\mathbb{H}_B$, we aim to find $u_h\equiv S_h(\xi) \in 
\mathbb{U}_h$, which is the solution of the weighted least-squares problem:
$$
u_h = \arg\!\!\min_{w_h\in\mathbb U_h}\frac{1}{2}
\left\|\sqrt{\omega(\xi)}\big(f-Bw_h\big)\right\|^2_{L^2(\Omega)}\,.
$$
The optimality condition of such a minimizer is:
\begin{equation}\label{eq:lsq_form}
\Big(\omega(\xi)(f-Bu_h),Bw_h\Big)_{L^2(\Omega)} = 0, \quad \forall w_h\in\mathbb U_h\,.
\end{equation}
In particular, notice that we can directly identify the test space in~\eqref{eq:introPGxi} as $\mathbb{V}_h(\xi) = \omega(\xi) B \mathbb{U}_h = \big\{ v \in L^2(\Omega)\,\big|\, v = \omega(\xi) B w_h \text{ for some } w_h\in \mathbb{U}_h\big\}$.
\par
To establish the connection with the general mixed system~\eqref{eq:StateEq}, we set $r=\omega(\xi)(f-Bu_h)$ so that~\eqref{eq:lsq_form} is equivalent to:
% we get that the mixed system associated with~\eqref{eq:lsq_form} aims to find $r\in\mathbb V$ and $u_h\equiv S_h(\xi)\in\mathbb U_h$ such that
\begin{subequations}
\label{eq:mixed_LSQ}
\begin{empheq}[left=\left.,right=\right.,box=]{alignat=3} 
 & \biginnerprod{\varpi(\xi)\,r,v}_{L^2(\Omega)} + \big(Bu_{h},v\big)_{L^2(\Omega)} && = (f,v)_{L^2(\Omega)}, && \quad\forall v\in {\mathbb{V}},\\
 & \big(Bw_h,r\big)_{L^2(\Omega)} && = 0, && \quad \forall w_h\in \mathbb{U}_h\,.
\end{empheq}
\end{subequations}
Thus, in this case the bilinear forms $a(\xi;\cdot,\cdot)\in\mathcal L(\mathbb V\times\mathbb V;\mathbb R)$ and 
$b(\cdot,\cdot)\in\mathcal L(\mathbb{H}_B\times\mathbb V;\mathbb R)$ in~\eqref{eq:StateEq} are given by 
\begin{subequations}
\label{eq:bilinear_LSQ}
\begin{empheq}[left=\left.,right=\right.,box=]{alignat=3} 
  a(\xi;v_1,v_2):= & \, \biginnerprod{ \varpi(\xi)\,v_1,v_2}_{L^2(\Omega)},  && \quad\forall v_1,v_2\in\mathbb V = L^2(\Omega),\\
  b(w,v):= & \, (Bw,v)_{L^2(\Omega)}, && \quad\forall w\in\mathbb{H}_B, \forall v\in \mathbb V.
\end{empheq}
\end{subequations}

\begin{proposition}[Weighted least squares]
~\label{prop:LSQ_well}
Let $\varpi:L^2(\Omega)\to L^\infty(\Omega)$ be a differentiable map, such that 
for some positive constants $\varpi_{\min}$, $\varpi_{\max}$, $\varpi_\infty'$, and $\varpi_L$, the application $\varpi(\cdot)$ satisfies  
\begin{itemize}
    \item $\varpi_{\min}\leq \varpi(\xi)\leq \varpi_{\max}$, for all $\xi\in L^2(\Omega)$;
    \item $\|\varpi'(\xi)\|_{\mathcal L(L^2(\Omega),L^\infty(\Omega))}\leq \varpi_\infty'$, for all $\xi\in L^2(\Omega)$;
    \item $\|\varpi'(\xi_1)-\varpi'(\xi_2)\|_{\mathcal L(L^2(\Omega),L^\infty(\Omega))}\leq \varpi_L \|\xi_1-\xi_2\|_{L^2(\Omega)}$, for all $\xi_1,\xi_2\in\mathbb R$.
\end{itemize}
Then, the following statements hold true:
\begin{itemize}
    \item[(i)] The bilinear forms in~\eqref{eq:bilinear_LSQ} satisfy the $\inf$-$\sup$ conditions~\eqref{eq:inf-sup}, and thus the mixed problem~\eqref{eq:mixed_LSQ} is well-posed. 
    \item[(ii)] The state operator $S_h(\cdot)$ ($=u_h$) of the mixed problem~\eqref{eq:mixed_LSQ} is uniformly bounded on $\mathbb X = L^2(\Omega)$ and differentiable.
    \item[(iii)] The derivative $S_h'(\cdot)$ is uniformly bounded on $\mathbb X = L^2(\Omega)$ and Lipschitz continuous.
\end{itemize}
\end{proposition}
\begin{proof} 
See Appendix~\ref{sec:LSQ_well}
\end{proof}
\begin{remark}[Neural control of weighted least squares]
Proposition~\ref{prop:LSQ_well} guarantees that the conditions of Propositions~\ref{prop:stateEq} and~\ref{prop:state_deriv} are satisfied, hence Corollary~\ref{cor:constrainedProb} applies to the neural optimization of the above weighted least-squares formulation.
% .  all the hypothesis of Theorem~\ref{thm:reduced} will be satisfied for the weighted least-squares formulation~\eqref{eq:lsq_form}. See Corollary~\ref{cor:constrainedProb}. 
\end{remark}

\subsection{Weighted Galerkin formulations}
Consider a Hilbert space $\mathbb U=\mathbb{V}$ on~$\Omega\subset \mathbb{R}^d$ and a bilinear form $b\in\mathcal L(\mathbb V\times \mathbb V;\mathbb R)$ satisfying (for some constant $\beta>0$) the following conditions
\begin{subequations}
\label{eq:b-well}
\begin{empheq}[left=\left.,right=\right.,box=]{alignat=3} 
\label{eq:b-coercive}
& \sup_{v\in\mathbb V}\frac{b(w,v)}{\|v\|_\mathbb V}\geq \beta \|w\|_\mathbb V \,,\quad\forall w\in \mathbb V\,,\\%\notag\\
& \Big\{v\in\mathbb V: b(w,v)=0, \forall w\in\mathbb V\Big\}=\{0\}.
\end{empheq}
\end{subequations}
Given $f\in\mathbb V^*$, the well-known Babu\v{s}ka--Brezzi theory (see,~e.g.,~\cite{ErnGueBOOK2021b}) ensures the existence of an unique $u\in\mathbb V$ such that
\begin{equation}\label{eq:b-problem}
b(u,v)=f(v)\,,\quad \forall v \in \mathbb V\,.
\end{equation}

Now, given a weight function $\omega:L^2(\Omega)\to \mathbb{W}_+$ (the space $\mathbb{W}_+$ will be clarified later), a control $\xi\in \mathbb{X} = L^2(\Omega)$, and a conforming discrete subspace $\mathbb U_h\subset\mathbb V$, we consider the following weighted-Galerkin discretization of problem~\eqref{eq:b-problem}: 
\begin{empheq}[left=\left\{,right=\right.,box=]{alignat=3} 
\notag
&\text{Find } u_h\equiv S_h(\xi)\in \mathbb{U}_h :
\\
\label{eq:w-Galerkin}
&\qquad b\big(u_h,\omega(\xi){v_h}\big) = f\big(\omega(\xi){v_h}\big)\,,\quad\forall v_h\in\mathbb U_h\,.
\end{empheq}
Notice that one can directly identify the test space in~\eqref{eq:introPGxi} as $\mathbb{V}_h(\xi) = \omega(\xi)\mathbb{U}_h = \big\{ v \in \mathbb{V}\,\big|\, v = \omega(\xi)  w_h \text{ for some } w_h\in \mathbb{U}_h\big\}$. 
We will show next that problem~\eqref{eq:w-Galerkin} admits also an equivalent mixed formulation of the type~\eqref{eq:StateEq}, and therefore it fits the abstract setting of Section~\ref{sec:Abstract}. %into the class of problems considered within this work. 

First, we need to provide sense to the weighted object $\omega(\xi){v_h}\in\mathbb V$. Thus, we further consider an abstract Banach space $\mathbb{W}\equiv \mathbb{W}(\Omega)$ of measurable functions on $\Omega$, such that for any $\mathrm w\in\mathbb{W}$, the multiplication operator 
$M_\mathrm w:\mathbb V\to\mathbb V$ given by 
$$M_\mathrm w v:=\mathrm w v,\quad \forall v\in \mathbb V\,,$$ 
is a well-defined linear and continuous map.

\begin{example}[Multiplication in~$H^1$]
Let $\mathbb V=H^1(\Omega)$. Then it is easy to see that the Sobolev space $\mathbb{W}=W^{1,\infty}(\Omega)$ is a space of functions for which the multiplicative operator $M_\mathrm w:H^1(\Omega)\to H^1(\Omega)$ is a well-defined linear and continuous map, for all $\mathrm w\in W^{1,\infty}(\Omega)$. The latter is also true for Hilbert spaces $\mathbb V\subset L^2(\Omega)$ containing at most first-order (weak) derivatives in $L^2(\Omega)$ (e.g., first-order graph spaces).  
\end{example}

A particular subset of interest for us will be
$$
\mathbb{W}_{+}:=\Big\{\mathrm w\in\mathbb{W} \,\Big|\,  % \tfrac{1}{\mathrm{w}}\in\mathbb{W} \text{ and }
\exists \mathrm w_{\min}>0 \text{ for which }
 \mathrm w_{\min}
 \le \mathrm w(x)\le\tfrac{1}{\mathrm w_{\min}} , \forall x\in\Omega\Big\}\,.
$$
Notice that $\tfrac{1}{\mathrm{w}}\in \mathbb{W}_{+}$ iff $\mathrm{w}\in \mathbb{W}_{+}$.
% We will work under the assumption that $\mathrm w\in\mathbb{W}_+$ if and only if  $\frac{1}{\mathrm w}\in\mathbb{W}_+$. For such elements, we make use of the notation 
We can then define
$M^{-1}_\mathrm w:=M_{\frac{1}{\mathrm w}}$, which is justified by the fact that
\begin{equation}\label{eq:MM^-}
M^{-1}_\mathrm w(M_\mathrm w v) = v 
= M_\mathrm w(M^{-1}_\mathrm w v) \,,\quad\forall v\in\mathbb V.
\end{equation}
The adjoint operators of $M_\mathrm w$ and $M^{-1}_\mathrm w$ will be denoted by $M_\mathrm w^*$ and $M_\mathrm w^{-*}$ respectively.
Using the relations~\eqref{eq:MM^-} it is straightforward to see that the adjoint operators satisfy 
\begin{equation}\label{eq:M*M^-*}
M^{-*}_\mathrm w(M_\mathrm w^* \ell) = \ell
= M_\mathrm w^*(M^{-*}_\mathrm w \ell) \,,\quad\forall \ell\in\mathbb V^*.
\end{equation}

We translate problem~\eqref{eq:w-Galerkin} into operator notation by means of the operator $B\in \mathcal{L}(\mathbb V;{\mathbb V}^*)$ such that
$\mathbb V\ni w\mapsto Bw:=b(w,\cdot)\in\mathbb V^*$. Notice that such an operator is invertible thanks to conditions~\eqref{eq:b-well}.
Problem~\eqref{eq:w-Galerkin} translates into finding $u_h\equiv S_h(\xi)\in \mathbb V_h$ such that 
$$
\left<Bu_h,M_{\varpi(\xi)}^{-1}v_h\right>
=\left<f,M_{\varpi(\xi)}^{-1}v_h\right>\,,
\quad\forall v_h\in\mathbb V_h\,.
$$
Hence, by means of the adjoint relation we get 
\begin{equation}\label{eq:B^*r}
\left<M_{\varpi(\xi)}^{-*}(f-Bu_h),v_h\right>=0\,,
\quad\forall v_h\in\mathbb V_h\,.
\end{equation}
Since $B$ is invertible, so is $B^*:\mathbb V\to {\mathbb V}^*$ defined by 
$\mathbb V\ni v\mapsto b(\cdot,v)\in \mathbb V^*$. Therefore, there exists a unique $r\in\mathbb V$ such that $B^*r=M_{\varpi(\xi)}^{-*}(f-Bu_h)$ in ${\mathbb V}^*$.
Thus, multiplying this last equation by $M^*_{\varpi(\xi)}$, using~\eqref{eq:M*M^-*},~\eqref{eq:B^*r}, and the definition of $r\in\mathbb V$, we arrive to the mixed form
\begin{subequations}
\label{eq:Galerkin-mixed}
\begin{empheq}[left=\left\{,right=\right.,box=]{alignat=3} 
%\label{eq:StateEq_a}
 \left<B^*r,M_{\varpi(\xi)}v\right> & + b(u_{h},v) && = f(v), && \quad\forall v\in \mathbb{V},\\
%\label{eq:StateEq_b}
 b(v_h,r) & && = 0, && \quad \forall v_h\in{\mathbb{V}_h}\,.
\end{empheq}
\end{subequations}
Observe that~\eqref{eq:Galerkin-mixed} has the structure of~\eqref{eq:StateEq} for
$\hat{\mathbb V}:=\mathbb V=\mathbb U$; 
$\mathbb U_h:=\mathbb V_h$; and
$$
a(\xi;r,v):= \left<B^*r,M_{\varpi(\xi)}v\right>
    = b\big(\varpi(\xi)v,r\big).
$$
The next proposition establishes a sufficient condition for the well-posedness of~\eqref{eq:Galerkin-mixed}, or equivalently~\eqref{eq:w-Galerkin}.
\begin{proposition}[Weighted Galerkin]
\label{prop:WPG-well}
~\\
Let $b\in\mathcal L(\mathbb V\times\mathbb V;\mathbb R)$ be a bilinear form satisfying the ($\inf$-$\sup$) conditions~\eqref{eq:b-well}. Consider a conforming discrete subspace $\mathbb V_h\subset \mathbb V$ and let
\begin{alignat}{2}
\label{eq:wGalK}
K:=\Big\{v\in\mathbb V : b(v_h,v)=0\,,\forall v_h\in\mathbb V_h\Big\}.
\end{alignat}
Let $\varpi:L^2(\Omega)\to \mathbb{W}_+$ be a weight function such that 
\begin{equation}\label{eq:a-ker}
\big|b(\varpi(\xi)v,v)\big| \geq \alpha_h(\xi)\|v\|^2_\mathbb V\,,\quad
\forall v\in K\,,
\end{equation}
for some positive function $\alpha_h(\cdot)>0$. Then, the following statements hold true:
\begin{itemize}
\item[(i)] For any $f\in\mathbb V^*$ and $\xi\in L^2(\Omega)$, problems~\eqref{eq:w-Galerkin} and~\eqref{eq:Galerkin-mixed} are well-posed. 
\item[(ii)] If there exist uniform constants $\alpha>0$ and $\varpi_\infty>0$ such that $\alpha_h(\xi)\geq \alpha$ and $\|\varpi(\xi)\|_\mathbb{W}\leq \varpi_\infty$ for all $\xi\in L^2(\Omega)$, then the solution $u_h\equiv S_h(\cdot)$ to problems~\eqref{eq:w-Galerkin} and~\eqref{eq:Galerkin-mixed} is uniformly bounded
on $\mathbb X = L^2(\Omega)$.
\item[(iii)] Additionally, if $\varpi(\cdot)$ is differentiable, then $S_h(\cdot)$ is also differentiable. Moreover, if $\varpi'(\cdot)$ is uniformly bounded and Lipschitz-continuous, then also $S_h'(\cdot)$ is uniformly bounded and Lipschitz-continuous.
\end{itemize}
\end{proposition}
\begin{proof}
See Appendix~\ref{sec:WPG-well}.
\end{proof}
\begin{remark}[Neural control of weighted Galerkin]
Proposition~\ref{prop:WPG-well} guarantees that the conditions of
Propositions~\ref{prop:stateEq} and~\ref{prop:state_deriv} are satisfied, hence Corollary~\ref{cor:constrainedProb} applies to the neural optimization of the above weighted least-squares formulation.
\end{remark}
\begin{remark}[Inconvenient condition for weighted Galerkin]
\label{rem:wGalerkinCond}
While for the weighted least-squares method the conditions on the weight are explicit (recall Proposition~\ref{prop:LSQ_well}), for weighted Galerkin the condition~\eqref{eq:a-ker} is problem dependent. Furthermore, Example~\ref{ex:wGalLapl} shows it may require inconvenient constraints on~$\xi$. It seems therefore much more convenient to have neural control of least-squares formulations, or of dual minimal-residual formulations, as we will see in Section~\ref{ssec:wDDMRes}. 
%for the weighted Petrov-Galerkin formulation~\eqref{eq:w-Galerkin}; see also Corollary~\ref{cor:constrainedProb}. However, we emphasize that it is not a trivial task to find a weight fulfilling all the requirements of Proposition~\ref{prop:WPG-well} for a given weak form. For now, we leave this problem open. %
\end{remark}
\begin{example}[Weighted Galerkin for Laplacian]
\label{ex:wGalLapl}
%
%hypothesis of Theorem~\ref{thm:reduced} will be satisfied
%
Let us illustrate the difficulty of condition~\eqref{eq:a-ker} using the elementary Laplacian.
Let $\mathbb{V} = H_0^1(\Omega)$, $b(u,v) = \int_\Omega \grad u \cdot \grad v$ for all~$u,v\in \mathbb{V}$. Let $\varpi \in W^{1,\infty}(\Omega)$ such that $\varpi(x)\ge \mathrm{w}_{\min} > 0$ for all~$x\in \Omega$.
\par
In particular, let $\varpi(x) = \mathrm{w}_{\min} + c\,y \cdot (x-x_0)$ for some $c\in \mathbb{R}$ and $y,x_0\in \mathbb{R}^d$ such that $y \cdot (x-x_0)\ge 0$ for all~$x\in \Omega$. Then 
\begin{alignat}{2}
\label{eq:bwvv}
  b(\varpi\,v,v) &= \int_\Omega \Big( \varpi |\grad v| + v \grad \varpi \cdot \grad v\Big)
\\
\notag
  &= \mathrm{w}_{\min} \int_\Omega |\grad v|^2 + c\,y \cdot \int_\Omega \Big( (x-x_0) |\grad v|^2 + v\grad v\Big)
  %\big( \varpi |\grad v|^2 + v\, (a\cdot \grad v)\big)
\end{alignat}
Consider any $v\in \mathbb{V}$ so that
$
y \cdot \int_\Omega \Big( (x-x_0) |\grad v|^2 +
v\grad v\Big)
< 0\,.
$
Then there is a $c>0$ such that $b(\varpi\,v,v) = 0$. This shows that \eqref{eq:a-ker} can not be satisfied in general without additional conditions on~$\varpi$.
\par
Indeed, from~\eqref{eq:bwvv} a sufficient condition can be obtained. First notice that, for any $\varpi\in W^{1,\infty}(\Omega)$ such that $\varpi(x) \ge \mathrm{w}_{\min}$ for all $x\in \Omega$, 
\begin{alignat*}{2}
  b(\varpi\,v,v) 
  &\ge \mathrm{w}_{\min} \norm{\grad v}_{L^2(\Omega)}^2 - \norm{\grad \varpi}_{L^{\infty}(\Omega)} \norm{v}_{L^2(\Omega)} \norm{\grad v}_{L^2(\Omega)} 
\\ & \ge\big( \mathrm{w}_{\min} - C_{\Omega}\norm{\grad \varpi}_{L^2(\Omega)}  \big) \norm{\grad v}_{L^2(\Omega)}^2\,,
\end{alignat*}
where a Poincar\'e inequality was used. Therefore, the constraint $C_{\Omega}\norm{\grad \varpi}_{L^2(\Omega)} < \mathrm{w}_{\min}$ is sufficient to guarantee~\eqref{eq:a-ker}. Unfortunately, since~$\varpi = \varpi(\xi)$, such a condition translates into a constraint on~$\grad \xi$, which may be very inconvenient to impose in practice.
\end{example}

% \begin{example}[Laplacian]
% \begin{alignat*}{2}
%  \innerprod{\grad v,\varpi \grad v} +
%  \innerprod{\grad v,v \grad \varpi}
%  &\ge ... -  \norm{\grad v}\norm{v} \norm{\grad \varpi}_\infty
% \\
%  &\ge ... -  C|\Omega|^{-1}\norm{\grad v}^2 \norm{\grad \varpi}_\infty
% \end{alignat*}
% Hence $\norm{\grad \varpi}_{\infty} \le $ 
% \end{example}

% \begin{remark}[Weighted Petrov--Galerkin]
% \end{remark}

% \begin{example}[Lowest order PG for Strong advection]
%  \begin{alignat*}{2}
%  \sup_{v_h} \frac{\int_\Omega\beta\cdot \grad u_h\, \varpi v_h}{\norm{v_h}}
%  &=  \sup_{v_h} \frac{\sum_K\int_K\beta_K\cdot \grad u_h\, \varpi v_h}{\norm{v_h}}
% \\
%  &\ge \sup_{v_h} \frac{\sum_K C_{\varpi,K} \int_K\beta_K\cdot \grad u_h\, v_h}{\norm{v_h}} 
% \\
% &\ge C_{\varpi}
%  \end{alignat*}
% \end{example}

%
\subsection{Weighted discrete-dual minimal residual formulations}
\label{ssec:wDDMRes}
%
%Let $\mathbb U, \mathbb V$ be Hilbert spaces and $b\in\mathcal L(\mathbb U\times \mathbb V;\mathbb R)$. As usual, given $f\in \mathbb V^*$, the problem of finding $u\in\mathbb U$ such that 
%$b(u,v)=f(v)$ for all $v\in\mathbb V$, will represent the weak-form of a partial differential equation (PDE). We will work under the assumption that this PDE problem is well-posed.
%
%Let $\mathbb X$ be a Hilbert space of control parameters. 
%
Let $\mathbb U_h\subset \mathbb U$ and $\mathbb V_h\subset\mathbb V$ be discrete subspaces, and assume:
%of $\mathbb U$ and $\mathbb V$, respectively.
%We make two important assumptions on these discrete spaces:
%
\begin{subequations}
\label{eq:discrete_infsup}
\begin{empheq}[left=\left\{,right=\right.,box=]{alignat=3} 
%\label{eq:StateEq_a}
 & \dim(\mathbb V_h)>\dim(\mathbb U_h),\\
\label{eq:discrete_infsup_b}
 & \exists\, \beta_h>0\, :\, \inf_{w_h\in\mathbb U_h}\sup_{v_h\in\mathbb V_h}
 \frac{b(w_h,v_h)}{ \|w_h\|_\mathbb U\|v_h\|_\mathbb V}\geq\beta_h\,.
\end{empheq}
\end{subequations}
%
%\par
For each $\xi\in\mathbb X$, we consider an equivalent (weighted) inner product $(\cdot,\cdot)_{\mathbb V,\xi}$ on $\mathbb V$, i.e., such that its induced norm 
$$\mathbb V\ni v\mapsto \bignorm{v}_{\mathbb V,\xi}:=\sqrt{\biginnerprod{v,v}_{\mathbb V,\xi}} \quad\hbox{satisfies~\eqref{eq:norm_equiv}.}$$
%
%
% We propose the following discrete mixed formulation to compute an approximation $u_h\in\mathbb U_h$ of the original PDE:
The minimal-residual method that we consider is then: 
Given $\xi\in\mathbb X$, find $r_h\in\mathbb V_h$ and $u_h\equiv S_h(\xi)\in\mathbb U_h$ such that
\begin{subequations}
\label{eq:discrete_mixed}
\begin{empheq}[left=\left.,right=\right.,box=]{alignat=3} 
%\label{eq:StateEq_a}
 & \biginnerprod{r_h,v_h}_{\mathbb V,\xi} + b(u_h,v_h) && = f(v_h)\,, && \quad\forall v_h\in\mathbb V_h\,,\\
%\label{eq:StateEq_b}
 & b(w_h,r_h) && =0\,, && \quad\forall w_h\in\mathbb U_h\,.
\end{empheq}
\end{subequations}
This has the structure of~\eqref{eq:StateEq} for
$\hat{\mathbb V}:=\mathbb V_h$ and $a(\xi;r,v) := \innerprod{r,v}_{\mathbb{V},\xi}$.
\par
As shown in~\cite[Theorem~4.1]{MugZeeSINUM2020}, the mixed formulation~\eqref{eq:discrete_mixed} is equivalent to minimizing the residual as measured by a discrete-dual norm:
% the following residual minimization problem in the discrete-dual space $(\mathbb V_h)^*$ 
%
\begin{equation}\label{eq:ddminres}
u_h = \arg\!\!\min_{w_h\in\mathbb U_h}\left( \sup_{v_h\in \mathbb V_h}
\frac{|f(v_h)-b(w_h,v_h)|}{ \|v_h\|_{\mathbb V_h,\xi}}\right)\,.
\end{equation}
Because $\biginnerprod{\cdot,\cdot}_{\mathbb V,\xi}$ and $\norm{\cdot}_{\mathbb{V}_{h,\xi}}$ depend on~$\xi$, we refer to the above as a~\emph{weighted} discrete-dual minimal residual formulations.
\begin{proposition}\label{prop:ddminres}
Let the bilinear form $b(\cdot,\cdot)\in\mathcal L(\mathbb U\times\mathbb V;\mathbb R)$ and $(\mathbb U_h,\mathbb V_h)$ satisfy~\eqref{eq:discrete_infsup}.
Consider a parametrized set of equivalent inner-products 
$$
\big\{(\cdot,\cdot)_{\mathbb V,\xi}\in\mathcal L(\mathbb V\times\mathbb V;\mathbb R) : \xi\in\mathbb X\big\}\,,
$$ 
whose induced norms $\|\cdot\|_{\mathbb V,\xi}$ satisfy~\eqref{eq:norm_equiv} for some equivalence constants $C_{1,\xi}>0$ and $C_{2,\xi}>0$. Let $A:\mathbb X\to \mathcal L(\mathbb V,\mathbb V^*)$ be defined by $A(\xi)v:=(v,\cdot)_{\mathbb V,\xi}\in\mathbb V^*$, for all $\xi\in\mathbb X$ and 
$v\in\mathbb V$. Then, the following statements hold true:
\begin{itemize}
    \item[(i)] The mixed discrete formulation~\eqref{eq:discrete_mixed} is well-posed.
    \item[(ii)] If there exist uniform constants $\tilde C_1>0$ and $\tilde C_2>0$ such that $C_{1,\xi}\geq \tilde C_1$ and $C_{2,\xi}\leq \tilde C_2$ for all $\xi\in\mathbb X$, then the solution $u_h\equiv S_h(\cdot)$ to problems~\eqref{eq:discrete_mixed} and~\eqref{eq:ddminres} is uniformly bounded on~$\mathbb X$.
    \item[(iii)] Additionally, if $A(\cdot)$ is differentiable, then $S_h(\cdot)$ is also differentiable. Moreover, if $A'(\cdot)$ is uniformly bounded and Lipschitz-continuous, then also $S_h'(\cdot)$ is uniformly bounded and Lipschitz continuous.
\end{itemize}
\end{proposition}
\begin{proof}
See Appendix~\ref{sec:ddminres_proof}.
\end{proof}
\begin{remark}[Neural control of weighted residual minimization]
Proposition~\ref{prop:ddminres} guarantees that the conditions of
Propositions~\ref{prop:stateEq} and~\ref{prop:state_deriv} are satisfied, hence Corollary~\ref{cor:constrainedProb} applies to the neural optimization of the above weighted minimal-residual formulation.
%
% that all the hypothesis of Theorem 2.A will be satisfied
% for the weighted discrete-dual minimal-residual formulation~\eqref{eq:discrete_mixed}; see also Corollary 2.10. 
%For a particular situation of a family of (weighted) inner product satisfying the hypothesis of Proposition~\ref{prop:ddminres}, we refer to the next Example~\ref{ex:weighted-H1}.  
\end{remark}
\begin{example}[Weighted $H^1(\Omega)$ inner-product] 
\label{ex:weighted-H1}
Consider a differentiable weight function $\omega:L^2(\Omega)\to L^\infty(\Omega)$, such that for some given constants $\omega_{\max}>\omega_{\min}>0$, and for all $\xi\in L^2(\Omega)$, we have $\omega_{\min}\leq \omega(\xi)\leq \omega_{\max}$. We further assume that $\omega'(\cdot)$ is uniformly bounded and Lipschitz-continuous.

Given $\xi\in L^2(\Omega)$, define the following weighted $H^1(\Omega)$ inner-product:
$$
(v_1,v_2)_{H^1,\xi}:=\int_\Omega \omega(\xi)\nabla v_1\cdot \nabla v_2 +
\int_\Omega v_1 v_2\,.
$$
Observe that 
$$
\min\{1,\omega_{\min}\}\|v\|_{H^1}^2\leq (v,v)_{H^1,\xi}\leq \max\{1,\omega_{\max}\}\|v\|_{H^1}^2\,,\quad\forall v\in H^1(\Omega)\,.
$$
Hence, statement (ii) of Proposition~\ref{prop:ddminres} is satisfied with 
$\tilde C_1=\sqrt{\min\{1,\omega_{\min}\}}$ and $\tilde C_2=\sqrt{\max\{1,\omega_{\max}\}}$. 

On the other hand, given $\xi\in L^2(\Omega)$, the operator $A(\xi)$ is defined by the following action: 
$$A(\xi)v=\big(\omega(\xi)\nabla v,\nabla(\cdot)\big)_{L^2(\Omega)}+(v,\cdot)_{L^2(\Omega)}\,,
\quad \forall v\in H^1(\Omega)\,.$$
Therefore, is easy to see that $A(\cdot)$ satisfies the statement (iii) of Proposition~\ref{prop:ddminres}. Indeed, observe that $A(\cdot)$ is differentiable and $[A'(\xi)\eta]v=\big([\omega'(\xi)\eta]\nabla v,\nabla(\cdot)\big)_{L^2(\Omega)}$ for any direction $\eta\in L^2(\Omega)$. Moreover, $A'(\cdot)$ is uniformly bounded and Lipschitz-continuous, since $\omega'(\cdot)$ is uniformly bounded and Lipschitz continuous. 

Of course, for any $v_1,v_2\in H^1(\Omega)$, we may have chosen the following equivalent inner-products where we can prove similar results:
$$
\begin{array}{rl}
(v_1,v_2)_{H^1,\xi} := & \big(\nabla v_1,\nabla v_2\big)_{L^2(\Omega)} + \big(\omega(\xi)v_1,v_2\big)_{L^2(\Omega)}\,,\\\\
(v_1,v_2)_{H^1,\xi}:= & \big(\omega(\xi)\nabla v_1,\nabla v_2\big)_{L^2(\Omega)} + \big(\omega(\xi)v_1,v_2\big)_{L^2(\Omega)}\,.
\end{array}
$$
Also, for $H_0^1(\Omega)$, we could consider just $\biginnerprod{\omega(\xi)\grad v_1\cdot\grad v_2}_{L^2(\Omega)}$.
\end{example}

\section{Numerical results}
\label{sec:Numerics}
In this section, we consider numerical examples for the advection--reaction PDE in 1-D and 2-D. We consider both weighted least squares and weighted residual minimization.%
\footnote{Weighted Galerkin is not considered in view of Remark~\ref{rem:wGalerkinCond}.}
\par
We construct weight functions~$\omega:L^2(\Omega)\rightarrow L^\infty(\Omega)$ that are based on algebraic expressions, i.e., for which $\omega(\xi)(x) = \omega(\xi(x))$ for~$x\in\Omega$. These are convenient expressions, but the price to pay is that $\omega':L^2(\Omega)\rightarrow \mathcal{L}(L^2(\Omega);L^\infty(\Omega))$ can not be Lipschitz. We do not believe this to have a major impact, and we leave the construction of more complicated weight functions for future investigation. While using algebraic weight functions, we have not observed any undesirable numerical effects. In fact, our results in Section~\ref{ssec:numExp:convANN} do demonstrate quasi-optimal convergence, as expected in our current theory.
%\par
%
%and we believe that this shows a limitation of our current theory. 
% numerical results have not shown any undesirable do however corroborate with the We have not noticed
%
%We are able to satisfy all  assumptions 
\subsection{Quantities of interest (point values)}
\label{ssec:numExp:QoI}

\subsubsection{Weighted least-squares approach}
\label{sec:point_value}
Let $\Omega=(0,1)\subset \mathbb R$ and $r>0$. Consider the advection-reaction problem
% Purpose: Create coarse approximation that has zero overshoot (small error for known quantities) for reaction-advection:
\begin{equation}\label{eq:ode2}
\left\{
\begin{array}{rl}
  u' + r\,u & = r \quad \text{in } \Omega\,,
  \\
  u(0) & = 0\,.
\end{array}\right.
\end{equation}
Since the exact solution to~\eqref{eq:ode2} is $u(x) = 1-\exp(-rx)$, we observe that $u(x)\to 1$ when $r\to +\infty$, for all $x>0$. Hence, for $r>0$ sufficiently large, the exact solution has a boundary layer in the neighbourhood of $x=0$.

Let $\mathbb{U}_h\subset H_{(0}^1(\Omega):=\{w\in H^1(\Omega): w(0)=0\}$ be the conforming subspace of continuous piecewise linear functions on the uniform mesh of $N$ elements of size $h=1/N$. We use the weighted least squares method from~\eqref{eq:lsq_form},
%
%We define the weighted least-squares formulation of~\eqref{eq:ode2} to be
%
% \begin{equation}\label{eq:lsq2}
% \left\{
% \begin{array}{l}
%   \text{Find } u_h 
%   \equiv S_h(\xi)\in \mathbb{U}_h:
%   \\\displaystyle
%   \int_0^1  \omega(\xi)\left({r(1-u_h)-u_h'}\right) \, w_h' = 0\,,
% \quad \forall w_h\in \mathbb{U}_h\,,
% \end{array}\right.
% \end{equation}
%
with weight function:
%is such that
%
\begin{alignat}{2}\label{eq:weight_lsq2}
 \omega\big(\xi(x)\big) := 1+ \frac{M}{1+\exp(-\xi(x))}\,, \qquad M>0.
\end{alignat}

% Take small~$N$, e.g., $N=8$, take proportionally large~$r$, e.g., $r = 10N = 80$, 

% Take~$\varpi(\xi)$ as previous numerical example\textbf{?} \textbf{(Will this work?)}
% \\
It is well known that the standard least-squares solution (i.e., the one with $\omega(\xi)\equiv 1$) will exhibit overshoots around the boundary layer. Aiming to remedy this situation, we choose a cost functional that measures the distance to the exact solution at the point value $x=h$. In fact, we consider  
\begin{alignat*}{2}
j(\xi) := \frac{1}{2} \Big( u(h) - u_{h,\xi}(h) \Big)^2 + \frac{\alpha}{2} 
\|\xi\|^2_{L^2}\,,\qquad \alpha\geq 0.
\end{alignat*}
% with $z = 1/N$ and $\bar{q} = u(z)$. 
% \par
% Take some~$n$ so that the overshoot is significantly reduced, when~$\alpha = 0$.

Let~$\mathcal{M}_8$ be the set of neural network functions with one hidden layer, $8$-neurons, and $\ReLU$ activation, i.e., 
\begin{alignat}{2}\label{eq:M8}
\mathcal{M}_8 := \bigg\{
%\xi:\mathbb{R} \rightarrow \mathbb{R} 
%\,\Big|\,
\eta_8(x) = \sum_{j=1}^8 c_j \ReLU(W_j x +b_j) 
\,\Big|\, c_j,W_j,b_j\in\mathbb R
\bigg\}\,.
\end{alignat}
We then consider the neural optimization of~$j(\cdot)$; see Definition~\eqref{def:qminprob}.
% Given $\delta>0$, our neural control problem will be
% \begin{equation}\label{eq:j_lsq2}
% \left\{
% \begin{array}{ll}
%  \text{Find } \xi_8 \in \mathcal{M}_8:
%  \\
%   j(\xi_8) \le \displaystyle\inf_{\eta_8\in \mathcal{M}_8} j(\eta_8) + \delta\,.
% \end{array}
% \right.
% \end{equation}
%
\par
\begin{figure}[!t]
\includegraphics[width=0.52\linewidth]{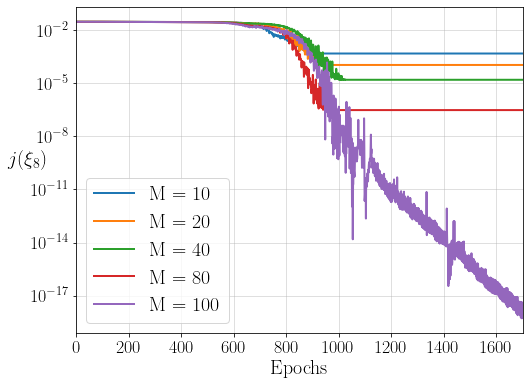}
\includegraphics[width=0.49\linewidth]{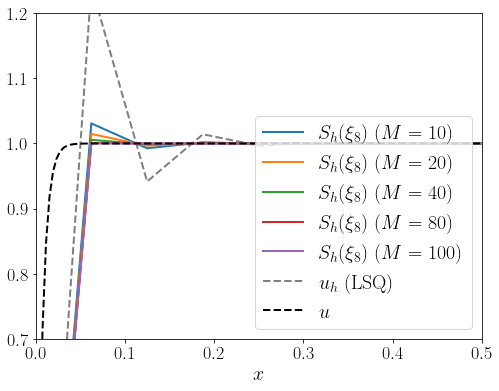}
\caption{Point value control for weighted least-squares. Minimization of the cost functional for several values of $M$ (left). Overshoot control of the discrete solutions (right).}
\label{fig:M-analysis1}
\end{figure}

For our first experiment, we choose a finite element space $\mathbb U_h$ consisting of $N=16$ elements of size $h=1/16$. We set $r=160$ and $\alpha=0$. We compute least-squares approximations for several configurations of the weight function~\eqref{eq:weight_lsq2}, varying the $M$ constant.
Figure~\ref{fig:M-analysis1} (left) shows that the weight needs to have enough room for variability ($M=100$) in order to pull down the cost functional to zero. Figure~\ref{fig:M-analysis1} (right) shows that our strategy is effective in reducing the overshoots of the finite element solution.

\begin{figure}[!t]
\includegraphics[width=0.52\linewidth]{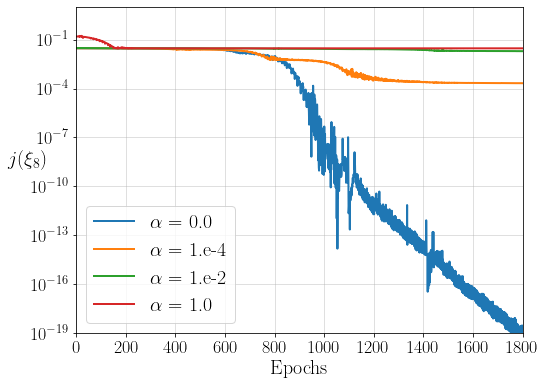}
\includegraphics[width=0.48\linewidth]{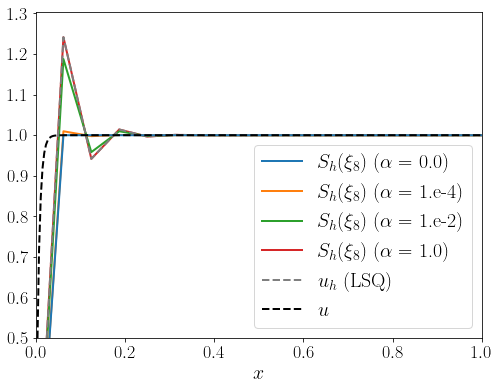}
\caption{Point value control for weighted least-squares. Minimization of the cost functional for several values of $\alpha$ (left). Overshoot control of the discrete solutions (right).}
\label{fig:alpha-analysis1}
\end{figure}

For the second experiment of this section, we fix $M=100$ and we investigate variations of the $\alpha$-parameter. Figure~\ref{fig:alpha-analysis1} (left) suggest that the $L^2$-norm of $\xi$ has to be able to reach high values (case when $\alpha=0$) in order to pull down to zero the cost functional. This is also related to allowing the weight to have more variability. Figure~\ref{fig:alpha-analysis1} (right) shows the impact of $\alpha$ reducing the overshoots of the finite element solution (the smaller $\alpha$, the better).
%

% Create figure with:
% \begin{itemize}
% \item Exact solution.
%     \item LSQ solution~$u_h$ with constant weight (hence no training needed) (equivalent to $\alpha =\infty$).
%     \item Trained approximation~$u_h = S_h(\xi_n)$, with~$\alpha = 0$
%     \item Trained approximations, with say 3 increasing values of~$\alpha$. Unsure what values of~$\alpha$ to pick, but something so that the overshoot varies between the other two approximations.
% \end{itemize}

\subsubsection{Weighted discrete-dual residual minimization approach}
This experiment has exactly the same configuration of the previous experiment in Section~\ref{sec:point_value}, except that $S_h(\xi)$ is computed with the discrete-dual minimal residual methodology. First, the approximation (\emph{trial}) space $\mathbb U_h\subset L^2(\Omega)$ corresponds to the space of piecewise constants functions over the mesh. Additionally, we make use of a discrete \emph{test} space
$\mathbb V_h\subset H^1_{0)}(\Omega):=\{v\in H^1(\Omega)\,:\, v(1)=0\}$ consisting in conforming piecewise linear functions over the refined uniform mesh of $2N=32$ elements. The weighted discrete-dual residual minimization formulation that computes $S_h(\xi)$ is as follows: $\text{Find } r_h\in\mathbb V_h \text{ and } u_h \equiv S_h(\xi)\in \mathbb{U}_h$ such that
\begin{equation}\label{eq:w-ddminres_advreac}
\left\{
\begin{array}{lll}
 \displaystyle
  \int_0^1\omega(\xi) r_h'v_h' - \int_0^1u_h(v_h'-r\,v_h) & = r\displaystyle\int_0^1 v_h\,,
  & \forall v_h\in \mathbb{V}_h\,,\\
  -\displaystyle \int_0^1w_h(r_h'-r\,r_h) & = 0\,,
& \forall w_h\in \mathbb{U}_h\,.
\end{array}\right.
\end{equation}

\begin{figure}[!t]
\includegraphics[width=0.51\linewidth]{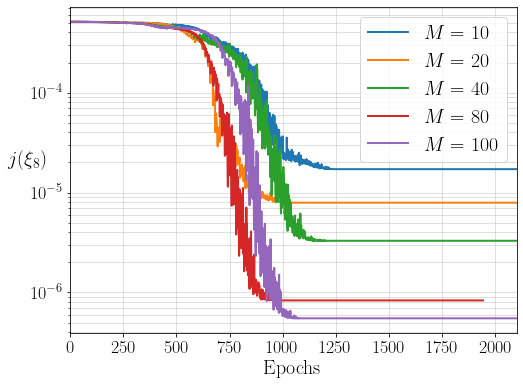}
\includegraphics[width=0.5\linewidth]{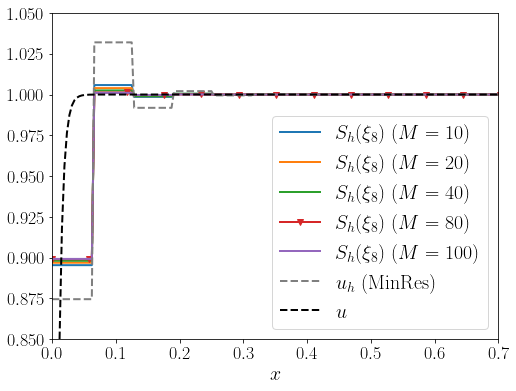}
\caption{Point value control for weighted discrete-dual residual minimization. Optimization of the cost functional for several values of $M$ (left). Overshoot control of the discrete solutions (right).}
\label{fig:M-analysis_ddminres}
\end{figure}

As in the previous Section~\ref{sec:point_value}, the computation of $S_h$ is carried out for several configurations of the weight function $\omega(\xi)$ (see~\eqref{eq:weight_lsq2}), varying its $M$ constant.
Figure~\ref{fig:M-analysis_ddminres} (left) shows that larger values of $M$ allow to pull down faster the cost functional in the training procedure. Figure~\ref{fig:M-analysis_ddminres} (right) shows how the overshoots of the finite element solutions are controlled.

\begin{figure}[!t]
\includegraphics[width=0.51\linewidth]{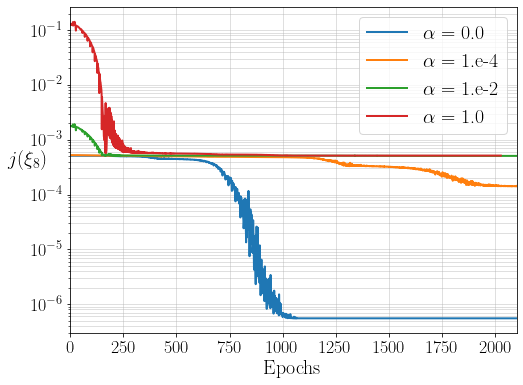}
\includegraphics[width=0.50\linewidth]{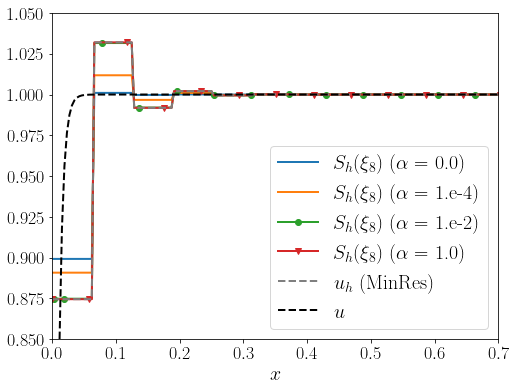}
\caption{Point value control for weighted discrete-dual residual minimization. Optimization of the cost functional for several values of $\alpha$ (left). Overshoot control of the discrete solutions (right).}
\label{fig:alpha-analysis_ddminres}
\end{figure}

The second experiment investigates variations of the $\alpha$-parameter. Figure~\ref{fig:alpha-analysis_ddminres} (left) suggest that the smaller $\alpha$, the better for faster minimization of $j(\cdot)$.   Figure~\ref{fig:alpha-analysis_ddminres} (right) shows the impact of $\alpha$ reducing the overshoots of the finite element solution. 

\subsection{Convergence of artificial neural networks}
\label{ssec:numExp:convANN}
%\subsubsection{Abundance of data (on the flux)}
\label{sec:sanity1}
Let $\Omega:=(0,1)\subset\mathbb R$ be a one-dimensional domain and consider the simple advection problem
\begin{equation}\label{eq:ode1}
\left\{
\begin{array}{rl}
 u' = & f 
 \quad \text{in } \Omega\,,
 \\
 u(0) = & 0\,,
 \end{array}
 \right.
\end{equation}
with $f(x) := \pi\sin (\pi x)$. Notice the exact solution to~\eqref{eq:ode1} is $u(x) = 1 - \cos (\pi x)$.
\par
Let $H_{(0}^1(\Omega):=\{w\in H^1(\Omega): w(0)=0\}$ and let $\mathbb{U}_h\subset H_{(0}^1(\Omega)$ be the finite element subspace of continuous piecewise linear functions on a uniform mesh consisting of $N$ elements of size $h=1/N$. We consider the weighted least-squares formulation
% of~\eqref{eq:ode1}
%
\begin{equation}\label{eq:lsq_ode1}
\left\{
\begin{array}{l}
  \text{Find } u_h 
  \equiv S_h(\xi)\in \mathbb{U}_h:
  \\\displaystyle
  \int_0^1  \omega(\xi)\left({f-u_h'}\right) \, w_h' = 0\,,
\quad \forall w_h\in \mathbb{U}_h\,,
\end{array}\right.
\end{equation}
where the weight function is such that
\begin{alignat}{2}\label{eq:omega_sanity1}
 \omega\big(\xi(x)\big) := \frac{1}{2}+ \frac{2}{1+\exp(-\xi(x))}\,.
\end{alignat}

% Equivalently, with~$\mathbb{V}= L^2(0,1)$,
% \begin{alignat*}{2}
%   \text{Find } (r, u_h)\in 
%   \mathbb{V}\times\mathbb{U}_h:
%   \qquad &
%   \\
%   \int_0^1 \varpi(\xi_n)^{-1} r\,v\,
%   + \int_0^1 u_h' \, v & = \int_0^1 f\, v
% \qquad && \forall v\in \mathbb{V}
% \\
%   \int_0^1 w_h' \, r  &= 0
% \qquad && \forall w_h\in \mathbb{U}_h
% \end{alignat*}

Let~$\mathcal{M}_n$ be the set of neural network functions with one hidden layer, $n$-neurons, and $\ReLU$ activation, i.e., 
\begin{alignat*}{2}
\mathcal{M}_n := \bigg\{
%\xi:\mathbb{R} \rightarrow \mathbb{R} 
%\,\Big|\,
\eta_n(x) = \sum_{j=1}^n c_j \ReLU(W_j x +b_j) 
\,\Big|\, c_j,W_j,b_j\in\mathbb R
\bigg\}\,.
\end{alignat*}
Consider the cost functional
%(note $S_h(\xi_n)' = u_h'$ is constant in each element):
\begin{alignat}{2}\label{eq:cost1}
  j(\xi) := 
  %\equiv j_N(\xi_n) :=  
  \frac{1}{2}  \int_{0}^1 
\bar\omega(x)
  \Big(f(x) - u_{h,\xi}'(x)\Big)^2\,\dd x\,,
  % j(\xi_n) \equiv j_N(\xi_n) := \sum_{i=1}^N \frac{1}{2} \bar{\varpi}(\tfrac{x_{i-1}+x_i}{2}) \int_{x_{i-1}}^{x_i} \Big(f(x) - S_h(\xi_n)'(x)\Big)^2\,\dd x
\end{alignat}
with $\bar\omega(x) = 1 + \sin(\pi x / 2)$.
% Thus, given $\delta_n>0$, our neural control problem will be
% $$
% \left\{
% \begin{array}{ll}
%  \text{Find } \xi_n \in \mathcal{M}_n:
%  \\
%   j(\xi_n) \le \displaystyle\inf_{\eta_n\in \mathcal{M}_n} j(\eta_n) + \delta_n\,.
% \end{array}
% \right.
% $$
%

\begin{figure}[!t]
\includegraphics[width=0.5\linewidth]{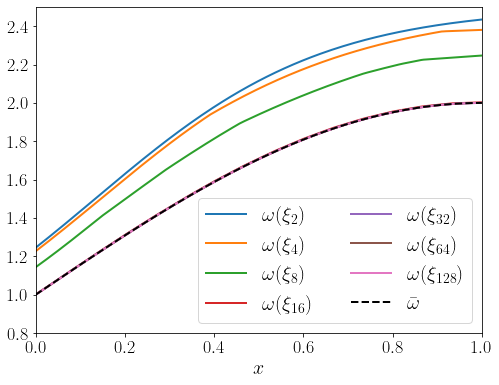}
\includegraphics[width=0.5\linewidth]{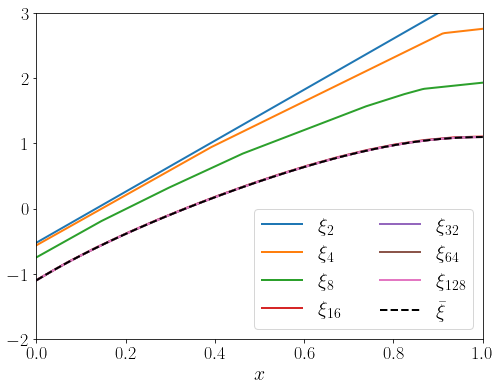}
\caption{Convergence of $\omega(\xi_n)\to\bar{\omega}$ (left) and $\xi_n\to\bar{\xi}$
(right), as $n\to+\infty$.}
\label{fig:test1}
\end{figure}

Since the minimization of the cost functional and the discrete problem~\eqref{eq:lsq_ode1} are both weighted least-squares formulations of the same problem~\eqref{eq:ode1}, we expect that $\omega(\xi_n)\to \bar{\omega}$ as $n\to+\infty$, which is confirmed in~Figure~\ref{fig:test1} (left). Additionally, solving for $\xi_n$ we get (see~Figure~\ref{fig:test1} (right))
$$
\xi_n(x) \longrightarrow \bar{\xi}(x)=-\ln \bigg(\frac{2}{\sin(\pi x/2)+1/2}-1\bigg),\qquad\hbox{as } n\to+\infty.
$$
To initialize the minimization algorithm, we have chosen 
$\xi_n^{(0)}\in\mathcal M_n$ as the neural network function that (linearly) interpolates $\bar{\xi}$ on a uniform mesh of $n-1$ subintervals of $\Omega$ (i.e., having $n$ uniformly distributed nodal points). The space $\mathbb U_h$ has been fixed to $N=16$ uniform elements. 

% On another experiment, we fix the number of neurons to be $n=?$ and we enrich the space $\mathbb U_h$ to see observe the quality 
% of the solution $u_h=S_h(\xi_n)$ as $1/h=N \to +\infty$. Results are depicted in Fig.~\ref{fig:test2} for the error $\|u-u_h\|_{L^2}$. Notice that we asymptotically reach the optimal convergence rate.

\begin{figure}[!t]
\begin{center}
\includegraphics[width=0.5\linewidth]{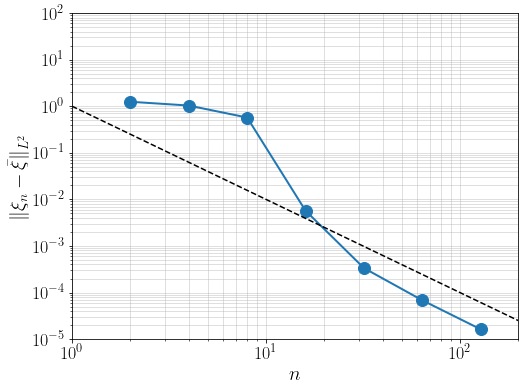}
\caption{$L^2$ error of $\xi_n$ as $n \to +\infty$.}
\label{fig:test2}
\end{center}
\end{figure}

In Figure~\ref{fig:test2}, we plot the error $\|\bar\xi-\xi_n\|_{L^2}$, which confirms quasi-optimal convergence behaviour; indeed the asymptotic rate is $O(n^{-1/2})$, which is expected for our single-hidden-layer ReLU neural network approximations (continuous piecewise-linear polynomials).

\subsection{$L^1$-based controls}
\label{ssec:numExp:L1}
We now consider numerical experiments that incorporate a stabilization mechanism. We note that the employed cost functionals use an $L^1$-type norm, and hence do not fit within the currently presented theory. However our numerics show that desirable quasi-minimizers have been computed. 

\subsubsection{Minimizing the total variation}
\begin{figure}[!t]
\includegraphics[width=0.52\linewidth]{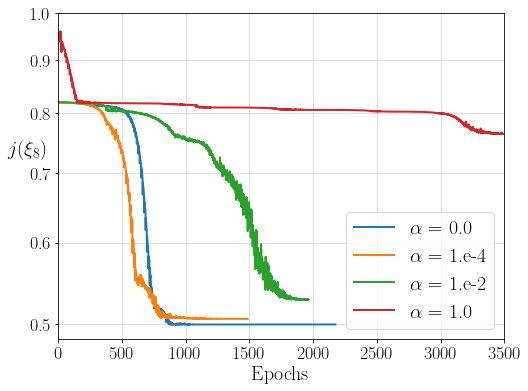}
\includegraphics[width=0.49\linewidth]{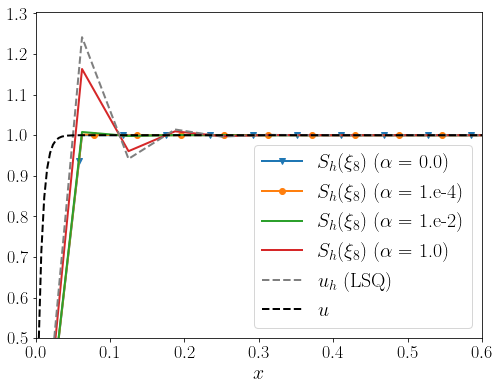}
\caption{Total variation control. Minimization of the cost functional for several values of $\alpha$ (left). Overshoot control of the discrete solutions (right).}
\label{fig:tv_alpha}
\end{figure}

In this section we work exactly with the same problem of the previous Section~\ref{sec:point_value}, but we introduce a modification in the cost functional. Instead of minimizing the distance to the exact solution of a particular point value (supervised training), we take an unsupervised approach by minimizing the total variation of $u_h$ (i.e., the $L^1$-norm of $u_h'$). Hence, we consider the cost functional:
\begin{alignat*}{2}
j(\xi) :=\big\|u_{h,\xi}'\big\|_{L^1} + \frac{\alpha}{2}\|\xi\|_{L^2}^2\,,\qquad \alpha\geq 0.
\end{alignat*}
For a fixed value of $M=100$, Figure~\ref{fig:tv_alpha} (left) shows the behavior of the cost functional for different values of $\alpha$, indicating that this value has to be chosen small enough to speed up the minimization process.
Figure~\ref{fig:tv_alpha} (right) shows the quality of overshoot reduction for several values of $\alpha$. 
%

% Purpose: Create coarse approximation that has zero overshoot (by having small $L^1$ residual).
% \\
% Same as previous numerical experiment, except:
% \begin{alignat*}{2}
% & j(\xi_n) := \int_0^1 \Big| \operatorname{res}(x) \Big|\,\dd x
% + \alpha \int_0^1 \xi_n(x)^2 \, \dd x
% \\
% & \text{with }
% \operatorname{res}(x) := f(x) - S_h(\xi_n)'(x) - r S_h(\xi_n)(x)
% \end{alignat*}

\subsubsection{Minimizing the $L^1$ residual (1D domain)}
\label{sec:viscosity_1D}
This experiment is inspired by the example of Guermond~\cite[Section~4.6.2]{GueSINUM2004}. 
As usual $\Omega=(0,1)\subset\mathbb R$.
The idea is to interpret the following overconstrained problem:
\begin{equation}\label{eq:ode3}
\left\{
\begin{array}{rl}
  u' + u & = 1 \quad \text{in } \Omega\,,
  \\
  u(0)=u(1) & = 0\,,
\end{array}\right.
\end{equation}
as the limiting case of a \emph{vanishing viscosity regime} (i.e., an equivalent problem having an extra $-\varepsilon u''$ term that vanishes as $\varepsilon\to 0^+$). Of course, the exact solution that we want to approach ($u(x)=1-e^{-x}$) only satisfies one of the boundary conditions. However, any discrete solution in a $H_0^1(\Omega)$-conforming space must satisfy both constrains.  
In this case, it is well-known that the standard least-squares solution to this problem does not deliver satisfactory results. To remedy this drawback, we propose a cost functional that mimics the $L^1$ residual minimization as proposed in~\cite{GueSINUM2004}.
Thus, our (\emph{unsupervised}) cost functional will be
$$
j(\xi):=\big\|\underbrace{1-u_{h,\xi}-u_{h,\xi}'}_{\hbox{\tiny residual}}\big\|_{L^1} + \frac{\alpha}{2}\|\xi\|^2_{L^2}\,,\qquad\alpha\geq 0\,.
$$

\begin{figure}[!t]
\begin{center}
\includegraphics[width=0.5\linewidth]{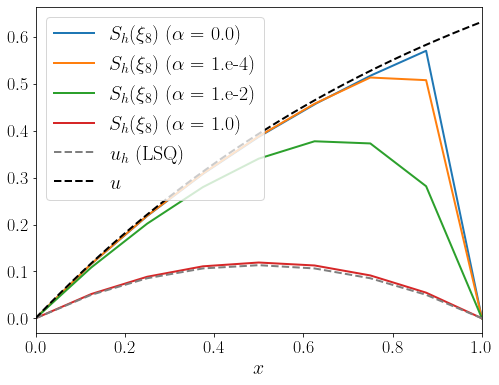}
\caption{Discrete weighted least-squares solutions, with $L^1$ residual minimization control, for several values of the $\alpha$-parameter.}
\label{fig:viscosity}
\end{center}
\end{figure}

We consider the weighted least-squares formulation for $u_{h,\xi}$, solved on a uniform mesh of $N=8$ elements.
%(i.e., formulation~\eqref{eq:lsq2} with $r=1$ and $H_0^1(\Omega)$-conforming $\mathbb U_h$). 
For a fixed $M=1000$ constant in the weight function~\eqref{eq:weight_lsq2}, we compute the discrete solution for several values of the $\alpha$-parameter. 
Large values of $\alpha$ allow for small values of $\|\xi\|_{L^2}$, and thus 
the weight becomes almost constant (close to the standard least-squares approach).
On the other hand, small values of $\alpha$ allow for more variability of the weight, and thus, we observe that we can recover a discrete solution mimicking the vanishing viscosity case (see~Fig.~\ref{fig:viscosity}). 

% Purpose: Create approximation (``unsupervised'') that mimic $L^1$ residual minimization (makes LSQ consistent and convergent to the viscosity solution).

% Following Guermond's example~\cite[Section~4.6.2]{GueSINUM2004}, but in 1-D.
% Same as previous numerical example, but two(!) Dirichlet BCs (to mimic the limiting diffusion situation): 
% \begin{alignat*}{2}
%  & u' + r\,u =f \qquad \text{in } (0,1)\,,
%   \\
%  & u(0) = 0\,,
%  \\
%  & u(1) = 0
% \end{alignat*}
% with~$r=1$ and $f(x) = 1$.
% %
% \\
% Create figure with:
% \begin{itemize}
% \item Exact solution: $u(x) =  1-\exp(-x)$.
%     \item LSQ solution~$u_h$ with constant weight (hence no training needed): It is known that LSQ does not converge to the right solution!
% %
%     \item Trained approximation~$u_h = S_h(\xi_n)$, with~$\alpha = 0$
%     \item Trained approximations, with say 4 increasing values of~$\alpha$. Unsure what values of~$\alpha$ to pick, but something so that it is between the other two approximations.
% \end{itemize}

\subsubsection{Minimizing $L^1$ residual (2D domain)}
\begin{figure}[!t]
\includegraphics[width=0.331\linewidth]{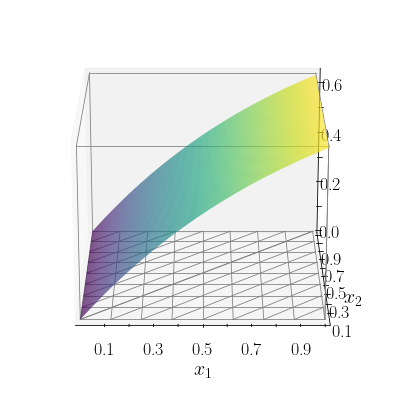}
\includegraphics[width=0.331\linewidth]{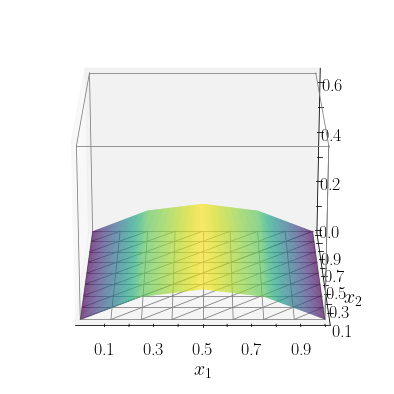}
\includegraphics[width=0.331\linewidth]{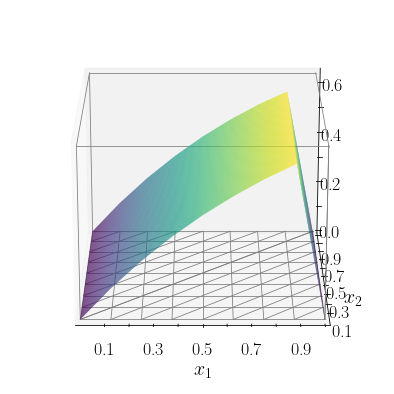}
\caption{Overconstrained weighted least-squares for advection-reaction, with $L^1$ residual minimization control. From left to right: exact solution; standard overconstrained least-squares, controlled weighted least-squares.}
\label{fig:viscosity2D}
\end{figure}

This is the two-dimensional version of the previous example in Section~\ref{sec:viscosity_1D}. Let $\Omega=(0,1)^2\subset\mathbb R^2$. For an advection field $\vec\beta=(1,0)$, we consider the over-constrained problem:
\begin{equation}\label{eq:viscosity_2D}
\left\{
\begin{array}{rlcl}
 \vec\beta\cdot\nabla u + u & =1 & \text{in } & \Omega\,,\\
 u & = 0 & \text{on } & \{(x_1,x_2)\in\partial\Omega: x_1=0 \text{ or } x_1 = 1\}\,.
% \partial_y u & =0 & \hbox{over} & \{y=0\}\cup\{y=1\}\,.
\end{array}
\right.
\end{equation}
We approach~\eqref{eq:viscosity_2D} using a coarse (and over-constrained) finite element space of piecewise linears functions of the form
$$
\mathbb U_h\subset \{w\in H_0^1(\Omega) : w(0,x_2)=w(1,x_2)= 0,\,\forall x_2\in[0,1]\}.
$$
We use the weighted least-squares method:
\begin{equation}\label{eq:lsq_viscosity_2D}
\left\{
\begin{array}{l}
  \text{Find } u_h 
  \equiv S_h(\xi)\in \mathbb{U}_h:
  \\\displaystyle
  \int_\Omega\omega(\xi)\big({1-u_h-\beta\cdot\nabla u_h}\big)\big( \beta\cdot\grad w_h + w_h \big)= 0\,,
\quad \forall w_h\in \mathbb{U}_h\,,
\end{array}\right.
\end{equation}
using the same weight~\eqref{eq:weight_lsq2} with $M=1000$. On the other hand, the cost functional $j(\cdot)$ for this case is defined as
$$
j(\xi):=\big\|1-u_{h,\xi}-\beta\cdot\nabla u_{h,\xi}\big\|_{L^1} + \frac{\alpha}{2}\|\xi\|_{L^2}\,,\quad \alpha\geq 0\,.
$$
The discrete neural network space where we minimize $j(\cdot)$ will be $\mathcal M_8$ (see~\eqref{eq:M8}). Results for the $\alpha=0$ case are depicted in Figure~\ref{fig:viscosity2D}. We observe a strong correlation with the results in~\cite[Figure~9]{GueSINUM2004}. 

%

%
% with~$r=1$ and $f(x_1,x_2) = 1$.
% %
% Take few elements in 2-D. Create 2D plots of $u$, $u_h$ and standard LSQ solution.

% Can take $\xi_n = \xi_n(x_1)$, instead of~$\xi_n(x_1,x_2)$.

% \subsubsection{Further Numerics}
% \textbf{Skip for now...}
% \begin{itemize}
%     \item Discontinuous $\xi(x) = a_0+a_1\chi_{[z,1]}(x)$
%     \item $L^2 \rightarrow L^\infty$
%     \item penalty: $\alpha > 0$ (Lipschitz constant to be determined?)
% \end{itemize}

% \clearpage

% \subsection{Weighted Residual Minimization}

% \subsubsection{Numerics 1D: Minimize error in quantities (point values), given an irregular solution}

% \subsubsection{Numerics: Minimize $\mathrm{TV}$}
% Purpose: Create approximation (``unsupervised'') that minimizes Gibbs.

% %
% \subsubsection{Numerics 1D: Minimize error in quantities (point values), given an irregular solution}

% Purpose: Create coarse approximation that has small error for known quantities.

%\input{SecConclusion.tex}
%-----------------------------------------------------------------------------%
% \input{SecAcknowledgements.tex}

%-----------------------------------------------------------------------------%
%\clearpage
\appendix
\section{Proofs}

\subsection{Proof of Theorem~\ref{thm:quasi}}
\label{sec:quasi_proof}
\begin{enumerate}
\item[(i)] %We choose any sufficiently large $\alpha>0$ for which $j$ becomes strongly convex. 
Strong convexity of $j$ implies coercivity, i.e., 
$j(\xi)\to +\infty$ when $\|\xi\|_{\mathbb X} \to +\infty$. Moreover, $j$ is continuous in the strong topology since it is differentiable. Additionally, we know that convexity plus continuity implies that $j$ is weakly lower semicontinuous (see,~e.g.~\cite[Corollary~3.9]{BreBOOK2011}). We thus satisfy all the hypothesis of the \emph{theorem of existence of minimizers for coercive and sequentially weakly lower semicontinuous functionals}~\cite[Theorem~9.3-1]{CiaBOOK2013}. Moreover, strong convexity ensures that such a (global) minimizer $\bar\xi\in\mathbb X$ is unique. Besides, global differentiablity of $j$ implies the first-order necessary optimality condition $j'\big(\bar\xi\big)=0$.

\item[(ii)] 
We now that $j$ has a global lower bound. Thus, by the infimum property, for any $\delta_n>0$ there must exist $\bar\xi_n\in\mathcal M_n$ such that
\begin{equation}\label{eq:infimizing}
j(\bar\xi_n)<\inf_{\eta_n\in\mathcal M_n}j(\eta_n) + {\delta_n\over2}.
\end{equation}

% be an infimizing sequence, i.e., such that
% $$\lim_{m\to+\infty}j(\xi_m)=\inf_{\eta_n\in\mathcal M_n}j(\eta_n)\geq0.$$
% Coercivity of $j$ implies that $\{\xi_m\}$ is bounded. Reflexivity of $\mathbb X$ guaranties the existence of a subsequence (still denoted by $\{\xi_m\}$) converging weakly to some $\bar\xi_n\in\mathcal X_n$. Since $j$ is is weakly lower semicontinuous we get
% $$
% j(\bar\xi_n) \leq \liminf_{m\to+\infty}j(\xi_m)=\inf_{\eta_n\in\mathcal M_n}j(\eta_n),
% $$
% with equality if $\bar\xi_n\in\mathcal M_n$. If the case, the characterization of minimizers for differentiable functions in general sets (see Section~\ref{sec:opt_char})
% shows the condition~\eqref{eq:opt_char}.

\item[(iii)] 
%Let $\gamma>0$ be the strong convexity constant of $j$. 
Let $\bar\xi\in\mathbb X$ be the global minimizer and let $\bar\xi_n\in\mathcal M_n$ satisfy~\eqref{eq:qmin}. By characterization of strong convexity we have for all $t\in(0,1)$
$$
j\big(\bar\xi\big) \leq j\big(t\bar\xi + (1-t)\bar\xi_n\big)\leq
tj\big(\bar\xi\big) + (1-t) j\big(\bar\xi_n\big) - {\gamma\over2}t(1-t)\|\bar\xi-\bar\xi_n\|_{\mathbb X}^2\,. 
$$
Thus, for all $t\in(0,1)$ and $\eta_n\in\mathcal M_n
$ we get 
\begin{equation}\label{eq:coercivity}
{\gamma\over2}t\|\bar\xi-\bar\xi_n\|_{\mathbb X}^2
\leq j\big(\bar\xi_n\big)-j\big(\bar\xi\big) < 
j\big(\eta_n\big)-j\big(\bar\xi\big) + \frac{\delta_n}{2}.
\end{equation}
On the other hand, using the facts that $j'$ is $L$-Lipschitz and $j'\big(\bar\xi\big)=0$, we deduce~\cite[cf.~proof of Thm.~7.7-3, page~488]{CiaBOOK2013}
\begin{align}
j\big(\eta_n\big)-j\big(\bar\xi\big)= & \int_0^1\left<j'\big(s\eta_n + (1-s)\bar\xi\big),\eta_n-\bar\xi\right>ds\notag\\
= & \int_0^1\left<j'\big(s\eta_n + (1-s)\bar\xi\big)-j'\big(\bar\xi\big),\eta_n-\bar\xi\right>ds\notag\\
\leq & \,L \|\eta_n-\bar\xi\|_{\mathbb X}^2\int_0^1 s = {L\over2} \|\eta_n-\bar\xi\|_{\mathbb X}^2.
\label{eq:int_estimate}
\end{align}
\end{enumerate}
Hence, combining~\eqref{eq:coercivity} with~\eqref{eq:int_estimate}, taking the limit when $t\to 1$ and the infimum over all $\eta_n\in\mathcal M_n$, we get the estimate
$$
\gamma\,\|\bar\xi-\bar\xi_n\|^2_{\mathbb X}<
L\inf_{\eta_n\in\mathcal M_n}
\bignorm{\bar\xi-\eta_n}_{\mathbb X}^2
+ \delta_n\,,
$$
from which~\eqref{eq:quasi_opt_abstract} is deducted.

\subsection{Proof of Theorem~\ref{thm:reduced}}
\label{sec:reduced_proof}

We proceed to prove each one of the statements.
\begin{enumerate}
\item[(i)] Since $\mathbb Z$ and $\mathbb X$ are a Hilbert spaces, the quadratic maps 
$\mathbb Z\ni z\mapsto {1\over2}\|z\|_{\mathbb Z}^2$ and $\mathbb X\ni\xi\mapsto{1\over2}\|\xi\|^2_{\mathbb X}$ are differentiable. On the other hand, $S_h$ and $Q$ are also differentiable ($Q$ is linear), and thus $j_1$ is differentiable by means of the chain rule (see,~e.g.~\cite[Theorem~2.20]{TroBOOK2010}). Moreover,
$$
    j_1'(\eta)(\cdot) =  \big(\,QS_h(\eta)\,,\,QS_h'(\eta)(\cdot)\,\big)_{\mathbb Z}
    = \big(\,S_h'(\eta)^\star Q^\star QS_h(\eta)\,,\,\cdot\,\big)_{\mathbb X}\,.
$$
Thus, we conclude that $j_1$ is Lipschitz since
\begin{align*}
\big\|j_1'(\eta)-j_1'(\zeta)\big\|_{\mathbb X^*} = \, & \big\|S_h'(\eta)^\star Q^\star QS_h(\eta)-S_h'(\zeta)^\star Q^\star QS_h(\zeta)\big\|_{\mathbb X}\notag\\\\
= \, & \big\|S_h'(\eta)^\star Q^\star Q\big(S_h(\eta)-S_h(\zeta)\big)\big\|_{\mathbb X}\\
& + \big\|\big(S_h'(\eta)-S_h'(\zeta)\big)^\star Q^\star QS_h(\zeta)\big\|_{\mathbb X}\notag\\\\
\leq \, & \|Q\|_{\mathcal L(\mathbb U,\mathbb Z)}^2\big(M_{S'}^2  +L_{S'}M_S\big)\|\eta-\zeta\|_{\mathbb X}\,,
\end{align*}
where we have used the mean value theorem together with
\begin{itemize}
    \item the boundedness of $S_h'$, with bounding constant $M_{S'}$; 
    \item the Lipschitzness of $S_h'$, with Lipschitz constant $L_{S'}$; 
    \item the boundedness of $S_h$, with bounding constant $M_S$.
\end{itemize}
Finally, by making $L_1:=\|Q\|_{\mathcal L(\mathbb U,\mathbb Z)}^2\big(M_{S'}^2  +L_{S'}M_S\big)$, it is straightforward to see that $L_1+\alpha$ will be a Lipschitz constant for $j$.
\item[(ii)] Just observe that
\begin{align*}
\big<j'(\eta)-j'(\zeta),\eta-\zeta\big>_{\mathbb X^*,\mathbb X}
=\, & \big<j_1'(\eta)-j_1'(\zeta),\eta-\zeta\big>_{\mathbb X^*,\mathbb X} + \alpha \|\eta-\zeta\|_{\mathbb X}^2\\
\geq\, & (-L_1 +\alpha) \|\eta-\zeta\|_{\mathbb X}^2\,.
\end{align*}
Thus, $j$ is strongly convex whenever $\alpha>0$ is sufficiently large.
\end{enumerate}

% \begin{lemma}[Strong convexity] 
% \label{lem:StrongConvexity}
% \marginnote{@Ignacio: Lemma makes sense? Possible when $j_2 = \frac{1}{2}\norm{\cdot}_\infty^2$?}
% Let $j_1$ differentiable with $j_1'$ Lipschitz continuous, and let $j_2$ strongly convex. In other words, $\exists M_1 >0$:
% %
% \begin{alignat*}{2}
%  \Bignorm{ j'_1(\mu)-j_1'(\nu) }_{V^*} \le M_1 \bignorm{\mu-\nu}
% \end{alignat*}
% %
% and $\exists \gamma_2 >0$:
% %
% \marginnote{Useful? If only strictly convex, then $\ge (\norm{\mu}-\norm{\nu})^2 $.}
% %
% %
% \begin{alignat*}{2}
% \Bigdual{  \mu_2^*-\nu_2^*, \mu-\nu} 
% &\ge \gamma_2 \norm{\mu-\nu}^2
% &&\qquad  \forall \mu,\nu\in V\;,
% \\
% & &&\qquad \text{for any } \mu_2^*\in \partial j_2(\mu),\,\nu_2^*\in \partial j_2(\nu)\;,
% \end{alignat*}
% %
% Then, if $\gamma_2 > M_1$, $j = j_1 + j_2$ is strongly convex.
% %
% \end{lemma}
% %
% \begin{proof}
% %
% Let $\mu^* = j_1'(\mu) + \mu_2^*\in \partial j(\mu)$, $\nu^*=j_1'(\nu) + \nu_2^*\in \partial j(\nu)$, then
% %
% \begin{alignat*}{2}
% \Bigdual{  \mu^*-\nu^*, \mu-\nu} 
% &= \Bigdual{j_1'(\mu)-j_1'(\nu),\mu-\nu} +
% \Bigdual{  \mu_2^*-\nu_2^*, \mu-\nu} 
% \\
% &\ge 
%  \big({-} M_1 + \gamma_2\big) \bignorm{\mu-\nu}^2
%  \;.
% \end{alignat*}
% %
% \end{proof}

\subsection{Proof of Proposition~\ref{prop:stateEq}}
\label{sec:wellposed_proof}
The statements~(i)~and~(ii) are classical from Babu\v{s}ka–Brezzi theory (see, e.g.,  Ern~\&~Guermond~\cite[Theorem 49.13]{ErnGueBOOK2021b}).
To prove statment~(iii) first observe that 
$$
\sup_{v_2\in\hat{\mathbb{K}}}{a(\xi;v_1,v_2)\over\|v_1\|_\mathbb V\|v_2\|_\mathbb V}\geq 
{a(\xi;v_1,v_1)\over\|v_1\|_\mathbb V\|v_1\|_\mathbb V}\geq {a(\xi;v_1,v_1)\over\|v_1\|_{\mathbb V,\xi}\|v_1\|_{\mathbb V,\xi}}
(C_{1,\xi})^2=(C_{1,\xi})^2,$$
which confirms  $\alpha_h=(C_{1,\xi})^2$ in~\eqref{eq:inf-sup_a}. For the a~priori bound, since $a(\xi;\cdot,\cdot)$ is an equivalent inner-product on $\hat{\mathbb V}$, consider $\hat z\in\hat{\mathbb V}$ such that
$$
a(\xi;\hat z,\hat v)=b(u_{h},\hat v),\quad\forall \hat v\in\hat{\mathbb V}. 
$$
Hence,
\begin{equation}\label{eq:supremizer}
\sup_{\hat v\in\hat{\mathbb V}}
{b(u_{h},\hat v)\over \|\hat v\|_{\mathbb V,\xi}}
=\sup_{\hat v\in\hat{\mathbb V}} {a(\xi;\hat z,\hat v)\over \|\hat v\|_{\mathbb V,\xi}}
= {a(\xi;\hat z,\hat z)\over \|\hat z\|_{\mathbb V,\xi}}
= {b(u_{h},\hat z)\over \|\hat z\|_{\mathbb V,\xi}}.
\end{equation}
Moreover, 
\begin{equation}\label{eq:z-ortho}
    a(\xi,\hat r,\hat z)=a(\xi,\hat z,\hat r)=
    b(u_{h},\hat r)=0.
\end{equation}
Next, observe that
\begin{alignat}{2}
\|u_{h}\|_\mathbb U\, \leq\, 
& \, {1\over\beta_h}\,
\sup_{\hat v\in\hat{\mathbb V}}
{b(u_{h},\hat v)\over \|\hat v\|_\mathbb V}\,
\leq \,{C_{2,\xi}\over\beta_h}\,
\sup_{\hat v\in\hat{\mathbb V}}
{b(u_{h},\hat v)\over \|\hat v\|_{\mathbb V,\xi}}
\tag{by~\eqref{eq:inf-sup} and~\eqref{eq:norm_equiv}}\\
= \, & \, {C_{2,\xi}\over\beta_h}\,
{b(u_{h},\hat z)\over \|\hat z\|_{\mathbb V,\xi}}\,
= \, {C_{2,\xi}\over\beta_h}\,
{\big(f(\hat z)-a(\xi,\hat r,\hat z)\big)\over \|\hat z\|_{\mathbb V,\xi}}
\tag{by~\eqref{eq:supremizer} and~\eqref{eq:StateEq}}\\
\leq \, & \, {C_{2,\xi}\over C_{1,\xi}}\,{1\over\beta_h}\,
{f(\hat z)\over \|\hat z\|_{\mathbb V}}\,,
%\leq {C_2\over C_1}{1\over\beta_h}\|f\|_{\mathbb V^*}\,,
\tag{by~\eqref{eq:norm_equiv} and~\eqref{eq:z-ortho}} 
\end{alignat}
from which~\eqref{eq:apriori} can be easily deducted.

\subsection{Proof of Proposition~\ref{prop:equiv}}
\label{sec:equiv_proof}
% The conditions in \eqref{eq:inf-sup_a} imply that $A(\xi)$ and $A(\xi)^*$ are invertible.
% \begin{enumerate}
% \item[(i)] Since $A(\xi)^*$ is invertible, $A(\xi)^{-*} B w_h$ exists in~$\mathbb{V}$ for any~$w_h\in \mathbb{U}_h$.
% %
% \item[(ii)] Since $A(\xi)$ is invertible, solving \eqref{eq:mixed_operator_a} for~$r$, and substituting this in~\eqref{eq:mixed_operator_b} yields:
% \begin{alignat*}{2}
%   B^* A(\xi)^{-1} \Big(f-Bu_h\Big) = 0\qquad  \text{in } (\mathbb{U}_h)^*\,,
% \end{alignat*}
% i.e., 
% \begin{alignat*}{2}
%   \Bigdual{B^* A(\xi)^{-1} Bu_h, w_h}  = \bigdual{B^* A(\xi)^{-1} f, w_h} \qquad  \forall w_h\in \mathbb{U}_h\,.
% \end{alignat*}
% This is equivalent  to~\eqref{eq:introPGxi} with~$\mathbb{V}_h(\xi)$ as in~\eqref{eq:Vhxi}.
% \item[(iii)] 
% We proceed by directly verifying the discrete inf-sup condition:
% \begin{alignat*}{2}
%  & \sup_{v_h\in \mathbb{V}_h(\xi)}
%  \frac{b(w_h,v_h)}{\norm{v_h}_{\mathbb{V}}}
%  = \sup_{z_h\in \mathbb{U}_h}
%  \frac{\bigdual{B w_h, A(\xi)^{-*} B z_h)}}{\norm{A(\xi)^{-*} B z_h}_{\mathbb{V}}}
% \\
%  & \quad = \sup_{z_h\in \mathbb{U}_h}
%  \frac{a\big(\xi; A(\xi)^{-*} Bz_h, A(\xi)^{-*} B w_h)\big)}{\norm{A(\xi)^{-*} B z_h}_{\mathbb{V}}}
% \end{alignat*}

%\item[(i)] 
    Let $(r,u_h)\in \hat{\mathbb V}\times\mathbb U_h$ solves the state problem~\eqref{eq:StateEq}, or equivalently~\eqref{eq:mixed_operator} in operator form. Testing with elements in $v_h\in\mathbb V_h(\xi)$ we get
    \begin{alignat}{2}
    \left<f,v_h\right> = &
    \left<A(\xi)r,v_h\right>+\left<Bu_h,v_h\right> \tag{by~\eqref{eq:mixed_operator_a}}\\
    = & \left<r,A(\xi)^*v_h\right>+\left<Bu_h,v_h\right>\tag{using the adjoint property}\\
    = & \left<r,Bw_h\right>+\left<Bu_h,v_h\right> \tag{by definition of $\mathbb V_h(\xi)$}\\
    = & \left<Bu_h,v_h\right>\,. \tag{by~\eqref{eq:mixed_operator_b}}
    \end{alignat}
    Thus, \eqref{eq:introPGxi} is satisfied.
    
Conversely, assume $u_h\in \mathbb{U}_h$ satisfies the Petrov--Galerkin problem~\eqref{eq:introPGxi} with~$A(\xi)^{*}\mathbb{V}_h(\xi) =  B \mathbb{U}_h$. In particular, 
$$\bigdual{f-Bu_h,v_0} =0 \qquad \forall v_0 \in \ker A(\xi)^* = \Big\{ v_0\in \hat{\mathbb{V}}\,\Big|\, A(\xi)^*v_0 = 0\Big\} \subset \mathbb{V}_h(\xi) \,.$$
Hence, by orthogonality,%
\footnote{That is, $(\ker(A^*))^\perp = \operatorname{ran} A$; see, e.g., \cite[Lemma~C.34]{ErnGueBOOK2021b}.}
there exists an $r \in \hat{\mathbb{V}}$ such that $A(\xi)r  = f-B u_h$, which is~\eqref{eq:mixed_operator_a}. 
\par
Next, let $w_h \in \mathbb{U}_h$. Since~$A(\xi)^*$ is surjective (by~\eqref{eq:inf-sup_a2}), there exists a~$v_{w_h} \in \hat{\mathbb{V}}$ such that $A(\xi)^*v_{w_h} = B w_h$, in other words, $v_{w_h} \in \mathbb{V}_h(\xi)$. Therefore, 
$$
 \bigdual{B w_h,r} = \bigdual{ A(\xi)^*v_{w_h},r} 
 = \bigdual{ A(\xi)r, v_{w_h}} 
 = \bigdual{ f-Bu_h, v_{w_h}}=0\,,
$$
which verifies~\eqref{eq:mixed_operator_b}.

%     Conversely, assume that equation~\eqref{eq:introPGxi} is satisfied. 
%     Hence, $f-Bu_h\in\ \mathbb V_h(\xi)^\bot$. On another hand, observe that
%     $\ker A(\xi)^*\subset \mathbb V_h(\xi)$. In particular, 
%     $f-Bu_h\in (\ker A(\xi)^*)^\bot$, which means that $f-Bu_h$ is in the image of $A(\xi)$. Thus, there is $r\in \hat{\mathbb V}$ such that $A(\xi)r=f-Bu_h$, which is equation~\eqref{eq:mixed_operator_a}. 
%     %Finally, given $w_h\in\mathbb U_h$, we decompose $Bw_h=$  
%     Then, for $v_h \in \mathbb{V}_h(\xi)$, 
% \begin{alignat*}{2}
%  \bigdual{A(\xi) r,v_h} = 0
% \end{alignat*}
% hence, we obtain~\eqref{eq:mixed_operator_b}: 
% \begin{alignat*}{2}
%  \bigdual{r,B w_h} = 0
%  \qquad \forall w_h\in \mathbb{U}_h\,.
% \end{alignat*}
% \end{enumerate}

\subsection{Proof of Proposition~\ref{prop:state_deriv}}
\label{sec:state_deriv_proof}
Let us start proving statements~(i),~(ii)~and~(iii) at the same time.

Recall the definition of the kernel space $\hat{\mathbb{K}}:=\ker B^*\subset\hat{\mathbb V}$. For any $\xi\in\mathbb X$, consider the restricted operator $A(\xi)\big|_{\hat{\mathbb{K}}}:~\hat{\mathbb{K}}\to\hat{\mathbb{K}}^*$, as well as the restriction $f\big|_{\hat{\mathbb{K}}}\in~\hat{\mathbb{K}}^*$. Observe that the $\inf$-$\sup$ condition~\eqref{eq:inf-sup} ensures that 
$A(\xi)\big|_{\hat{\mathbb{K}}}$ is a boundedly invertible linear operator. Thus, given a direction $\eta\in\mathbb X$ and $t\in\mathbb R$, from the first equation of the mixed system~\eqref{eq:mixed_operator} (restricted to $\hat{\mathbb{K}}$) we obtain that
\begin{subequations}
\begin{empheq}[left=,right=,box=]{alignat=3} 
 A(\xi+t\eta)\big|_{\hat{\mathbb{K}}}R_h(\xi+t\eta) & =  f\big|_{\hat{\mathbb{K}}}  \label{eq:restricted_1}\\
% \notag\\
   A(\xi)\big|_{\hat{\mathbb{K}}}R_h(\xi) & =  f\big|_{\hat{\mathbb{K}}}\,\,.\label{eq:restricted_2}
\end{empheq}
\end{subequations}
In particular, continuity of $A(\cdot)$ implies continuity of $R_h(\cdot)$. Moreover, using the $\inf$-$\sup$ condition~\eqref{eq:inf-sup}, it is clear that
\begin{equation}\label{eq:r-estimate}
\|R_h(\cdot)\|_\mathbb V\leq {\|f\|_{\hat{\mathbb V}^*}
\over\alpha_h(\cdot)}\,.
\end{equation}
Next, adding the term $A(\xi)\big|_{\hat{\mathbb{K}}} R_h(\xi+t\eta)$ on both sides of equation~\eqref{eq:restricted_1}, rearrange it, and subtracting equation~\eqref{eq:restricted_2} we get 
$$
R_h(\xi+t\eta)-R_h(\xi)=\left[A(\xi)\big|_{\hat{\mathbb{K}}}\right]^{-1}\left(A(\xi)\big|_{\hat{\mathbb{K}}}-A(\xi+t\eta)\big|_{\hat{\mathbb{K}}}\right)R_h(\xi+t\eta),
$$
from which, if $A'(\xi)\eta$ exists, we imply that $R_h(\cdot)$ has a G\^ateaux derivative and
\begin{equation}\label{eq:r'}
R_h'(\xi)\eta = - \left[A(\xi)\big|_{\hat{\mathbb{K}}}\right]^{-1}A'(\xi)\eta\Big|_{\hat{\mathbb{K}}}\,R_h(\xi).
\end{equation}
Finally, if $A(\cdot)$ is G\^ateaux-differentiable at $\xi$, then using the $\inf$-$\sup$ condition~\eqref{eq:inf-sup}, the boundedness of the linear operator $A'(\xi)$, and the estimate~\eqref{eq:r-estimate}, we imply
\begin{equation}\label{eq:r'estimate}
\|R_h'(\xi)\eta\|_{\mathbb V}\leq 
%{1\over\alpha_h}\sup_{\hat v\in\hat{\mathbb{K}}}{a\left(\xi;\hat r'(\xi)\eta,\hat v\right)\over\|\hat v\|_{\mathbb V}}\leq
{\|A'(\xi)\eta\|_{\mathcal L(\hat{\mathbb V},\hat{\mathbb V}^*)}\|R_h(\xi)\|_{\mathbb V}\over\alpha_h}
\leq {\|A'(\xi)\|\|f\|_{\mathbb V^*}\over\alpha_h^2}\|\eta\|_{\mathbb X}\,,
\end{equation}
which proves that $R_h(\cdot)$ is  G\^ateaux-differentiable at $\xi$. Besides, if $A'(\cdot)$ and $\alpha_h^{-1}(\cdot)$ are uniformly bounded on $\mathbb X$, then $R_h'(\cdot)$ is uniformly bounded on $\mathbb X$.

Now is the turn of $S_h$.
From the mixed system~\eqref{eq:mixed_operator} we deduce
$$
\begin{array}{rl}
 BS_h(\xi+t\eta) =  & f-A(\xi+t\eta)R_h(\xi+t\eta)\\
  BS_h(\xi) = & f-A(\xi)R_h(\xi).
\end{array}
$$
Since $B$ is boundedly invertible onto its closed range we get
$$
S_h(\xi+t\eta)-S_h(\xi)=B^{-1}\Big(
[A(\xi)-A(\xi+t\eta)]R_h(\xi+t\eta) 
+A(\xi)[R_h(\xi)-R_h(\xi+t\eta)]\Big).
$$
Therefore, if $A'(\xi)\eta$ exists, then we already know that $R_h'(\xi)\eta$ exists, and thus
\begin{equation}\label{eq:u_h'}
S_h'(\xi)\eta = B^{-1}\Big(-
[A'(\xi)\eta]R_h(\xi) - A(\xi)R_h'(\xi)\eta\Big).
\end{equation}
Moreover, if $A(\cdot)$ is G\^ateaux-differentiable, then using the $\inf$-$\sup$ condition~\eqref{eq:inf-sup} and the estimate~\eqref{eq:r'estimate} we get
\begin{equation}\label{eq:u_h'-estimate}
\begin{array}{rl}
\|S_h'(\xi)\eta\|_{\mathbb U}
\leq & \displaystyle
{1\over\beta_h}
\|B[S_h'(\xi)\eta]\|_{\hat{\mathbb V}^*}\\\\
\leq & \displaystyle {\|A'(\xi)\|\|R_h(\xi)\|_{\mathbb V}+\|A(\xi)\|_{\mathcal L(\hat{\mathbb V},\hat{\mathbb V}^*)}\|R_h'(\xi)\|_{\mathcal L(\mathbb X,\hat{\mathbb V})}\over\beta_h}\|\eta\|_{\mathbb X}\\\\
\leq & \displaystyle
{\|A'(\xi)\|\|f\|_{\mathbb V^*}\over\alpha_h\beta_h}
\left(1+{\|A(\xi)\|_{\mathcal L(\hat{\mathbb V},\hat{\mathbb V}^*)}\over\alpha_h}\right)
\|\eta\|_{\mathbb X}\,,
\end{array}
\end{equation}
which proves that $S_h(\cdot)$ is G\^ateaux-differentiable. Besides, it is clear from~\eqref{eq:u_h'-estimate} that $\|S_h'(\cdot)\|_{\mathcal L(\mathbb X,\mathbb U)}$ will be uniformly bounded on $\mathbb X$, whenever $\|A(\cdot)\|_{\mathcal L(\hat{\mathbb V},\hat{\mathbb V}^*)}$ and $\|A'(\cdot)\|$ are uniformly bounded on $\mathbb X$, as well as $\alpha_h^{-1}(\cdot)$. 

(iv) Let us prove Lipschitzness. Using~\eqref{eq:r'}, observe that for any
$\xi_1,\xi_2,\eta\in\mathbb X$ we have
\begin{align}\notag
A(\xi_2)\big|_{\hat{\mathbb{K}}}\big(R_h'(\xi_1)-R_h'(\xi_2)\big)\eta = & 
\big[A'(\xi_2)-A'(\xi_1)\big]\eta\Big|_{\hat{\mathbb{K}}}
R_h(\xi_2) + \big[A(\xi_2)-A(\xi_1)\big]\Big|_{\hat{\mathbb{K}}}R_h'(\xi_1)\eta\\\notag
& +A'(\xi_1)\eta\Big|_{\hat{\mathbb{K}}}
\big[R_h(\xi_2)-R_h(\xi_1)\big].
\end{align}
Hence,
\begin{subequations}
\begin{empheq}[left=,right=,box=]{alignat=3} 
\label{eq:r'-Lips_1}
\|R_h'(\xi_1)-R_h'(\xi_2)\|_{\mathcal L(\mathbb X,\hat{\mathbb V})}\leq &
{\|R_h(\xi_2)\|_{\mathbb V}\over\alpha_h(\xi_2)}\|A'(\xi_1) - A'(\xi_2)\|\\
& + \label{eq:r'-Lips_2}
{\|R_h'(\xi_1)\|_{\mathcal L(\mathbb X,\hat{\mathbb V})}\over\alpha_h(\xi_2)}\|A(\xi_1) - A(\xi_2)\|_{\mathcal L(\hat{\mathbb V},\hat{\mathbb V}^*)}\\
& + \label{eq:r'-Lips_3}
{\|A'(\xi_1)\|\over\alpha_h(\xi_2)}\|R_h(\xi_1) - R_h(\xi_2)\|_{\mathbb V}\,.
\end{empheq}
\end{subequations}
Recall that under our hypothesis, $\alpha_h^{-1}(\cdot)$, $R_h(\cdot)$, $R_h'(\cdot)$, and $A'(\cdot)$, they are all uniformly bounded on $\mathbb X$. Therefore, the first term on the right hand side (expression~\eqref{eq:r'-Lips_1}) is Lipschitz by the Lipschitz assumption on $A'(\cdot)$;
the second term (expression~\eqref{eq:r'-Lips_2}) is Lipschitz as a consequence of the mean value theorem on $A(\cdot)$ and the uniform boundedness of $A'(\cdot)$; while the last term 
(expression~\eqref{eq:r'-Lips_3}) is Lipschitz by the mean value theorem on $R_h(\cdot)$ and the uniform boundedness of $R_h'(\cdot)$.

Finally, to prove the Lipschitzness of $S_h'(\cdot)$, we use~\eqref{eq:u_h'} to write
\begin{align}\notag
B\big(S_h'(\xi_1)\eta-S_h'(\xi_2)\eta\big)= &
[A'(\xi_2)\eta]\big(R_h(\xi_2)-R_h(\xi_1)\big) + A(\xi_2)\big[R_h'(\xi_2)\eta
-R_h'(\xi_1)\eta\big]\\
\notag & +
\big[(A'(\xi_2)-A'(\xi_1))\eta\big]R_h(\xi_1) +\big[A(\xi_2)- A(\xi_1)\big]R_h'(\xi_1)\eta\,.
\end{align}
Hence,
\begin{subequations}
\begin{empheq}[left=,right=,box=]{alignat=3} 
\label{eq:u_h'-Lips_1}
\big\|S_h'(\xi_1)-S_h'(\xi_2)\big\|_{\mathcal L(\mathbb X,\mathbb U)}
\leq &
{\|A'(\xi_2)\|\over\beta_h}
\|R_h(\xi_1)-R_h(\xi_2)\|_\mathbb V\\
\label{eq:u_h'-Lips_2}
& + {\|A(\xi_2)\|_{\mathcal L(\hat{\mathbb V},\hat{\mathbb V}^*)}\over\beta_h}
\|R_h'(\xi_1)-R_h'(\xi_2)\|_{\mathcal L(\mathbb X,\hat{\mathbb V})}\\
\label{eq:u_h'-Lips_3}
& + {\|R_h(\xi_1)\|_\mathbb V\over\beta_h}
\|A'(\xi_1)-A'(\xi_2)\|\\
\label{eq:u_h'-Lips_4}
& + {\|R_h'(\xi_1)\|_{\mathcal L(\mathbb X,\hat{\mathbb V})}\over\beta_h}
\|A(\xi_1)-A(\xi_2)\|_{\mathcal L(\hat{\mathbb V},\hat{\mathbb V}^*)}\,.
\end{empheq}
\end{subequations}
We recall again that $R_h(\cdot)$, $R_h'(\cdot)$, $A(\cdot)$, and $A'(\cdot)$, they are all uniformly bounded on $\mathbb X$. Therefore, the Lipschitzness of $S_h'(\cdot)$ is implied by the following facts:
the Lipschitzness of the first term on right hand side (expression~\eqref{eq:u_h'-Lips_1}) is a consequence of the mean value theorem applied to $R_h(\cdot)$ and the uniform boundedness of $R_h'(\cdot)$; the Lipschitzness of the second term (expression~\eqref{eq:u_h'-Lips_2}) is due to the previously proved Lipschitzness of $R_h'(\cdot)$; the Lipschitzness of the third term (expression~\eqref{eq:u_h'-Lips_3}) is implied by the assumed Lipschitzness of $A'(\cdot)$;
and the Lipschitzness of the last term (expression~\eqref{eq:u_h'-Lips_4}) is  consequence of the mean value theorem applied to $A$ and the uniform boundedness of $A'(\cdot)$.

% \begin{alignat*}{2}
%   \norm{r'_1-r'_2}
%   &\le C \norm{A_1 r'_1- A_1 r'_2}
% \\
%   &\le C\bigg( \norm{A_1 r'_1- A_2 r'_2}
%   + \norm{A_2 r'_2- A_1 r'_2}
%   \bigg)
% \\
%   &= C\bigg( \norm{{-}A'_1 r_1 + A'_2 r_2}
%   + \norm{A_2 r'_2- A_1 r'_2}
%   \bigg)
% \\
%   &\le C\bigg( \norm{(A'_1 - A'_2) r_1}
%   +\norm{A'_2( r_1-r_2)}
%   + \norm{(A_2- A_1) r'_2}
%   \bigg)\end{alignat*}

\subsection{Proof of Proposition~\ref{prop:LSQ_well}}
\label{sec:LSQ_well}
Let us prove item by item.
\begin{itemize}
    \item[(i)] Observe that in this case, the bilinear form $a(\xi,\cdot,\cdot)$ defines a weighted inner product in $L^2(\Omega)$, for which its induced norm $\|v\|_{\mathbb V,\xi}:=\sqrt{(\varpi(\xi)v,v)_{L^2}}$ satisfies
$$
\sqrt{\varpi_{\min}}\,\|v\|_{L^2}\leq \|v\|_{\mathbb V,\xi}\leq \sqrt{\varpi_{\max}}\,\|v\|_{L^2}\,,\quad\forall v\in \mathbb V=L^2(\Omega).
$$
Hence, the first $\inf$-$\sup$ condition in~\eqref{eq:inf-sup} is satisfied with 
$\alpha_h=\varpi_{\min}$; see Proposition~\ref{prop:stateEq}(iii) and Footnote~\ref{ftnt:a-surjectivity}.

On the other hand, we are under the assumption that the operator $B:\mathbb{H}_B\to\mathbb V^*$ is boundedly invertible. Hence, there must be a uniform constant $\beta>0$ such that
$$
\sup_{v\in\mathbb V}{b(w_h,v)\over\|v\|_\mathbb V}=\|Bw_h\|_{\mathbb V^*}
\geq \beta\|w_h\|_{\mathbb{H}_B}\,,\quad\forall w_h\in\mathbb U_h\,,
$$
which implies the second $\inf$-$\sup$ condition in~\eqref{eq:inf-sup}.
\item[(ii)]
Uniform boundedness of $S_h(\cdot)$ is a consequence of Proposition~\ref{prop:stateEq}(iii). Indeed, in our particular case we get
$$
\|S_h(\xi)\|_{\mathbb{H}_B}\leq {\varpi_{\max}\over\varpi_{\min}}
{1\over\beta}\|f\|_{L^2}\,,\quad\forall\xi\in L^2(\Omega).
$$
To show differentiability of $S_h(\cdot)$, let us recall the operator $A:\mathbb X\to\mathcal L(\mathbb V,\mathbb V^*)$ defined in section~\ref{sec:analysis}, which in this particular case, given $\xi\in L^2(\Omega)$, it takes the form
$$
A(\xi)v:=(\varpi(\xi)v,\cdot)_{L^2}\,,\quad \forall v\in L^2(\Omega).
$$
Furthermore, we have the uniform bound
\begin{equation}\label{eq:A_bound_LSQ}
\|A(\xi)\|=\sup_{v\in L^2(\Omega)}{\|\varpi(\xi)v\|_{L^2}\over \|v\|_{L^2}} 
=\|\varpi(\xi)\|_{L^\infty}\leq \varpi_{\max}\,.
\end{equation}

Since $\varpi(\cdot)$ is differentiable, it is straightforward to check that $A(\cdot)$ is also differentiable, and given $\xi,\eta\in L^2(\Omega)$, we have
$$
[A'(\xi)\eta]v = \big([\varpi'(\xi)\eta]v,\cdot\big)_{L^2}\,,\quad\forall v\in L^2(\Omega).
$$
Moreover, we can verify
\begin{equation}\label{eq:A'bound_LSQ}
\|A'(\xi)\|=\sup_{\eta\in L^2(\Omega)}{\|\varpi'(\xi)\eta\|_{L^\infty}\over \|\eta\|_{L^2}}=\|\varpi'(\xi)\|_{\mathcal L(L^2(\Omega),L^\infty(\Omega))}\leq \varpi'_\infty\,.
\end{equation}
Thus, the differentiablity of $S_h(\cdot)$ is a consequence of Proposition~\ref{prop:state_deriv}(ii).

\item[(iii)]
Uniform boundedness of $S_h'(\cdot)$ is a consequence of Proposition~\ref{prop:state_deriv}(iii), using the fact that $A(\cdot)$, $A'(\cdot)$, and $\alpha_h^{-1}\equiv\varpi_{\min}^{-1}$, are all uniformly bounded (see the above expressions~\eqref{eq:A_bound_LSQ}~and~\eqref{eq:A'bound_LSQ}). 

On the other hand, the Lipschitz-continuity of $S_h'(\cdot)$ relies on the Lipschitz-continuity of $A'(\cdot)$ (by Proposition~\ref{prop:state_deriv}(iv)). 
The latter is true since
$$
\norm{ A'(\xi_1)-A'(\xi_2) } =
\sup_{\eta\in L^2(\Omega)}
\frac{\norm{ \varpi'(\xi_1)\eta-\varpi'(\xi_2)\eta}_{L^\infty} }{\norm{\eta}_{L^2}} \leq \varpi_L \|\xi_1-\xi_2\|_{L^2} \,.
$$
\end{itemize}

\subsection{Proof of Proposition~\ref{prop:WPG-well}}
\label{sec:WPG-well}

\begin{itemize}
    \item[(i)] We verify the hypothesis of Proposition~\ref{prop:stateEq}(i).
Since $\mathbb V$ is infinite dimensional, we first need to show that
$$
\{v_2\in K : b(\varpi(\xi)v_2,v_1)=0\,,\forall v_1\in K\}=\{0\},
$$
which is an immediate consequence of~\eqref{eq:a-ker} taking $v_1=v_2$.
To show the $\inf$-$\sup$ conditions~\eqref{eq:inf-sup}, from one hand observe that
$$
\sup_{v_2\in K} {b(\varpi(\xi)v_2,v_1)\over\|v_2\|_\mathbb V}
\geq  {b(\varpi(\xi)v_1,v_1)\over\|v_1\|_\mathbb V}\geq \alpha_h(\xi)
\|v_1\|_\mathbb V\,,\quad\forall v_1\in K.
$$
On the other hand, we have that $b(\cdot,\cdot)$ satisfies~\eqref{eq:b-coercive}. Thus, in particular
$$
\sup_{v\in\mathbb V}{b(v_h,v)\over\|v\|_\mathbb V}\geq \beta
\|v_h\|_\mathbb V\,,\quad\forall v_h\in\mathbb V_h\,.
$$
    \item[(ii)] We use the a~priori bound of Proposition~\ref{prop:stateEq}(ii). In this case $\beta_h=\beta$, $\alpha_h^{-1}(\xi)\leq \alpha^{-1}$, and 
    $a(\xi;v_1,v_2)=b\big(\varpi(\xi)v_2,v_1\big)$, for all $v_1,v_2\in\mathbb V$.
    It is easy to see that $\|a(\xi;\cdot,\cdot)\|_{\mathcal L(\mathbb V\times\mathbb V;\mathbb R)}\leq \|b\|_{\mathcal L(\mathbb V\times\mathbb V;\mathbb R)}\|\varpi(\xi)\|_\mathbb{W}$. Thus, we get
    $$
    \|S_h(\cdot)\|_\mathbb V\leq {1\over\beta} \left(1+{\varpi_\infty\over\alpha}\right)\|f\|_{\mathbb V^*}.
    $$
    \item[(iii)] Now we apply Proposition~\ref{prop:state_deriv}. The operator $A:\mathbb X\to\mathcal L(\mathbb V,\mathbb V^*)$ takes the form $A(\xi)v=b\big(\varpi(\xi)\,\cdot\,,v\big)\in\mathbb V^*$, for all $\xi\in \mathbb X= L^2(\Omega)$ and $v\in\mathbb V$. Moreover,
    $$
    \|A(\xi)\|_{\mathcal L(\mathbb V,\mathbb V^*)}\leq \|b\|_{\mathcal L(\mathbb V\times\mathbb V;\mathbb R)}\|\varpi(\xi)\|_\mathbb{W}\leq 
    \|b\|_{\mathcal L(\mathbb V\times\mathbb V;\mathbb R)}\,\varpi_\infty\,.
    $$
    On the other hand, it is immediate to see that if $\varpi$ is differentiable, then $A$ is differentiable and $[A'(\xi)\eta]v=b\big([\varpi'(\xi)\eta]\,\cdot\,,v\big)\in\mathbb V^*$, for any direction $\eta\in L^2(\Omega)$. 
    Moreover,
    $$
    \|A'(\xi)\eta\|_{\mathcal L(\mathbb V,\mathbb V^*)}\leq \|b\|_{\mathcal L(\mathbb V\times\mathbb V;\mathbb R)}\|\varpi'(\xi)\|_{\mathcal L(L^2(\Omega),\mathbb{W})}\|\eta\|_{L^2}\,.
    $$
    Hence, $A'(\cdot)$ is uniformly bounded and Lipschitz-continuous whenever $\varpi'(\cdot)$ is. By Proposition~\ref{prop:state_deriv}, differentiability of $S_h(\cdot)$ is implied by differentiability of $A(\cdot)$; uniform boundedness of $S_h'(\cdot)$ is implied by uniform boundedness of $A(\cdot)$, $A'(\cdot)$ and $\alpha_h^{-1}(\cdot)$; while Lipschitzness of $S_h'(\cdot)$ is implied by Lipschitzness of $A'(\cdot)$.
\end{itemize}

\subsection{Proof of Proposition~\ref{prop:ddminres}}
\label{sec:ddminres_proof}
\begin{itemize}
    \item[(i)] Making the identification $\hat{\mathbb V}\equiv\mathbb V_h$ and
    $a(\xi;\cdot,\cdot)\equiv(\cdot,\cdot)_{\mathbb V,\xi}\,$, we observe that the well-posedness of~\eqref{eq:ddminres} is a direct consequence of Proposition~\ref{prop:stateEq}, using the fact that $(\cdot,\cdot)_{\mathbb V,\xi}$ is an equivalent inner-product, together with assumption~\eqref{eq:discrete_infsup_b}.  
    \item[(ii)] Using the hypothesis of this statement and the estimate~\eqref{eq:apriori} in Proposition~\ref{prop:stateEq}(iii), we get the uniform bound  
    $$
    \|S_h(\xi)\|_\mathbb U\leq {1\over\beta_h}{\tilde C_2\over\tilde C_1}\|f\|_{\mathbb V^*}\,,\quad\forall \xi\in \mathbb X\,.
    $$
    \item[(iii)] Direct application of Proposition~\ref{prop:state_deriv}, 
    noticing also that $\alpha_h^{-1}(\xi)\leq \tilde C_1^{-2}$ and
    \begin{alignat}{2}
    \notag
    \|A(\xi)\|_{\mathcal L(\mathbb V,\mathbb V^*)} = &
    \sup_{v_1\in\mathbb V}{\|(v_1,\cdot)_{\mathbb V,\xi}\|_{\mathbb V^*}\over\|v_1\|_\mathbb V}\\
    \notag
    \leq & \sup_{v_1\in\mathbb V}{\tilde C_2^2\over\|v_1\|_{\mathbb V,\xi}}
    \left(\sup_{v_2\in\mathbb V}{|(v_1,v_2)_{\mathbb V,\xi}|\over\|v_2\|_{\mathbb V,\xi}}\right)\\
    \notag
    = & \, \tilde C_2^2\,.
    \end{alignat}
\end{itemize}
%\input{SecOld.tex}
%-----------------------------------------------------------------------------%
%
%
%=============================================================================%
\clearpage
\small
\bibliography{BibFileKris}

\begin{thebibliography}{10}

\bibitem{AinDonSISC2021}
{\sc M.~Ainsworth and J.~Dong}, {\em Galerkin neural networks: {A} framework
  for approximating variational equations with error control}, SIAM J. Sci.
  Comput., 43 (2021), pp.~A2474--A2501.

\bibitem{BarHoyHicBrePNAS2019}
{\sc Y.~Bar-Sinai, S.~Hoyer, J.~Hickey, and M.~P. Brenner}, {\em Learning
  data-driven discretizations for partial differential equations}, Proceedings
  of the National Academy of Sciences, 116 (2019), pp.~15344--15349.

\bibitem{BerNysNC2018}
{\sc J.~Berg and K.~Nystr\"om}, {\em A unified deep artificial neural network
  approach to partial differential equations in complex geometries},
  Neurocomputing, 317 (2018), pp.~28--41.

\bibitem{BerNysJCMDS2021}
{\sc J.~Berg and K.~Nystr\"om}, {\em Neural networks as smooth priors for
  inverse problems for {PDEs}}, Journal of Computational Mathematics and Data
  Science, 1 (2021), p.~100008.

\bibitem{BocGunBOOK-CH2016}
{\sc P.~Bochev and M.~Gunzburger}, {\em Chapter 12 - {L}east-squares methods
  for hyperbolic problems}, in Handbook of Numerical Methods for Hyperbolic
  Problems, R.~Abgrall and C.-W. Shu, eds., vol.~17 of Handbook of Numerical
  Analysis, Elsevier, 2016, pp.~289--317.

\bibitem{BohFeiCAMWA2021}
{\sc J.~Bohn and M.~Feischl}, {\em Recurrent neural networks as optimal mesh
  refinement strategies}, Comput. Math. Appl., 97 (2021), pp.~61--76.

\bibitem{BorSchBOOK2012}
{\sc A.~Borz\`{i} and V.~Schulz}, {\em Computational Optimization of Systems
  Governed by Partial Differential Equations}, Siam series on Computational
  Science and Engineering, Society for Industrial and Applied Mathematics,
  2012.

\bibitem{BreMugZeeCAMWA2021}
{\sc I.~Brevis, I.~Muga, and K.~G. {van der Zee}}, {\em A machine-learning
  minimal-residual {(ML-MRes)} framework for goal-oriented finite element
  discretizations}, Comput. Math. Appl., 95 (2021), pp.~186--199.
\newblock Recent Advances in Least-Squares and Discontinuous Petrov–Galerkin
  Finite Element Methods.

\bibitem{BreBOOK2011}
{\sc H.~Brezis}, {\em Functional Analysis, Sobolev Spaces and Partial
  Differential Equations}, Universitext, Springer, New York, 2011.

\bibitem{BurErnMOC2005}
{\sc E.~Burman and A.~Ern}, {\em Stabilized {Galerkin} approximation of
  convection-diffusion-reaction equations: discrete maximum principle and
  convergence}, Math. Comp., 74 (2005), pp.~1637--1652.

\bibitem{CaiCheLiuJCP2021}
{\sc Z.~Cai, J.~Chen, and M.~Liu}, {\em Least-squares {ReLU} neural network
  {(LSNN)} method for linear advection-reaction equation}, J.~Comput. Phys.,
  443 (2021), p.~110514.

\bibitem{ChaWicZhuRab2021}
{\sc A.~Chakraborty, T.~Wick, X.~Zhuang, and T.~Rabczuk}, {\em
  Multigoal-oriented dual-weighted-residual error estimation using deep neural
  networks}.
\newblock arXiv:2112.11360v2, 2021.

\bibitem{CiaBOOK2013}
{\sc P.~G. Ciarlet}, {\em Linear and Nonlinear Functional Analysis with
  Applications}, SIAM, Philadelphia, 2013.

\bibitem{DemGopBOOK-CH2017}
{\sc L.~Demkowicz and J.~Gopalakrishnan}, {\em Discontinuous
  {Petrov–Galerkin} ({DPG}) method}, in Encyclopedia of Computational
  Mechanics, Second Edition, E.~Stein, R.~de~Borst, and T.~J.~R. Hughes, eds.,
  Wiley, 2017.
\newblock Part~2~Fundamentals.

\bibitem{DemGopMugZitCMAME2012}
{\sc L.~Demkowicz, J.~Gopalakrishnan, I.~Muga, and J.~Zitelli}, {\em Wavenumber
  explicit analysis of a {DPG} method for the multidimensional {Helmholtz}
  equation}, Comput. Methods Appl. Mech. Engrg., 213-216 (2012), pp.~126--138.

\bibitem{DisHesRayJCP2020}
{\sc N.~Discacciati, J.~S. Hesthaven, and D.~Ray}, {\em Controlling
  oscillations in high-order discontinuous galerkin schemes using artificial
  viscosity tuned by neural networks}, J.~Comput. Phys., 409 (2020), p.~109304.

\bibitem{ECCP2020}
{\sc W.~E}, {\em Machine learning and computational mathematics}, Commun.
  Comput. Phys., 28 (2020), pp.~1639--1670.

\bibitem{EYuCMS2018}
{\sc W.~E and B.~Yu}, {\em The {Deep} {Ritz} {Method}: {A} deep learning-based
  numerical algorithm for solving variational problems}, Commun. Math. Sci., 6
  (2018), pp.~1--12.

\bibitem{ErnGueBOOK2021b}
{\sc A.~Ern and J.-L. Guermond}, {\em Finite Elements II. Galerkin
  Approximation, Elliptic and Mixed PDEs}, vol.~73 of Texts in Applied
  Mathematics, Springer Nature, Switzerland, 2021.

\bibitem{EvaHugSanCMAME2009}
{\sc J.~A. Evans, T.~J. Hughes, and G.~Sangalli}, {\em Enforcement of
  constraints and maximum principles in the variational multiscale method},
  Comput. Methods Appl. Mech. Engrg., 199 (2009), pp.~61--76.

\bibitem{GueSINUM2004}
{\sc J.~L. Guermond}, {\em A finite element technique for solving first-order
  {PDEs} in {$L^p$}}, SIAM J. Numer. Anal., 42 (2004), pp.~714--737.

\bibitem{GuhKutPetAA2020}
{\sc I.~G\"uhring, G.~Kutyniok, and P.~Petersen}, {\em Error bounds for
  approximations with deep {ReLU} neural networks in {$W^{s,p}$} norms},
  Analysis and Applications, 18 (2020), pp.~803--859.

\bibitem{HesUbiJCP2018}
{\sc J.~Hesthaven and S.~Ubbiali}, {\em Non-intrusive reduced order modeling of
  nonlinear problems using neural networks}, J.~Comput. Phys., 363 (2018),
  pp.~55--78.

\bibitem{HigHigSIREV2019}
{\sc C.~F. Higham and D.~J. Higham}, {\em Deep learning: An introduction for
  applied mathematicians}, SIAM Rev., 61 (2019), pp.~860--891.

\bibitem{HinPinUlbUlbBOOK2009}
{\sc M.~Hinze, R.~Pinnau, M.~Ulbrich, and S.~Ulbrich}, {\em Optimization with
  {PDE} Constraints}, Springer, 2009.

\bibitem{JohKnoCMAME2007}
{\sc V.~John and P.~Knobloch}, {\em On spurious oscillations at layers
  diminishing {(SOLD)} methods for convection–diffusion equations: {Part~I}
  – {A~review}}, Comput. Methods Appl. Mech. Engrg., 196 (2007),
  pp.~2197--2215.

\bibitem{KarKevLuPerWanYanNRP2021}
{\sc G.~E. Karniadakis, I.~G. Kevrekidis, L.~Lu, P.~Perdikaris, S.~Wang, and
  L.~Yang}, {\em Physics-informed machine learning}, Nature Reviews Physics, 3
  (2021), pp.~422--440.

\bibitem{KerPruChaLafCMAME2017}
{\sc K.~Kergrene, S.~Prudhomme, L.~Chamoin, and M.~Laforest}, {\em A new
  goal-oriented formulation of the finite element method}, Computer Methods in
  Applied Mechanics and Engineering, 327 (2017), pp.~256--276.
\newblock Advances in Computational Mechanics and Scientific Computation—the
  Cutting Edge.

\bibitem{KhaBalJosSarHegKriGan2021}
{\sc B.~Khara, A.~Balu, A.~Joshi, S.~Sarkar, C.~Hegde, A.~Krishnamurthy, and
  B.~Ganapathysubramanian}, {\em {NeuFENet}: {N}eural finite element solutions
  with theoretical bounds for parametric pdes}.
\newblock arXiv:2110.01601, 2021.

\bibitem{LioBOOK1971}
{\sc J.~L. Lions}, {\em Optimal Control of Systems Governed by Partial
  Differential Equations}, Springer-Verlag, Berlin, 1971.

\bibitem{LiuCaiCheCAMWA2021}
{\sc M.~Liu, Z.~Cai, and J.~Chen}, {\em Adaptive two-layer {ReLU} neural
  network: {I.} {B}est least-squares approximation}, Comput. Math. Appl., 113
  (2022), pp.~34--44.

\bibitem{MinRic2021}
{\sc P.~Minakowski and T.~Richter}, {\em Error estimates for neural network
  solutions of partial differential equations}.
\newblock arXiv:2107.11035v1, 2021.

\bibitem{MisMINE2018}
{\sc S.~Mishra}, {\em A machine learning framework for data driven acceleration
  of computations of differential equations}, Mathematics in Engineering, 1
  (2018), pp.~118--146.

\bibitem{MisMolIMANA2022a}
{\sc S.~Mishra and R.~Molinaro}, {\em Estimates on the generalization error of
  physics-informed neural networks for approximating a class of inverse
  problems for {PDEs}}, IMA J. Numer. Anal., 42 (2022), pp.~981--1022.

\bibitem{MisMolIMANA2022b}
\leavevmode\vrule height 2pt depth -1.6pt width 23pt, {\em Estimates on the
  generalization error of physics-informed neural networks for approximating
  {PDEs}}, IMA J. Numer. Anal.,  (2022), pp.~1--43.
\newblock To appear.

\bibitem{MugZeeSINUM2020}
{\sc I.~Muga and K.~G. van~der Zee}, {\em Discretization of linear problems in
  {Banach} spaces: {Residual} minimization, nonlinear {Petrov--Galerkin}, and
  monotone mixed methods}, SIAM J. Numer. Anal., 58 (2020), pp.~3406--3426.

\bibitem{MulZei2021}
{\sc J.~M\"uller and M.~Zeinhofer}, {\em Error estimates for the deep {Ritz}
  method with boundary penalty}.
\newblock arXiv:2103.01007, 2021.

\bibitem{OdeRedBOOK2011}
{\sc J.~T. Oden and J.~N. Reddy}, {\em An Introduction to the Mathematical
  Theory of Finite Elements}, Dover, Mineola, New York, 2011.
\newblock Unabridged republication of the edition published by John Wiley and
  Sons, New York, 1976.

\bibitem{PenAlbBugCanDeDurGarKarLytPerPetKuhACME2021}
{\sc G.~C.~Y. Peng, M.~Alber, A.~B. Tepole, W.~R. Cannon, S.~De,
  S.~Dura-Bernal, K.~Garikipati, G.~Karniadakis, W.~W. Lytton, P.~Perdikaris,
  L.~Petzold, and E.~Kuhl}, {\em Multiscale modeling meets machine learning:
  {W}hat can we learn?}, Arch. Comput. Methods Eng., 28 (2021), pp.~1017--1037.

\bibitem{PetMOC2017}
{\sc D.~Peterseim}, {\em Eliminating the pollution effect in {Helmholtz}
  problems by local subscale correction}, Math. Comp., 86 (2017),
  pp.~1005--1036.

\bibitem{PetRasVoiFOCM2021}
{\sc P.~Petersen, M.~Raslan, and F.~Voigtlaender}, {\em Topological properties
  of the set of functions generated by neural networks of fixed size}, Found.
  Comput. Math., 21 (2021), pp.~375--–444.

\bibitem{PouJCMDS2022}
{\sc J.~Pousin}, {\em Least squares formulations for some elliptic second order
  problems, feedforward neural network solutions and convergence results},
  Journal of Computational Mathematics and Data Science, 2 (2022), p.~100023.

\bibitem{RaiPerKarJCP2019}
{\sc M.~Raissi, P.~Perdikaris, and G.~E. Karniadakis}, {\em Physics-informed
  neural networks: {A} deep learning framework for solving forward and inverse
  problems involving nonlinear partial differential equations}, J.~Comput.
  Phys., 378 (2019), pp.~686--707.

\bibitem{RanVexSICON2005}
{\sc R.~Rannacher and B.~Vexler}, {\em A priori error estimates for the finite
  element discretization of elliptic parameter identification problems with
  pointwise measurements}, SIAM J. Control Optim., 44 (2005), pp.~1844--1863.

\bibitem{RayHesJCP2018}
{\sc D.~Ray and J.~S. Hesthaven}, {\em An artificial neural network as a
  troubled-cell indicator}, J.~Comput. Phys., 367 (2018), pp.~166--191.

\bibitem{RotSchWicSNAS2022}
{\sc J.~Roth, M.~Schr\"oder, and T.~Wick}, {\em Neural network guided adjoint
  computations in dual weighted residual error estimation}, SN Applied
  Sciences, 4 (2022), p.~62.

\bibitem{SchRayHesJCP2021}
{\sc L.~Schwander, D.~Ray, and J.~S. Hesthaven}, {\em Controlling oscillations
  in spectral methods by local artificial viscosity governed by neural
  networks}, J.~Comput. Phys., 431 (2021), p.~110144.

\bibitem{ShiZhaKar2020}
{\sc Y.~Shin, Z.~Zhang, and G.~E. Karniadakis}, {\em Error estimates of
  residual minimization using neural networks for linear pdes}.
\newblock arXiv:2010.08019, 2020.

\bibitem{SirSpiJCP2018}
{\sc J.~Sirignano and K.~Spiliopoulos}, {\em {DGM:} a deep learning algorithm
  for solving partial differential equations}, J.~Comput. Phys., 375 (2018),
  pp.~1339--1364.

\bibitem{TeiNatVenGarCMAME2019}
{\sc G.~Teichert, A.~Natarajan, A.~{Van der Ven}, and K.~Garikipati}, {\em
  Machine learning materials physics: {I}ntegrable deep neural networks enable
  scale bridging by learning free energy functions}, Comput. Methods Appl.
  Mech. Engrg., 353 (2019), pp.~201--216.

\bibitem{TroBOOK2010}
{\sc F.~Tr{\"{o}}ltzsch}, {\em Optimal Control of Partial Differential
  Equations: Theory, Methods and Applications}, vol.~112 of Graduate Studies in
  Mathematics, American Mathematical Society, Providence, 2010.

\bibitem{UriParOmeCMAME2022}
{\sc C.~Uriarte, D.~Pardo, and Ángel Javier~Omella}, {\em A finite element
  based deep learning solver for parametric {PDEs}}, Comput. Methods Appl.
  Mech. Engrg., 391 (2022), p.~114562.

\bibitem{WanSheLonDonCICP2020}
{\sc Y.~Wang, Z.~Shen, Z.~Long, and B.~Dong}, {\em Learning to discretize:
  {Solving} {1D} scalar conservation laws via deep reinforcement learning},
  Commun. Comput. Phys., 28 (2020), pp.~2158--2179.

\bibitem{XuCCP2020}
{\sc J.~Xu}, {\em Finite neuron method and convergence analysis}, Commun.
  Comput. Phys., 28 (2020), pp.~1707--1745.

\bibitem{XuDarJCP2022}
{\sc K.~Xu and E.~Darve}, {\em Physics constrained learning for data-driven
  inverse modeling from sparse observations}, J.~Comput. Phys., 453 (2022),
  p.~110938.

\bibitem{YarNN2017}
{\sc D.~Yarotsky}, {\em Error bounds for approximations with deep {ReLU}
  networks}, Neural Networks, 94 (2017), pp.~103--114.

\end{thebibliography}
\bibliographystyle{siam}
%=============================================================================%
%
%
%-----------------------------------------------------------------------------%
%
%
%=============================================================================%
\end{document}